\documentclass[12pt]{article}
\usepackage{authblk}    
\usepackage[pagebackref=true]{hyperref}
\hypersetup{colorlinks=true,linkcolor=blue,citecolor=red}
\usepackage[a4paper,margin=1in]{geometry} 
\usepackage{amssymb}
\usepackage{amsmath}

\usepackage{mathtools}  

\usepackage{xcolor}

\usepackage{graphicx} 

\usepackage{gnuplot-lua-tikz} 

\usepackage{hyperref}   
\usepackage{amsthm}     

\usepackage[normalem]{ulem}     

\usepackage{verbatim}     

\usepackage{physics} 	

\usepackage{mathrsfs} 

\usepackage{subcaption} 

\usepackage{booktabs}	

\usepackage{upgreek}

\usepackage{nomencl}    
\makenomenclature
\usepackage[title]{appendix} 

\newtheorem{theorem}{Theorem}

\usepackage{textcomp} 
\usepackage{fancyhdr} 

\fancypagestyle{firstpage}{
    \fancyhf{} 
    \fancyfoot[L]{\small\copyright\ 2026. All rights reserved.}
}

\newcommand{\cc}{\mathrm{c}}
\newcommand{\st}{\mathrm{s}}
\newcommand{\un}{\mathrm{u}}
\newcommand{\hr}{\mathrm{s}}
\newcommand{\pr}{\mathrm{p}}
\newcommand{\vr}{\mathrm{v}}
\newcommand{\iu}{\mathrm{i}\,}
\newcommand{\torus}{T}
\newcommand{\rtorus}{\uptau}
\newcommand{\Tpi}{\mathbb{T}_{2\pi}}
\newcommand{\R}{\mathbb{R}}
\newcommand{\Z}{\mathbb{Z}}
\newcommand{\T}{\mathbb{T}}
\newcommand{\ee}{\,\mathrm{e}}
\newcommand{\TM}{\mathcal{T}} 
\newcommand{\tM}{\mathcal{\tau}} 
\newcommand{\SM}{\mathcal{S}} 
\newcommand{\sM}{\mathcal{\sigma}} 
\newcommand{\IM}{\mathcal{F}} 
\DeclareMathAlphabet{\mathpzc}{OT1}{pzc}{m}{it} 
\newcommand{\iM}{\mathpzc{f}}                    
\newcommand{\Fh}{f_{h}}
\newcommand{\Pext}{P_{\mathrm{ext}}}
\newcommand{\Text}{T_{\mathrm{ext}}}
\newcommand{\Tint}{T_{\mathrm{int}}}
\newcommand{\hTint}{\hat{T}_{\mathrm{int}}}
\newcommand{\Next}{N_{\mathrm{ext}}}
\newcommand{\Nint}{N_{\mathrm{int}}}
\newcommand{\bZ}{\boldsymbol{Z}}
\newcommand{\bX}{\boldsymbol{X}}
\newcommand{\bz}{\boldsymbol{z}}
\newcommand{\bx}{\boldsymbol{x}}
\newcommand{\by}{\boldsymbol{y}}
\newcommand{\bJ}{I}

\newcommand{\bI}{\boldsymbol{I}}
\newcommand{\bphi}{\phi}

\newcommand{\bnu}{\nu}

\newcommand{\nup}{\nu_{\pr}}
\newcommand{\nuv}{\nu_{\vr}}
\newcommand{\fr}{\omega}

\newcommand{\lambdaL}{\lambda^0}
\newcommand{\bnuL}{\nu^0}
\newcommand{\hL}{h^0}
\newcommand{\nupL}{\nup ^0}
\newcommand{\nuvL}{\nu_{\vr}^0}
\newcommand{\nuL}{\nu^0}
\newcommand{\frL}{\omega^0}

\newcommand{\bzeta}{\mathcal{K}}

\newcommand{\bK}{\boldsymbol{K}}
\newcommand{\bk}{\boldsymbol{k}}

\newcommand{\tT}{\boldsymbol{T}}

\newcommand{\NHIM}{\Lambda}         
\newcommand{\SecGlobal}{\mathscr{S}}
\newcommand{\SecLocal}{\Sigma}
\newcommand{\SecMinus}{\Sigma_{-}}
\newcommand{\SecPlus}{\Sigma_{+}}
\newcommand{\SecPlusMinus}{\Sigma_{\pm}}
\newcommand{\NHIMh}{\NHIM_{h}}      
\newcommand{\NHIMhred}{\NHIM^{\SecLocal}_{h}} 
\newcommand{\NHIMc}{\NHIM_{c}}      
\newcommand{\NHIMcred}{\NHIM^{\SecLocal}_{c}}
\newcommand{\LE}{\chi}

\newcommand{\domain}{\mathcal{A}}	

\newcommand{\PauModify}[1]{#1} 
\newcommand{\AModify}[1]{#1}  
\newcommand{\MModify}[1]{#1}  
\newcommand{\error}{\epsilon}
\newcommand{\errorBNF}{\epsilon_M}
\newcommand{\flow}{\Phi}

\newcommand{\redSpace}{U_h^\Sigma} 
\newcommand{\frm}{\mathcal{F}} 

\newcommand{\tin}{t_\mathrm{in}}
\newcommand{\tout}{t_\mathrm{out}}

\newcommand{\nftrans}{\Psi} 

\theoremstyle{plain}
\newtheorem{thm}{Theorem}[section]
\newtheorem{lem}[thm]{Lemma}
\newtheorem{prop}[thm]{Proposition}
\newtheorem{corollary}[thm]{Corollary}

\newtheorem{result}{Numerical Result}

\theoremstyle{definition}
\newtheorem{defn}{Definition}[section]

\theoremstyle{remark}
\newtheorem{rem}{Remark}

\title{Effective shadowing of transition maps for a Hamiltonian saddle-center-center equilibrium}

\author[1]{Rodney L. Anderson} 

\author[2]{Amadeu Delshams} 

\author[3]{Marian Gidea} 

\author[4]{Pablo Roldan\thanks{Corresponding author: Pablo.Roldan@upc.edu}} 

\affil[1]{Jet Propulsion Laboratory, California Institute of Technology}
\affil[2]{Laboratory of Geometry and Dynamical Systems and IMTech, UPC}
\affil[3]{Department of Mathematical Sciences, Yeshiva University}
\affil[4]{Departament de Matem\`atiques, Universitat Polit\`ecnica de Catalunya}

\begin{document}
\maketitle

\thispagestyle{firstpage}

\begin{abstract}
\MModify{We introduce a  general methodology} to construct globally unstable trajectories near the normally hyperbolic invariant manifold (NHIM) of a saddle-center-center Hamiltonian equilibrium. The proposed approach builds upon the 
classical geometric mechanism of searching for trajectories near heteroclinic connections between invariant tori of the NHIM. 
\MModify{The main ingredient is the Birkhoff normal form near the NHIM, which allows us to control the dynamics of the hyperbolic and center components. We construct transition maps on the NHIM that encode the dynamics along heteroclinic connections. By combining the transition maps with the inner dynamics, given by the restriction of the Hamiltonian flow to the NHIM, we show that these connections can be efficiently (in terms of the travel time) shadowed by true trajectories of the full Hamiltonian system.} Finally, the efficiency of this methodology is \MModify{illustrated}  in the case of the NHIM associated to the libration point  $L_{1}$ of the Spatial Restricted 3-Body Problem (SR3BP) for the Sun-Earth system.
\end{abstract}
\tableofcontents

\nomenclature[coordsnf]{$(x,y,I,\phi)=(x,y, I_\pr, I_\vr, \phi_\pr, \phi_\vr) \in U$}{normal form coordinates} 
\nomenclature[coordsorig]{$(X,Y,Z,P_X,P_Y,P_Z)\in\R^6$}{original coordinates}
\nomenclature[nftrans]{$\nftrans^{(M)}\colon U^{(M)} \mapsto \R^6$}{near-identity canonical transformation}
\nomenclature[freqh]{$\lambda(\bI)=\pdv{\mathcal{Z}}{I_\hr} = \lambda^0 + O(\bI)$}{saddle exponent}
\nomenclature[freqp]{$\nup(\bI)=\pdv{\mathcal{Z}}{I_\pr} = \nupL + O(\bI)$}{planar frequency}
\nomenclature[freqv]{$\nuv(\bI)=\pdv{\mathcal{Z}}{I_\vr} = \nuvL + O(\bI)$}{vertical frequency}
\nomenclature[H]{$H(X,Y,Z,P_X,P_Y,P_Z)$}{3-degree-of-freedom Hamiltonian in original coordinates}
\nomenclature[Hnf]{$\mathcal{Z}^{(M)} + \mathcal{R}^{(M)}$}{Hamiltonian in normal form coordinates (order $M$)}
\nomenclature[Hnft]{$\mathcal{Z}^{(M)}(\bI)$}{truncated Hamiltonian in normal form coordinates (order $M$)}
\nomenclature[H2]{$H_2 = \lambda^0 I_\hr + \nupL I_\pr + \nuvL I_\vr=\lambda^0 I_\hr + \nuL\cdot I$}{quadratic part of Hamiltonian in normal form}
\nomenclature[nfinv]{$\bI=(I_\hr, I_\pr, I_\vr)=(xy, I_\pr, I_\vr)$}{normal form invariants}
\nomenclature[ptNHIM]{$\bX=(0,0, \bJ,\bphi)$}{point in $\NHIM$ (4d)}
\nomenclature[ptNHIMred]{$\bx=(0,0, \bJ,\bphi)$ with $I_\pr=f_h(0,I_\vr), \phi_\pr=0$}{point in $\NHIMhred$ (2d)}
\nomenclature[ptmfld]{$\bZ=(x,y,I,\bphi)$, with $xy=0$}{point in $W^{\un,\st}\left(\NHIMhred\right)$ (5d)}
\nomenclature[ptmfldred]{$\bz=(x,y, \bJ,\bphi)$, with $xy=0$, $I_\pr=f_h(0,I_\vr), \phi_\pr=0$ }{point in $W^u(\NHIM)$ or $W^s(\NHIM)$ (3d)}
\nomenclature[ptU]{$\bK=(x,y,I,\phi)$}{general point in $U$ (6d)}
\nomenclature[ptUr]{$\bk=(x,y,I,\phi)$, with $I_\pr=f_h(0,I_\vr), \phi_\pr=0$}{point in $U_h^\Sigma=U\cap\{H=h\}\cap\SecLocal$ (4d)}
\nomenclature[mapFlow]{$\flow_t$}{flow of the Hamiltonian}
\nomenclature[mapiM]{$\iM(\bx)$}{reduced inner map (2d)}
\nomenclature[mapTM]{$\TM(\bX_-)$, $\tM(\bx_-)$}{full, reduced transition map}
\nomenclature[mapSM]{$\SM(\bX_-)$, $\sM(\bx_-)$}{full, reduced scattering map}
\nomenclature[seclocal]{$\SecLocal=\{(x,y,\bJ,\bphi)\in U : \phi_\pr=0\}$}{local Poincar\'e section}
\nomenclature[secglobal]{$\SecGlobal=\{(X,Y,Z,P_X,P_Y,P_Z)\in \R^6 :Y=0\}$}{global Poincar\'e section}
\nomenclature[secglobalpoint]{$\bzeta=(X,Y,Z,P_X,P_Y,P_Z)$}{general point in the original coordinates (6d)}
\nomenclature[secplus]{$\SecPlus=\{(x,y,\bJ,\bphi)\in U : y=\pm\delta\}$}{incoming boundary section}
\nomenclature[secminus]{$\SecMinus=\{(x,y,\bJ,\bphi)\in U : x=\pm\delta\}$}{outgoing boundary section}
\nomenclature[orbguiding]{$\left\{\bx_i\in\NHIMhred\right\}_{i=1,\dotsc,N}$}{\parbox[t]{0.63\textwidth}{guiding orbit, consisting of $N$ points ($N-1$ homoclinic segments)}}
\nomenclature[norm]{$\norm{\bK} = \norm{(x,y,I, \phi)} = \max\left(\abs{x}, \abs{y}, \abs{I_\pr}, \abs{I_\vr}, \abs{\phi_\pr}, \abs{\phi_\vr}\right)$}{maximum norm}

\printnomenclature

\section{Introduction}

In 1964, Arnold~\cite{Arnold64} introduced the method of ``transition chains'' formed by transverse heteroclinic trajectories between invariant tori to prove the existence of global instability in Hamiltonian systems, a phenomenon that later received the name \emph{Arnold diffusion}. 
Arnold’s seminal work was extended through the introduction of geometric methods ~\cite{DelshamsLS00}, based on the existence of a normally hyperbolic invariant manifold (NHIM), which acts as an organizing structure for the dynamics. Assuming that the stable and unstable manifolds of the NHIM intersect transversely along a homoclinic manifold, the dynamics can then be decomposed into an inner dynamics---corresponding to  the restriction of the dynamics to the NHIM---and  outer dynamics---corresponding to the dynamics along homoclinic orbits to the NHIM. The outer dynamics  can be 
encoded in several  ways. 
One way to encode the outer dynamics is via the \emph{scattering map}~\cite{DelshamsLS00}, which relates the asymptotic past of a homoclinic point to its asymptotic future. 
Another way is via the \emph{transition map}~\cite{DelshamsGR13}, which relates the projection onto the NHIM of the starting point of a homoclinic orbit segment to the projection onto the NHIM of the ending point of the segment. 
While the scattering map and the transition map are conjugated by the dynamics, the transition map offers the advantage of direct numerical computability. Computing the scattering map, by contrast, requires knowing both the transition map and the inner dynamics. In this paper, we focus on the transition map as our primary tool.

Combining the inner and the scattering maps, or the transition maps, one can construct diffusion pseudo-orbits. Finally, applying some shadowing lemmas, true diffusion trajectories of the system can be found.
Geometric methods have been used to prove the existence of Arnold diffusion in different perturbative settings, where the scattering map can be computed explicitly to leading order. This is the case, for example, for ``a priori unstable'' Hamiltonian systems, which are systems close to integrable that already possess a NHIM~\cite {DelshamsLS06a,DelshamsLS06b,DelshamsH09,DelshamsH11,DelshamsLS16,DelshamsS17, DelshamsS18,GideaL18,GideaLS20,ChenL22,DelshamsGS24,DelshamsZ25}. This is also the case for  ``a priori-chaotic'' systems, which are perturbations of globally stable systems that possess a NHIM with transverse intersection of their invariant manifolds~\cite{DelshamsLS00,DelshamsLS06b,GelfreichT17, AkingbadeGS23,ClarkeT23,Clarke24}. There are also examples  in celestial mechanics~\cite{DelshamsKRS19,CapinskiG23,ClarkeFG24,AlessiBGG25,ClarkeFG25,BurgosGS25,GuardiaPS25}.

In this  paper,  we consider a   \textbf{non-perturbative} setting that is often encountered in applications, where the system parameters are fixed and specified in advance. If the same geometric structures are present (NHIM, transverse intersection of the stable and unstable manifolds), we may still consider the inner and scattering maps or transition maps, and use them to construct long pseudo-orbits that explore large regions of interest in the phase space. 
However,  in the general case it is not possible to calculate these maps perturbatively; instead,  one can perform a  semi-analytical approximation~\cite{DelshamsGR26}.

Concretely, we consider  a 3-degree-of-freedom Hamiltonian $H$ with  a saddle-center-center equilibrium point $L$, for which the corresponding center manifold has stable and unstable manifolds that  intersect transversally.  We parameterize the center manifold by a system of action-angle variables.  In this context, we study  global instability, specifically, we  explicitly construct  orbits along which one of the action variables varies significantly. 
As we will describe later, a particularly relevant motivation and application is the SR3BP.

We now outline some of the main steps of the construction. 
Around the equilibrium point $L$, we introduce coordinates provided by the \emph{Birkhoff normal form} (BNF for short, see  Section~\ref{sec:BNF} for full details) up to a given order $M$. This construction relies on a general non-resonance condition on the characteristic exponents of $L$ up to order $M$.
In these BNF coordinates, denoted by
\begin{equation}
\label{eq:BNFcoordinates}
(x,y,\bJ,\bphi)=(x,y,I_\pr,I_\vr,\phi_\pr,\phi_\vr)\in U=[-\delta,\delta]^2\times [0,r]^2\times\T^2,
\end{equation}
the Hamiltonian $H$ splits as $H=\mathcal{Z}+\mathcal{R}=\mathcal{Z}^{(M)}+\mathcal{R}^{(M)}$. The integrable part,
\[
\mathcal{Z}=\mathcal{Z}(\bI)=\mathcal{Z}(I_\hr, I_\pr, I_\vr),
\]
depends only on the actions $\bI=(I_\hr, I_\pr, I_\vr)$, where $I_\hr=x y$ represents the saddle component, with $x$ and $y$ being the conjugate saddle variables, and $I=(I_\pr,I_\vr)$ are the two action variables, with conjugate angles $\bphi=\left(\phi_\pr,\phi_\vr\right)$.
The subscripts $\hr$, $\pr$ and $\vr$ stand for the \emph{saddle}, \emph{planar}, and \emph{vertical} components, respectively (the terminology is borrowed from the SR3BP framework). Consequently, the flow of the \emph{truncated normal form} $\mathcal{Z}(\bI)$ can be integrated explicitly (see Section \ref{sec:BNF}).

The term $\mathcal{R}=\mathcal{R}(M,\delta,r)$ is the \emph{remainder}. We estimate its size within the domain $U$ both analytically and numerically (see Sections~\ref{sec:GeneralErrorBNF}, \ref{sec:AnalErrorSEBNF} and \ref{sec:NumErrorSEBNF}), showing that it is typically negligible for practical and numerical purposes.

In BNF coordinates, the 4-dimensional (local) center manifold $\NHIM$ is normally hyperbolic and is given by $x=y=0$. A point
$\bX=(0,0,\bJ,\bphi)\in\NHIM$ is parameterized by the two action variables $\bJ=(I_\pr,I_\vr)$ and the two conjugate angles $\bphi=(\phi_\pr,\phi_\vr)$. Furthermore, its 5-dimensional stable and unstable invariant manifolds are given by equations $x=0$ and  $y=0$, respectively.

In our setting, neglecting the remainder $\mathcal{R}$ implies that the dynamics within $\NHIM$ consist purely of angle translations along invariant tori. This recurrent behavior can be leveraged to reduce the dimensionality of the system through two steps.
First, fixing the energy at $H=h$ isolates one of the action variables. For instance, the planar action can be expressed as $I_\pr=\Fh(I_\hr,I_\vr)$ using the energy equation, which in BNF variables reads $\mathcal{Z}(I_\hr,I_\pr,I_\vr)=h$. Let $U_h=U\cap\{H=h\}$ denote the BNF domain within the energy manifold, and $\NHIMh=\NHIM\cap\{H=h\}$ the NHIM within this same energy manifold.
Second, we eliminate one dimension by choosing an appropriate section within $U_h$, namely $\SecLocal  =\{\phi_\pr=0\}$. 

This allows us to restrict the dynamics to a 4-dimensional \emph{reduced phase space} defined as  $U_h^\Sigma\coloneq  U_h\cap\SecLocal$:
\begin{equation}
\label{eq:ReducedSpace}
\redSpace=\{\bk=\left(x,y,\bJ,\bphi\right): I_\pr=\Fh(I_\hr,I_\vr), \phi_\pr=0, (x,y)\in[-\delta,\delta]^2, I_\vr\in[0,r], \phi_\vr\in\T \},
\end{equation}
which is parameterized by $(x,y,I_\vr,\phi_\vr)$. 
The dynamics on $\redSpace$ is described by the first return map $\frm$ of the Hamiltonian flow, which is a symplectic diffeomorphism with respect to the symplectic form $dy\wedge dx+dI_\vr \wedge d\phi_\vr$. 

Embedded within $\redSpace$ is the two-dimensional \emph{reduced NHIM} 
\begin{equation}
\label{eq:NHIMred}
\NHIMhred\coloneq  \NHIMh\cap\SecLocal  =\left\{\bx=\left(0,0,\bJ,\bphi\right): I_\pr=\Fh(0,I_\vr),\phi_\pr=0, I_\vr\in[0,r], \phi_\vr\in\T \right\},
\end{equation}
which is parameterized solely by $(I_\vr,\phi_\vr)$.
We denote by $\iM$ the restriction of the first return map $\frm$ to $\NHIMhred$.

Then, a semi-analytical computation of the \emph{transition map}~\cite{DelshamsGR13} can be performed. In this paper, we employ the transition map $\TM$ and the reduced transition map $\tM$ instead of the scattering map $\SM$ and the reduced scattering map $\sM$ used in~\cite{DelshamsGR13, DelshamsGR26}. Although both maps are conjugate by the flow (see Section~\ref{sec:scattering_map_transition_map}) and play a similar role, they differ fundamentally in their temporal scales: the transition map encodes finite-time homoclinic  segments, whereas the scattering map encodes infinite-time homoclinic orbits (see Definition~\ref{defn:transition_map_flow}).
Consequently, the transition map implicitly preserves information on the transition time of each homoclinic segment, which is absent in the scattering map. This temporal information is crucial for tracking the total duration of pseudo-orbits. Note that obtaining sharp estimates on diffusion times is notoriously difficult (see, for instance, \cite{Bounemoura11,BertiBB03,Treschev04,DelshamsS18} and references therein).
Full details are given in Section~\ref{sec:BNF}.

Moreover, the transition map can be  directly obtained from the numerical procedure to compute homoclinic trajectories, while the scattering map is indirectly obtained from the transition map (see Section \ref{sec:transition_map_BNF}). We expect that numerical analysts will prefer the transition map, as it remains closer to computations.

Although transition maps are generally not globally defined, they are exact and symplectic on their domains. Consequently, they are generated—at least locally—by a generating function that depends on both old and new variables. In the context of the SR3BP, however, both the transition maps and their corresponding generating functions are globally defined, with the latter being suitably approximated using Fourier-polynomial interpolation (see Section~\ref{sec:semianalytical_results} for more details).



In our context, we obtain several transition maps $\{\tM_1,\dots,\tM_L\}$ associated to different families of homoclinic segments. Then, we consider the Iterated Function System (IFS) $\{\iM, \tM_1,\dots,\tM_L\}$, formed by the transition maps and the inner map. We construct orbits of the IFS along which the vertical action $\bJ_\vr$ undergoes some significant gain. 
Then we develop some \emph{efficient shadowing} strategies (see Section~\ref{sec:main}) that provide true trajectories that achieve the same gain (up to a small error). 
In particular, we consider a simple situation (motivated by applications) in which the IFS consists only of transition maps and not inner maps.
\MModify{Theorem~\ref{th:shadowing} and Theorem~\ref{th:shadowing_topological}} guarantee the existence of true orbits that shadow a given orbit of the IFS \MModify{$\{\tM_1,\dots,\tM_L\}$.}

The shadowing strategies presented in Section~\ref{sec:main} significantly improve existing shadowing theorems in the literature (such as, among others, \cite{FontichM01,DelshamsLS06a,GelfreichT17,GideaLlaveSeara20-CPAM}). While those classical results establish the existence of shadowing orbits, they generally lack a constructive procedure to locate them. Furthermore, they typically require a large number of iterations, which translates into excessively long integration times in the original Hamiltonian system. 

In contrast, by leveraging the BNF to approximate the motion near the NHIM, our shadowing lemmas in Section~\ref{sec:main} are \emph{constructive} and \emph{efficient}. Here, \emph{efficiency} implies that the travel time along the true trajectories closely matches that of the IFS orbits. This approach shares similarities with the shadowing method in~\cite{GelfreichT17}, which utilizes  Fenichel normal form.   

To distinguish the IFS orbits from the true trajectories of the system, we will refer to the former as \emph{guiding orbits}. 
In our setting, a guiding orbit of an iterated function system (IFS) with dynamics $\{\iM, \tM_1,\dots,\tM_L\}$ is simply a finite sequence $(\bx_i)_{i=1}^N$ 
of points in the two-dimensional reduced NHIM $\NHIMhred$, where each iterate satisfies one of the following conditions:
\begin{itemize}
\item \textbf{Inner map:} $\bx_{i+1}=\iM\left(\bx_i\right)$. This map models a ``complete revolution'' of the vertical angle $\phi_\vr$ from $0$ to $2\pi$ within $\NHIMhred$, during which the vertical action $\bJ_\vr$ remains constant.        
\item \textbf{Transition map:} $\bx_{i+1}=\tM_j\left(\bx_i\right)$ for some $j=1,\dots,L$. This map governs the evolution along a segment of a homoclinic orbit to $\NHIMhred$, along which $\bJ_\vr$ may vary.
\end{itemize}


\PauModify{
In Section~\ref{sec:multiple_shooting} we introduce a multiple shooting algorithm to find ``numerical shadowing trajectories'' that shadow a given guiding orbit of the IFS. This section is numerical, but that is usually enough for applications. The shadowing trajectory obtained is, in fact, the concatenation of several (high-precision) numerical integrations that match at the endpoints (up to a small numerical error).
}

We apply this methodology to the Spatial Restricted 3-Body Problem (SR3BP for short). This model describes the spatial dynamics of an infinitesimal body under the gravitational influence of two massive bodies that move on circular orbits about their center of mass. In this paper, the two massive bodies are the Sun and Earth, and the infinitesimal body represents a spacecraft. 
The motion of the spacecraft is described by a 3-degree-of-freedom Hamiltonian system. 
For the equilibrium point $L$, we consider the libration point $L_1$ located between the Sun and the Earth. 
In this case, the action variables $I_\pr$ and $I_\vr$ correspond to the planar and vertical amplitude---relative to the ecliptic---of a quasi-periodic orbit around $L_1$, respectively. We consider homoclinic orbits that go once around the Earth, and select two transition maps. 
 
Following the methodology described for the general case, we construct true trajectories that exhibit a large variation in the vertical amplitude $I_\vr$. These vertical \emph{drift trajectories} closely track a selected finite sequence of homoclinic segments associated with the transition maps. 
These are analogous to diffusion orbits in a perturbative setting, where the action variable changes by an amount of order one \MModify{with respect to the parameter, i.e.,}, independent of the perturbation size. 

High-inclination orbits at the libration points are desirable for a variety of
applications including studying the zodiacal dust cloud~\cite{Reach:1997,
May:2007} and they may also offer observations of the Sun's poles, not studied
since the Ulysses mission~\cite{Wenzel:1995}. Trajectories to these types of
orbits have previously been studied using optimized low-thrust
trajectories~\cite{Stuart:2019}.

Another compelling possibility is the construction of long trajectories that visit designated regions of phase space in a prescribed order. This strategy has been widely proposed for space mission design~\cite{GomezLMS01,KoonLMR08}, where homoclinic and heteroclinic trajectories are combined to enable low-propulsion transfers between different invariant objects in the solar system~\cite{LoWBHMCR01,CanaliasM06}.
To our knowledge, this paper represents the first time that true trajectories shadowing transition chains of two or more homoclinic or heteroclinic orbits have been explicitly computed.

\PauModify{
To illustrate, in Section~\ref{sec:semianalytical_results} we obtain guiding orbits with the desired vertical action gain with a simple simulation of these two or more dynamics (inner and transition maps), leveraging the classic Dijkstra algorithm to find the shortest path in a graph.
%

Then we apply our multiple shooting algorithm to shadow a guiding orbit within the Sun-Earth  SR3BP consisting of 16 iterations of the transition map (corresponding to 16 homoclinic trajectories). The algorithm converges quickly thanks to a careful selection of the initial guess and produces a shadowing trajectory with a numerical error of $10^{-12}$.
}

In Section~\ref{sec:two_level_diff_corr} we use an alternative numerical method (two-level differential correction) to independently validate the shadowing trajectory obtained with multiple shooting.
The differential corrector used here has been previously implemented in
LTool~\cite{Koon:2000b, Lo:2002b}, which has been subsequently used for
a range of missions including the Genesis mission~\cite{Lo:2001c,Wilson:2003}, the MoonRise Mission Concept~\cite{Parker:2015}, and some Europa Clipper studies~\cite{Anderson:2018}. 

The multiple shooting algorithm is able to reduce the total position and velocity error of the pseudo-trajectory from $\Delta x=184\,005$ km and $\Delta v=338\,000$ km/h, to virtually zero error in the shadowing trajectory ($\Delta x<1.5$ m and $\Delta v<0.00003$ m/s). This beats a state-of-the-art method like differential correction by 1-2 orders of magnitude. 
However, this is not surprising, given that the method of Section~\ref{sec:multiple_shooting} has been developed for the special purpose of shadowing chains of homoclinic trajectories, while the method of Section~\ref{sec:two_level_diff_corr} is more general purpose.
An additional advantage is that the multiple shooting method runs automatically, while the differential corrector needs frequent intervention to heuristically select the number and location of patch points.

\section{Preliminaries}\label{sec:scattering_map_transition_map}
\subsection{Scattering map}
We recall the definition of the scattering map following \cite{DLS08}.
We consider a general setting of  flow $\flow_t$ on a manifold $M$, and assume that there is a normally hyperbolic invariant manifold $\Lambda$ for the flow.
Associated to $\Lambda$ we have stable and unstable manifolds,  $W^s(\Lambda)$ and  $W^u(\Lambda)$,  which are foliated  by stable and unstable leaves
$W^s(\Lambda)=\bigcup_{\by\in\Lambda}W^s_{\by}$ and $W^u(\Lambda)=\bigcup_{\by\in\Lambda}W^u_{\by}$.
The canonical projections along the leaves are  defined as follows: $\Omega_+:W^s(\Lambda)\to \Lambda$ is given by
$\Omega_+(\bx)=\bx^\st$ provided that $\bx\in W^s(\bx^\st)$ and $\Omega_-:W^u(\Lambda)\to \Lambda$ is given by
$\Omega_-(\bx)=\bx^\un$ provided that $\bx\in W^u(\bx^\un)$. 

We will implicitly assume that the the manifold $M$, the flow $\flow_t$, as well as the NHIM $\Lambda$ and the related geometric structures described above are differentiable enough.

Suppose  that  $W^s(\Lambda)$ and  $W^u(\Lambda)$
intersect transversally along a
homoclinic manifold $\Gamma$, that is,
for each $\bx\in\Gamma  \subseteq W^u(\Lambda) \cap W^s(\Lambda)$, we have
\begin{equation*} \begin{split}
		T_{\bx}M=T_{\bx}W^u(\Lambda)+T_{\bx}W^s(\Lambda),\\
		T_{\bx}\Gamma=T_{\bx}W^u(\Lambda)\cap T_{\bx}W^s(\Lambda).
\end{split} \end{equation*} 
and moreover, the following strong transversality condition holds:
\begin{equation*} 
	\begin{split}
		T_{\bx}W^s(\Lambda)=T_{\bx}W^s(\bx_+)\oplus T_{\bx}(\Gamma),\\
		T_{\bx}W^u(\Lambda)=T_{\bx}W^u(\bx_-)\oplus T_{\bx}(\Gamma).
	\end{split}
\end{equation*}

Then, the restrictions $\Omega_+^\Gamma$, $\Omega_-^\Gamma$ of
$\Omega_+$, $\Omega_-$, respectively, to $\Gamma$ are  local
diffeomorphisms.
If $\Gamma$ is so  that
$\Omega_+^\Gamma$, $\Omega_-^\Gamma$ are diffeomorphisms onto their images, then it is referred to as a \emph{homoclinic channel}. 

When we flow the homoclinic channel $\Gamma$ to $\flow_t(\Gamma)$, $t\in\mathbb{R}$,
the corresponding wave maps are conjugate by the flow, i.e., 
\begin{equation}\label{eqn:conjugate_wave_maps}
\begin{aligned}
	\Omega_+^{\flow_t(\Gamma)} &= \flow_{t} \circ\Omega_+^\Gamma\circ  \flow_{-t} .
\end{aligned}
\end{equation}


\begin{defn} \label{defn:scattering_map_flow}
	Given a homoclinic channel $\Gamma$, the scattering map associated to
	$\Gamma$ is a  diffeomorphism from
	$ \textrm{Dom}(\SM^\Gamma):=\Omega^\Gamma_-(\Gamma)\subseteq\Lambda$ to
	$ \textrm{Im}(\SM^\Gamma):= \Omega^\Gamma_+(\Gamma)\subseteq\Lambda$ given by
	\[\SM^\Gamma=\Omega^\Gamma_+\circ (\Omega^\Gamma_-)^{-1}.\] 
\end{defn}

The scattering map depends on the choice of the homoclinic channel $\Gamma$. In general, its domain and range are proper subsets of $\Lambda$.

From Definition \ref{defn:scattering_map_flow} it follows that,
if $\bx\in \Gamma$ is a homoclinic point and
\[\SM^\Gamma(\bx^\un)=\bx^\st,\textrm{ where } \bx^{\un,\st}\in\Lambda,  \]
then
\begin{equation*}
	\begin{split}
		d(\flow_{-t_-}(\bx),\flow_{-t_-}(\bx^\un))& \to 0, \textrm { as } t_-\to\infty,\\
		d(\flow_{t_+}(\bx),\flow_{t_+}(\bx^\st))& \to 0, \textrm { as } t_+\to\infty.
	\end{split}
\end{equation*}

A remarkable property of the  scattering map is that it is (exact) symplectic, provided that the manifold $M$ and the flow $\flow_t$ are (exact) symplectic.
We refer to \cite{DLS08} for details.

When we flow the homoclinic channel $\Gamma$ to $\flow_t(\Gamma)$, it follows from \eqref{eqn:conjugate_wave_maps} that the corresponding scattering maps  are conjugate by the flow, so we have the following invariance property
\begin{equation}\label{eqn:conjugate}
	\SM^{\flow_t(\Gamma)}= \flow_{t} \circ\SM^\Gamma\circ  \flow_{-t},
\end{equation}
for $t\in\mathbb{R}$.
We note that on the right-hand side the above, the flow $\flow_t$ acts on $\Lambda$.

\subsection{Transition map}

We now recall the \emph{transition map} (see \cite{DLS08}).
\begin{defn} \label{defn:transition_map_flow}
	Given a homoclinic channel $\Gamma$ and  (sufficiently large) times $T_-,T_+$, the transition map  $\TM=\TM^\Gamma_{T_-,T_+}: \flow_{-T_-}(\textrm{Dom}(\SM^\Gamma))\to  \flow_{T_+}(\textrm{Im}(\SM^\Gamma))$ is
	given by
	\begin{equation} \label{eq:SMAndTM}
		\TM=\flow_{T_+}\circ \SM^\Gamma\circ\flow_{T_-}.
	\end{equation}
We write $\mathcal{T}(\bx_-)=\bx_+$. 
\end{defn}

\begin{rem}\label{rem:transition_alternative}
In Definition \ref{defn:transition_map_flow} the times $T_+,T_-$ are fixed.
More generally, we can consider $T_+,T_-$ as continuous functions on $\Gamma$. In Section~\ref{sec:transition_map_BNF}, given  a neighborhood  $U$ of the NHIM $\Lambda$, and $\Sigma_+$, $\Sigma_-$ two surfaces of section inside $U$, for any $\bx\in\Gamma$ we will
define $T_-=T_-(\bx)$ as the first time $t>0$ when  
$\flow_{-t}(\bx)\in \Sigma_-$, and $T_+=T_+(\bx)$ as the first time $t>0$ when  
$\flow_{t}(\bx)\in \Sigma_+$.
\end{rem}

Since both the scattering map  and the inner dynamics are (exact) symplectic, it follows that the transition map is also (exact)   symplectic.  

Using the conjugacy~\eqref{eqn:conjugate_wave_maps}, the transition map can also be written
\begin{equation} \label{eq:transition_map}
	\TM=\Omega_+^{\Gamma_+} \circ \flow_{T_+}\circ \flow_{T_-} \circ \left(\Omega_-^{\Gamma_-}\right)^{-1},
\end{equation}
where we have introduced $\Gamma_+ := \flow_{T_+}(\Gamma)$ and $\Gamma_- := \flow_{-T_-}(\Gamma)$. 
This will be used in Section~\ref{sec:transition_map_BNF} to compute the transition map, because the local wave maps $\Omega_\pm^{\Gamma_\pm}$ have particularly simple expressions in BNF coordinates.

Now we  consider  the case of a  discrete-time dynamical system, given by a diffeomorphism $\frm$ on $M$. Let $\Lambda$ be a NHIM for $\frm$, and let $f$ denote the restriction of $\frm$ to $\Lambda$. Let $\Gamma$ be a homoclinic channel. 

Then the scattering map $\sigma^\Gamma$  can be defined in a similar way to the flow case. 

The invariance property  \eqref{eqn:conjugate} becomes:
\[\sigma ^{f(\Gamma)} =f\circ \sigma^\Gamma \circ f^{-1}.\]

The transition map $\tM=\tM^\Gamma_{m,n}$ can be defined in a similar way to the flow case, i.e., 
\[\tM=f^m\circ\sigma^\Gamma\circ f^n,\]
for $m,n> 0$ sufficiently large.

Using the invariance of the scattering map we obtain
\[\tau^\Gamma_{m,n}=f^m\circ\sigma^\Gamma\circ f^n=f^{m-1}\circ \sigma^{\frm(\Gamma)}\circ f^{n+1}=\tau^{\frm(\Gamma)}_{m-1,n+1}.\]
It follows:
\[\tau^\Gamma_{m,n}=\tau^{\frm^{k}(\Gamma)}_{m-k,n+k},\]
for any $k\in\mathbb{Z}$. 


\subsection{Shadowing lemmas}
Let $\Lambda$ be a NHIM for $\frm$, $f=\frm_{\mid\Lambda}$ and 
$\sigma_{\alpha}$, $\alpha\in \{1,\ldots,L\}$, be a finite collection of scattering maps. (We omit the homoclinic channel from the notation.)

We assume that there is a neighborhood $\mathcal{U}$ of $\Lambda$  on which there is as system of coordinates $(x,y,\theta)$ with the following properties:
\begin{equation}\label{eqn:coordinates}
\begin{split}
\Lambda=&\{(x,y,\theta)\,:\, x=y=0\},\\
W^u(\Lambda)=&\{(x,y,\theta)\,:\, y=0\},\\
W^s(\Lambda)=&\{(x,y,\theta)\,:\, x=0\},\\
W^u(\bx)=&\{(x,y,\theta)\,:\, y=0,\, \theta=\theta(\bx)\} \textrm{ for } \bx\in\Lambda,\\
W^s(\bx)=&\{(x,y,\theta)\,:\, x=0,\, \theta=\theta(\bx)\} \textrm{ for } \bx\in\Lambda.
\end{split}
\end{equation}
\begin{rem}
In the context of  Sections \ref{sec:BNF} and \ref{sec:transition_map_BNF}, the coordinates $(x,y,\theta)$ represent the normal form coordinates $(x,y,I,\phi)$---with $\theta=(I,\phi)$---as described in \eqref{eq:ReducedSpace}. \end{rem}

We consider a sequence of  transition maps $\tau_i$, $i=1,...,N-1$, where 
\[\tau_{i}=f^{m_i}\circ\sigma_{\alpha(i)}\circ f^{n_i},\]
for some scattering map $\sigma_{\alpha(i)}$ with $\alpha(i)\in \{1,\ldots,L\}$ and  some $m_i,n_i>0$.  Each transition map $\tau_i$ depends on both the scattering map $\sigma_{\alpha_i}$ 
and the orders of iterates $m_i,n_i$. Hence, one can have $\tau_i=\tau_j$ for some $i\neq j$. 
\MModify{We denote by  $M_i:=m_i+n_i$ the number of iterates associated to the transition map  $\tau_i$.}

Let \[\bx_1,\bx_2\ldots,\bx_N\] be a guiding orbit in $\Lambda$ of the transition maps $\tau_1,\ldots, \tau_{N-1}$, in the sense that
\[\tau_i(\bx_{i})=\bx_{i+1}\, \textrm { for } i=1,\ldots,N-1.\]

\begin{rem}
\MModify{Above, we implicitly assume that
the range of $\tau_i$ intersects the domain of $\tau_{i+1}$, for $i=0,\ldots,N-1$. In the context of the SR3BP, this is automatically satisfied since the scattering map and the corresponding transition maps are globally defined. See Section~\ref{subsec:Series_transition_maps}.
}
\end{rem}


\begin{lem}\label{lem:shadowing_topological}

Let $\sigma_\alpha$ be a family of scattering maps with associated homoclinic channels $\Gamma_\alpha$, for  $\alpha \in\{1,\ldots, L\}$.

For every $\varepsilon>0$ there exists $M_{\min}<M_{\max}$ such that, if 
\[\tau_i= f^{m_i}\circ\sigma_{\alpha(i)}\circ f^{n_i}, i=1,\ldots,N-1,\] 
is a transition map
associated with the homoclinic channel $\Gamma_{\alpha(i)}$, and with $m_i,n_i$ satisfying
 $M_{\min} \le M_i=m_i+n_i \le M_{\max}$,  and if \[\bx_{i+1}=\tau_i(\bx_i),  \textrm { for }i=1,\ldots,N-1,\]
then there exists a point $\by_1$ whose orbit $\by_1,\ldots, \by_N$  $\varepsilon$-shadows $\bx_1,\bx_2,\ldots, \bx_N$ in the sense
\begin{equation*}\begin{split}
\by_{i+1}=&\frm^{M_i}(\by_i),\textrm { and }\\
d(\by_i,\bx_i)<&\varepsilon\, \textrm{ for } i=1,\ldots,N.
\end{split}
\end{equation*}
The orbit $\by_1,\ldots, \by_N$ is a true orbit of the system. The number of iterates to take $\by_i$ to $\by_{i+1}$ is $M_i=m_i+n_i$, which is the same as  the number of iterates associated to the transition map  $\tau_i$.

In particular, we can choose the number of iterates $m_i=m_*$, $n_i=n_*$, and $M_i=M_*$ to be independent of $i$.
\end{lem}

The proof of this lemma will be given in Section \ref{sec:shadowing_topological}.

\begin{rem}
We compare Lemma \ref{lem:shadowing_topological} with some similar shadowing lemma type of results. The paper \cite{GideaLlaveSeara20-CPAM} provides a result on shadowing  orbits generated  by scattering maps on a NHIM. The argument requires that the dynamics restricted to the NHIM satisfies Poincar\'e recurrence and, therefore, the shadowing orbits may take a long time. Another result in the paper is on shadowing  orbits obtained by alternately applying scattering maps and iterates of the inner dynamics. While this result does not require Poincar\'e recurrence, it requires that the number of iterates at each step is sufficiently large with respect to the number of iterates at all previous steps. A similar result to the latter appears in \cite{GelfreichT17}.  
Our Lemma~\ref{lem:shadowing_topological}  concerns shadowing orbits of transition maps and gives a uniform number of iterates of the inner dynamics. Since the transition map also encodes the effect of the inner dynamics, and we require that the input of each transition map coincides with the output of the preceding transition map, there is no need to invoke recurrence.
However, this lemma only allows for the shadowing of finite orbits.
\end{rem}

\section{Main results}
\label{sec:main} 

\MModify{Some of the standard shadowing results about finding true trajectories that follow given pseudo-orbits are not practical for concrete applications. On the one hand, the existence of the true trajectories is not constructive and, on the other hand, the theoretical shadowing results may involve a large number of iterations, which  translates into an enormous time in the original Hamiltonian.}

The main objective of this paper is to provide a \emph{constructive} and \emph{efficient} shadowing of the guiding orbits.

\subsection{Main shadowing results}

The first shadowing result is achieved by using Birkhoff Normal Form (BNF) variables within the domain $U$~\eqref{eq:U}, where the reduced NHIM $\NHIMhred$ takes the form~\eqref{eq:NHIMred} and the associated reduced transition maps are defined.  Outside $U$, the candidates for drift trajectories are \MModify{$N$ orbit segments }   that exit $U$  through the unstable (outgoing) section $\SecMinus=\{x=\pm\delta\}$ and return to $U$ through the stable (incoming) section $\SecPlus=\{y=\pm\delta\}$, traveling close to the selected homoclinic trajectories to the NHIM $\NHIM$. 
Note that the drift trajectories constructed in Theorem~\ref{th:shadowing} (below) are not arbitrarily close to the guiding orbit---as is the case in standard shadowing lemmas---but are instead only $\delta$-close. This restriction minimizes the drift time, matching the duration spent along the homoclinic segments that generate the guiding orbit.

Moreover, to obtain drift trajectories with the shortest travel time, we consider in Theorem~\ref{th:shadowing} guiding orbits formed exclusively by iterations of transition maps. Consequently, we seek trajectories that cross the section $\Sigma$ only once inside $U$, which implies they will not remain too close to the NHIM. Thus, to guarantee a proper return, we impose the \emph{homecoming condition}~\eqref{eq:HomecomingCondition}:
\begin{equation*}
(N-1)\LE < 2\pi\lambda/\nup,
\end{equation*}
where the saddle exponent $\displaystyle\lambda=\frac{\partial \mathcal{Z}}{\partial I_\hr}$ and the planar frequency
$\displaystyle\nu_\pr=\frac{\partial \mathcal{Z}}{\partial I_\pr}$ are introduced in equations~(\ref{eq:flow_NF}--\ref{eq:nu}), whereas $\ee^{\LE}$ is a uniform bound for the maximum of the derivatives of the exterior Poincar\'e maps $\Pext^j:\SecMinus\to\SecPlus$ associated with the transition maps $\tM_j$ for $1\leq j\leq L$.

\begin{theorem}[\MModify{Analytical shadowing of a guiding orbit of transition maps}]
\label{th:shadowing}
Let $\left(\bx_i\right)_{i=1}^N$ be a guiding orbit of the iterated function system (IFS) defined on the reduced NHIM $\NHIMhred$~\eqref{eq:NHIMred} generated \textbf{solely} by $L\geq 1$ reduced transition maps $\{\tM_1,\ldots,\tM_L\}$, expressed in BNF coordinates~\eqref{eq:BNFcoordinates}. 
Let $\delta$ be the size of the saddle variables~\eqref{eq:ReducedSpace} within the BNF domain $U$,
and assume that the  \emph{homecoming condition}~\eqref{eq:HomecomingCondition} holds.

Then, there exists a true trajectory $\flow_t({\bf x})$ of the system that $\delta$-shadows the guiding orbit. Specifically, there exists a strictly increasing sequence of times $0 < t_1 < \dots < t_N$ such that
\[
\norm{\flow_{t_i}({\bf x}) -   \bx_i } \leq \delta \quad \text{for all } i = 1, \dots, N.
\]
Moreover, if the guiding orbit $\left(\bx_i\right)_{i=1}^N$ is generated by a sequence of homoclinic segments 
$\left[\bz_i^-,\bz_{i+1}^+\right]$, where $\bz_i^-\in W^\un\left( \bx_i\right)\cap\NHIMhred$ and 
$\bz_{i+1}^+=\flow_{\tT^i}\left(\bz_i^-\right)\in W^\st\left( \bx_{i+1}\right)\cap\NHIMhred$, then each $t_i$ can be chosen as $t_i=\tT^1+\cdots+\tT^i$, yielding a total drift time of $t_N=\tT^1+\cdots+\tT^N$.
\end{theorem}

In Section~\ref{sec:analytical_shadowing}, we present a direct analytical proof for the construction of this shadowing trajectory. Our approach relies fundamentally on the transition map and the BNF; these tools enables us to obtain an optimal time estimate that matches the duration given by the sequence of homoclinic segments associated to the guiding orbit. 

General guiding orbits are composed of transition maps and inner maps. Because each inner map adds a revolution around the NHIM, these guiding orbits do not generally provide optimal drift time estimates.

Subsequently, in Section~\ref{sec:multiple_shooting} we introduce a complementary numerical method to find trajectories that shadow the guiding orbit. It is based on a multiple shooting scheme and an appropriate choice of initial guess. In fact, the method does not produce an exact shadowing trajectory, only a numerical approximation, but this is often enough for many practical applications. Moreover, this numerical approximation could theoretically be used as the seed for a computer assisted proof of a true shadowing orbit near it (see for example~\cite{KPC07}).

Theorem~\ref{th:shadowing} provides an efficient shadowing of a given guiding orbit. In particular, the shadowing orbits spend  little time near the NHIM (they cross the section $\Sigma$ only once inside $U$). However, the error in the shadowing equals the size $\delta$ of the BNF neighborhood, which can be large. 

Theorem~\ref{th:shadowing_topological}  below  provides an alternative shadowing result that minimizes the shadowing error. If we ignore the ``effects of the truncation error" $\errorBNF=\errorBNF(\delta)$ given by~\eqref{eq:errorBNF},
for any   ``mathematical error" $\varepsilon$ it gives a threshold value on the number of iterates of the inner map that are required by the transition maps. 
Then, any orbit generated by such transition maps can be shadowed 
within a $\varepsilon$-mathematical error (to which we have to add the effects of the truncation.) However, if  $\varepsilon$ is small, the shadowing orbits may spend a long time near the NHIM 
(they may cross the section $\Sigma$ multiple times  inside $U$).

\begin{theorem}[Topological shadowing of a guiding orbit of transition maps]
\label{th:shadowing_topological}


Let $\errorBNF=\errorBNF(\delta)$ denote the effect of truncation error~\eqref{eq:errorBNF} of the transformation to BNF variables.

Then, for every $\epsilon>0$ there exist $m_*,n_*$ sufficiently large and a constant $C>0$, such that, for any guiding orbit $\left( \bx_i\right)_{i=1}^N$ 
given by
\[  \bx_{i+1}=\tau_i( \bx_i), \textrm{ for } i=1,\ldots,N-1, \]
where each transition map $\tau_i$ is of the form
\[\tau_i=f^{m_*}\circ\sigma_{\alpha(i)} \circ f^{n_*}, \]
there exists a true trajectory $\flow_t({\bx})$ that $(\varepsilon+C\errorBNF)$-shadows the guiding orbit, 
\[
\norm{\flow_{t_i}({\bx}) -   \bx_i } \leq \epsilon + C \errorBNF \quad \text{for all } i = 1, \dots, N, 
\]
for some strictly increasing sequence of times $0 < t_1 < \dots < t_N$. 
\end{theorem}

The proof is given in Section~\ref{sec:shadowing_topological}.

\subsection{Main numerical result}

In Sections~\ref{sec:NumErrorSEBNF} and~\ref{sec:semianalytical_results}, we study the Sun-Earth spatial R3BP around the Lagrange equilibrium point $L_1$. Using semi-analytical methods, we find the shortest-time drift guiding orbit in this particular setting, and then apply the multiple shooting scheme of Section~\ref{sec:multiple_shooting} to shadow it.

\begin{result}\label{thm:main_numerical}
Consider the SR3BP with a mass ratio given by $\mu_{\text{Sun-Earth}} = 3.040423398444176 \cdot 10^{-6}$, described in Section~\ref{sec:AnalErrorSEBNF}. Let the energy value $c=3.00088$ be fixed.
\begin{itemize}
\item[(a)] Following~\cite{DelshamsGR16}, we obtain a semi-analytical approximation of the NHIM $\NHIMc$ around the equilibrium point $L_1$ using the BNF of order $M=16$. 
\item[(b)] Consider a neighborhood $U=U(\delta)$ of size $\delta=0.1$ in the saddle BNF variables. In Section \ref{sec:NumErrorSEBNF} we show that the error of the semi-analytical approximation of the NHIM in $U(\delta)$ is less than $3\times 10^{-7}$. 
\item[(c)]
By part (b), the level sets of $I_\vr$ are preserved by the inner map $\iM$ up to an error of $\varepsilon_1=3\times 10^{-7}$. 
In Section~\ref{subsec:Series_transition_maps} we  derive series representations for two different reduced transition maps $\tM_1, \tM_2$ that are globally defined on an annulus $\domain\subset \NHIMcred$, corresponding to the range $I_\vr\in[I_1,I_2]$, where $I_1=0.001$ and $I_2=0.007$.
\item[(d)]
In Section~\ref{subsec:Shortest-time_guiding} we find the orbit of the IFS $\{\iM,\tM_1,\tM_2\}$ that takes the least amount of time to drift from the lower boundary to the upper boundary in $\domain$. See Figure~\ref{fig:diffusion-shortest-path-delta-0.1}. It turns out that, in this particular scenario, the shortest-time drift orbit does not use the inner map $\iM$; it consists solely of iterates of the transition maps $\tM_1, \tM_2$.
\item[(e)] Let $\left(\bx_i\right)_{i=1}^{17}$ be the guiding orbit of the IFS $\{\tM_1, \tM_2\}$ consisting of the last $N=17$ points in the shortest-time drift orbit found in part (d).
In Section~\ref{sec:Shadowing_trajectory}, we apply the multiple shooting method of Section~\ref{sec:multiple_shooting} to construct a (numerical) shadowing trajectory of the SR3BP close to the guiding orbit. The result is given in Table~\ref{tab:OneCut} and shown in Figure~\ref{fig:shadowing_trajectory}.
\end{itemize}
\end{result}

It should be emphasized here that the trajectories found are \MModify{true trajectories of the flow of SR3BP}, that is, ``zero cost'' trajectories. If a slight increase in speed Delta v ($\Delta v$) (or a slight thrust) is allowed in the spacecraft, faster drift trajectories could be achieved, that is,  trajectories that take less time. These would thus be ``low cost'' drift trajectories.

\section{Birkhoff normal form (BNF)}\label{sec:BNF}

Around a saddle-center-center equilibrium point of a Hamiltonian system,
the parameterization method provides an efficient framework to compute the local center, center-stable and center-unstable manifolds, see, e.g., \cite{CabreFL05,HaroCFLM16}. Alternatively, the classical Birkhoff normal form (BNF),
(see, e.g., \cite{GiorgilliDFGS89,DelshamsG96,GuzzoCP26}), or its variants, such as resonant or saddle normal forms \cite{JorbaM99}, can be employed.  Notably, these normal form approaches offer the additional advantage of approximating the full flow in a neighborhood of the equilibrium point. 

Let $V$ be a neighborhood of a saddle-center-center point $L$ of a three degrees of freedom Hamiltonian $H(X,Y,Z,P_X,P_Y,P_Z)$ where the BNF up to some order $M$ is computed, and
let 
\begin{equation}\label{eq:nf}
\bK=(x,y,\bJ,\bphi)\in U\mapsto (X,Y,Z,P_X,P_Y,P_Z)=\bzeta=\nftrans(\bK)\in V=\nftrans(U),
\end{equation}
where $I=(I_\pr,I_\vr)$ and $\bphi=(\phi_\pr,\phi_\vr)$, be a symplectic transformation with respect to the symplectic forms $dx\wedge dy+d\phi_\pr\wedge dI_\pr+d\phi_\vr\wedge dI_\vr$ on $U$ and 
$dX\wedge dP_X+dY\wedge dP_Y+ dZ\wedge dP_Z$ on $V$, 
where 
\begin{equation}
\label{eq:U}
U=U(\delta,r)=\left\{\bK=(x,y,\bJ,\bphi)\in[-\delta,\delta]^2\times [0,r]^2\times\T^2\right\}.
\end{equation} 
The transformation  $\nftrans=\nftrans^{(M)}$ is such that  $H\circ\nftrans=\mathcal{Z}(\bI)+\mathcal{R}=\mathcal{Z}^{(M)}(\bI)+\mathcal{R}^{(M)}$,
where $\bI\coloneq(I_\hr,\bJ)$, $I_\hr\coloneq x y$, $\bJ\coloneq \left(I_\pr,I_\vr\right)$ are three first integrals in involution of the normal form $\mathcal{Z}(I_\hr, \bJ)$,   $\bphi\coloneq\left(\phi_\pr,\phi_\vr\right)$ are the conjugate angles to the two actions $\bJ$,
and $\mathcal{R}$ is the remainder. 

The normalized Hamiltonian expands as 
\[\mathcal{Z}(\bI)=\hL+\lambdaL I_\hr+\bnuL\cdot \bJ+O_2(\bI),
\textrm { with } \bnuL \cdot \bJ=\nupL I_\pr+\nuvL I_\vr,\]
where $\lambdaL$ is the saddle characteristic exponent at the equilibrium $L$, while $\nupL$, $\nuvL$ are the planar and vertical frequencies at $L$. 
For definiteness, we will assume that $\hL, \lambdaL, \nupL, \nuvL>0$.
Note that $\nftrans$ is a genuine coordinate change (i.e., a diffeomorphism) for $I_\pr>0$ and $I_\vr>0$. In particular, 
as is usual when considering polar coordinates, the BNF coordinates are degenerate at the equilibrium point $L$, which is given by $x=y=\bJ_\pr=I_\vr=0$ for arbitrary angles $\phi_\pr$ and $\phi_\vr$. Analogously for the planar and vertical Lyapunov orbits given by $x=y=I_\vr=0$ and $x=y=I_\pr=0$, respectively. Note also that, for simplicity, we are considering the same domain for the saddle variables $x$, $y$, and respectively for the actions $I_\pr$ and $I_\vr$.

The main feature of the BNF is that often, by choosing the order $M$ appropriately in terms of the parameters 
$\lambdaL, \nupL, \nuvL$, and the sizes $\delta,r$ of the neighborhood $U$ introduced in equation~\eqref{eq:U}, the remainder $\mathcal{R}$ is so small that it can be neglected for practical purposes, and therefore we have a complete and explicit description of the local motion around the equilibrium point.

In fact, the BNF of order $M$ can only be used in the absence of resonances between the characteristic exponents of the equilibrium point of order less than or equal to $M$, that is, when $\abs{k_0\lambdaL+\iu k_\pr\nupL+\iu k_\vr\nuvL}$ is not too small for integers $k_0$, $ k_\pr$, $ k_\vr$ such that $\abs{k_0}+\abs{k_\pr}+\abs{k_\vr}\leq M$, which is what will be assumed in this paper. In case of resonances of small order, a resonant normal form should be used. See more details in Section~\ref{sec:BNF}.

\begin{rem}
\MModify{
In the case of the Spatial Restricted 3-Body Problem (SR3BP), for the saddle-center-center equilibrium point considered here, one can take one  the collinear libration points $L_{1,2}$, and compute the BNF around that point.  (For $L_3$, the two center frequencies are close to a  \(1:1\) resonance, so one typically computes  a resonant Birkhoff normal form instead.)
Numerical verifications at the libration point $L_{1}$ of the SR3BP---rather than purely theoretical bounds---show that truncating the normal form to order $M=16$ guarantees an error $\errorBNF$~\eqref{eq:errorBNF} between the true dynamics and the dynamics induced by the BNF that is less than $3 \times 10^{-7}$ for $\delta=0.1$ and $r=0.05$, at least over short time intervals $t<3$. Further details are provided in Section~\ref{sec:NumErrorSEBNF}, where we illustrate that the analytical bounds for the remainder $R$ are overly pessimistic compared to the numerical estimates for $L_{1}$. This discrepancy is a well-known feature documented by several authors (see, e.g.,~\cite{EfthyGC04,CaraccioloL20}).}
\end{rem}

Let us now be more precise. If we neglect the remainder $\mathcal{R}$, the flow inside $U$ is simply given by the normal form $\mathcal{Z}(\bI)$, with $\bI=(I_\hr,\bJ)$, where $I_\hr=xy$, $\bJ=\left(I_\pr,I_\vr\right)$ are conserved quantities.
\PauModify{
Abusing notation, we will also denote the BNF flow by $\flow_{t}$.} 
It takes the explicit form
\begin{equation}\label{eq:flow_NF}
x(t)=\ee^{t\lambda} x(0),\quad
y(t)=\ee^{-t\lambda} y(0), \qquad
\bJ(t)=\bJ(0),\quad
\bphi(t)=\bphi(0)+t\bnu,
\end{equation}
with $\bphi=\left(\phi_\pr,\phi_\vr\right)$, where we have introduced the notation
\begin{align}
\lambda&\coloneq \frac{\partial \mathcal{Z}}{\partial I_\hr}=\frac{\partial \mathcal{Z}}{\partial I_\hr}(\bI)=\lambdaL +O(\bI)
=\frac{\partial \mathcal{Z}}{\partial I_\hr}(\bI(0)),\label{eq:lambda}\\
\bnu&\coloneq\frac{\partial \mathcal{Z}}{\partial \bJ}=\nabla_{\bJ} \mathcal{Z}(\bI)
=\bnuL+O(\bI)=\nabla_{\bJ} \mathcal{Z}(\bI(0)).\label{eq:nu}
\end{align} 
Notice that, restricting $\delta$ and $r$ if necessary, $\lambda>0$ for $-\delta^2 \leq I_\hr\leq \delta^2$, $0\leq I_\pr, I_\vr\leq r$. 

Thanks to the explicit expression~\eqref{eq:flow_NF} for the trajectories in $U$, the description of the most important invariants sets within $U$, the NHIM, their invariants manifolds, as well as the rest of trajectories, is straigthforward.

First, let us denote 
\[
\NHIM=W_{\mathrm{loc}}^{\cc}=\left\{\bX=(x=0, y=0, \bJ, \bphi\right\}=\left\{(0, 0, \bJ,\bphi)\right\}
=\{(0,0)\} \times[0, r]^{2}\times\T^2
\]
the local center manifold of $L$ provided by the normal form, which is a 4-dimensional NHIM. 
The dynamics inside $\NHIM$ is simply determined by $\dot\bJ=0,\dot\bphi=\bnu$, where $\nu=\nu(0,I)$, and given by $\Phi_t(\bX)=\left(0,0,\bJ,\phi+t\nu\right)$.

Any point $\bX=(0,0, \bJ,\bphi)\in\NHIM$ has associated \emph{local stable and unstable fibers} in $U$
\[
W_{\mathrm{loc}}^\st\left(\bX\right)=\left\{\bZ=(x=0, y, \bJ,\bphi),\ |y|\leq \delta\right\},
\ 
W_{\mathrm{loc}}^\un\left(\bX\right)=\left\{\bZ=(x, y=0, \bJ,\bphi),\ |x|\leq \delta\right\},
\]
which can be extended outside $U$ by the backward and forward flow, respectively, to form the global stable and unstable fibers $W^\st\left(\bX\right)$ and  $W^\un\left(\bX\right)$. 
Their unions on $\bX\in\NHIM$ form the stable and unstable \emph{invariant} manifolds of the NHIM~$\NHIM$
\begin{equation*}
W^\st\left(\NHIM\right)=\bigcup\limits_{\bX\in\NHIM} W^\st\left(\bX\right),
\qquad
W^\un\left(\NHIM\right)=\bigcup\limits_{\bX\in\NHIM} W^\un\left(\bX\right).
\end{equation*}

The dynamics inside $W_{\mathrm{loc}}^\st\left(\NHIM\right)$ is simply determined by \[\dot y=-\lambda y,\quad \dot\bphi=\bnu,\] where $\lambda=\lambda(0,I)$ and $\nu=\nu(0,I)$, and given by $\Phi_t(\bZ)=\left(0,\ee^{-\lambda t}y,\bJ,\phi+t\nu\right)$ for $Z=(0,y,\bJ,\bphi)\in W_{\mathrm{loc}}^\st\left(\NHIM\right)$. With the same $\lambda=\lambda(0,I)$ and $\nu=\nu(0,I)$, the dynamics inside $W_{\mathrm{loc}}^\un\left(\NHIM\right)$ is determined by \[\dot x=\lambda x,\quad \dot\bphi=\bnu,\] and given by $\Phi_t(\bZ)=\left(\ee^{\lambda t} x,0,\bJ,\phi+t\nu\right)$ for $Z=(x,0,\bJ,\bphi)\in W_{\mathrm{loc}}^\un\left(\NHIM\right)$.

\begin{figure}[ht]
\centering
\includegraphics[width=0.9\textwidth]{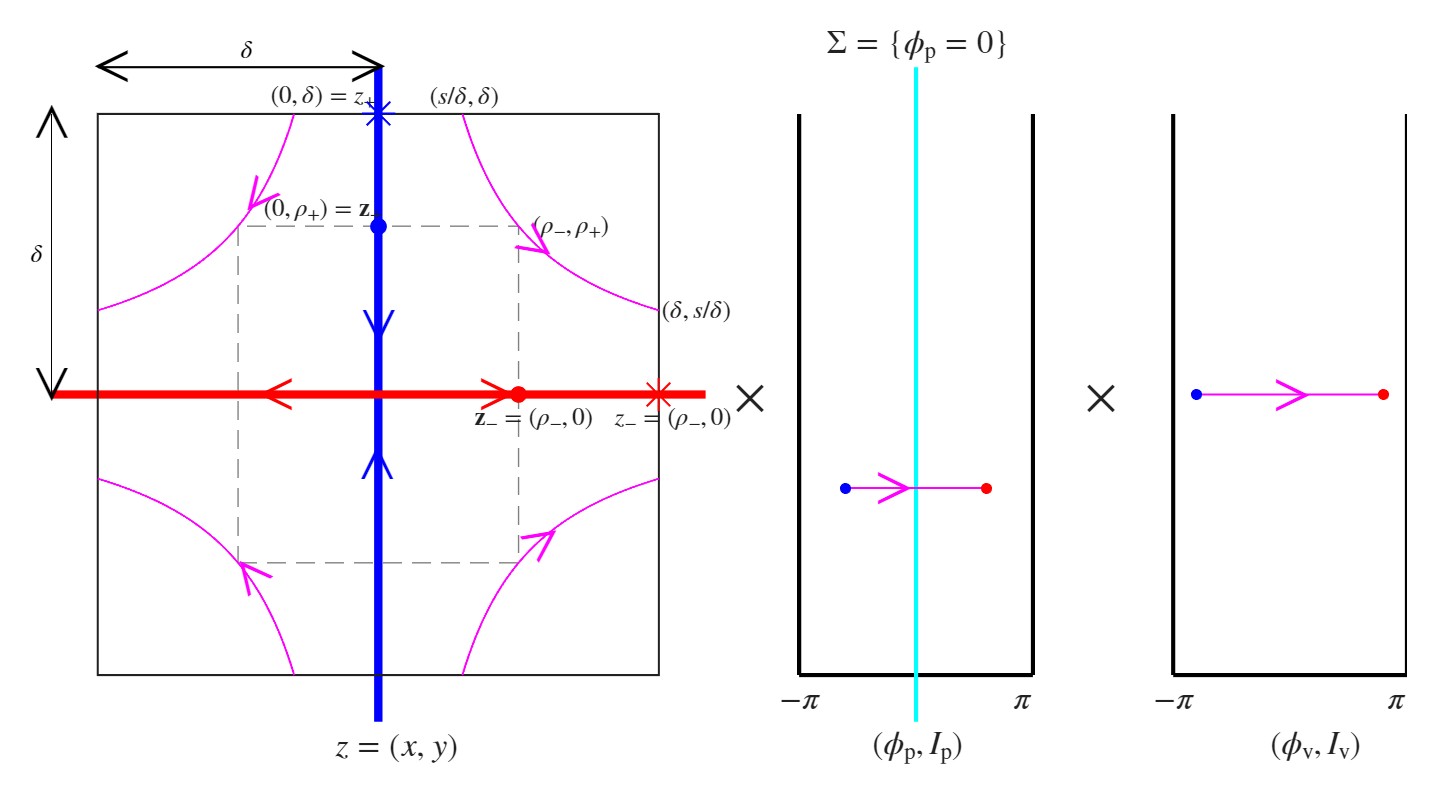}
\caption{$\bZ_+=(z_+,\bJ,\bphi)=(0,\delta,\bJ_+,\bphi_+)$,
$\bK_+=(s/\delta,\delta,\bJ_,\bphi_+)$,
$\bZ_-=(z_-,\bJ,\bphi)=(\delta,0,\bJ_-,\bphi_-)$,
$\bK_-=(\delta,s/\delta,\bJ_,\bphi_-)$,
$\bz_+=(\textbf{z}_+,\bJ,\bphi)=(0,\rho_+,\bJ,\phi)$,
$\bz_-=(\textbf{z}_-,\bJ,\bphi)=(\rho_-,0,\bJ,\bphi)$,
$\bk=(\rho_-,\rho_+,\bJ,\bphi)$.
In this figure the trajectory from $\bK_+$ to $\bK_-$ only crosses once the section $\SecLocal=\{\phi_\pr=0\}$, so that
$\bk=\bk_+=\bk_-$. See the text for more information.}
\label{fig:NearbyTrajectories}
\end{figure}

\begin{defn}[Sections]
    We define the following local Poincaré sections for the flow. The \emph{incoming boundary} is given by 
    $\SecPlus\coloneqq\{\bK=(x,y,\bJ,\bphi)\in \partial U : y=\pm\delta\}$, while the \emph{outgoing boundary} is given by 
    $\SecMinus\coloneqq\{\bK=(x,y,\bJ,\bphi)\in \partial U : x=\pm\delta\}$.
    \MModify{We also define the section $\SecLocal=\{\bK=(x,y,\bJ,\bphi)\in \partial U : \phi_p=0 \}$.}
\end{defn}

Finally, consider general points $\bK=(x,y,\bJ,\bphi)\in U$ outside the local stable and unstable manifolds $W_{\mathrm{loc}}^{\st,\un}\left(\NHIM\right)$, which implies that $xy\neq 0$. Inside $U$, the dynamics are governed by the differential equations \[\dot x=\lambda x,\quad \dot y=-\lambda y,\textrm{ and }\dot\bphi=\bnu.\] The corresponding flow is given by $\Phi_t(\bK)=\left(\ee^{\lambda t} x,\ee^{-\lambda t} y,\bJ,\bphi+t\bnu\right)$, which confirms that both $I_{\hr}=xy$ and $\bJ$ are constants of motion. For $s=xy\neq 0$, we will show that the domain $U$ is stratified into several zones, depending of the number of times the trajectories intersect the section $\SecLocal$. To avoid tracking sign combinations throughout, we first restrict our explicit formulas to the first quadrant $x>0, y>0$, and later adapt equations~\eqref{eq:OneCut} to the remaining three quadrants.

The trajectory $\Phi_t(\bK)$ enters and exits $U$ by connecting the incoming and outgoing boundary sections (see Figure~\ref{fig:NearbyTrajectories}): it starts at a point $\bK_+=(s/\delta,\delta,\bJ,\bphi_+)\in\SecPlus$, passes through $\bK$, and ends at $\bK_-=(\delta, s/\delta,\bJ,\bphi_-)\in\SecMinus$. The entry and exit angles are linked by
$\bphi_-=\bphi_+ +\Tint\bnu$,  where $\bnu=\bnu(s,\bJ)$.
The \emph{internal transit time} $\Tint$ within $U$, associated with the Poincaré map $\SecPlus\to\SecMinus$ of the Hamiltonian flow, is given by 
\begin{equation}
\label{eq:Tint}
\Tint=\Tint(\delta,s,\bJ)=\frac{1}{\lambda(s,\bJ)}\,\ln\dfrac{\delta^2}{s} 
=\frac{1}{\lambdaL}\,\ln\dfrac{\delta^2}{s} \left(1+O(s,\bJ)\right).
\end{equation}

Furthermore, the starting and ending points can be expressed from $\bK$ as $\bK_+=\Phi_{\Tint{}_{,+}}\left(\bK\right)$ and $\bK_-=\Phi_{\Tint{}_{,-}}\left(\bK\right)$, where the partial transit times are given by $\displaystyle \Tint{}_{,+}=\frac{1}{\lambda}\ln\frac{\delta}{y}$ and $\displaystyle \Tint{}_{,-}=\frac{1}{\lambda}\ln\frac{\delta}{x}$, with $\lambda=\lambda\left(s,\bJ\right)$. These partial times reconstruct the total transit time via $\Tint{}_{,+}+\Tint{}_{,-} = \Tint$. Notably, all three time durations are independent of the angular variables $\bphi_+$, $\bphi$, and $\bphi_-$.

When working later with the four-dimensional reduced phase space $\redSpace=U\cap\{H=h\}\cap\SecLocal$ introduced in~\eqref{eq:ReducedSpace}, we will require $\bK$ to be inside $\redSpace$ and outside $\SecPlusMinus$. We will denote this point as $\bk=(x,y,I,\phi)\in\SecLocal\setminus\SecPlusMinus$, with $I_{\pr}=f_h(0,I_{\vr})$ and $\phi_{\pr}=0$. Again, the trajectory through $\bk$ starts at $\bK_+$ and ends at $\bK_-$, but it may intersect $\SecLocal$ one or more times. 

To quantify this phenomenon, we introduce the following.
\begin{defn}[Transit times $t_+, t_-$] \label{defn:tplus_tminus}
Let $t_+=t_+\left(s,I,\phi_{\pr,+}\right)>0$ denote the first time at which the forward trajectory starting at $\bK_+=(s/\delta,\delta,\bJ,\bphi_+)\in\SecPlus$ intersects $\SecLocal$ 
at a point $\bk_+$. This corresponds to the transit time of the primary Poincaré map $\SecPlus\to\SecLocal$. 
Analogously, let $-t_-=-t_-\left(s,I,\phi_{\pr,-}\right)<0$ denote the first time at which the (backward) trajectory starting at $\bK_-=(\delta,s/\delta,\bJ,\bphi_-)\in\SecMinus$ intersects $\SecLocal$ at a point $\bk_-$. This represents the transit time of the primary Poincaré map $\SecLocal\to\SecMinus$. 
\end{defn}

Clearly, we have 
\begin{equation}
\label{eq:tplusminus}
0 < t_+ = \frac{2\pi-\phi_{\pr,+}}{\nup } < \frac{2\pi}{\nup },\quad  0 < t_- = \frac{\phi_{\pr,-}}{\nup } < \frac{2\pi}{\nup},
\ \text{ for }\phi_{\pr,\pm}\in(0,2\pi),
\end{equation}
where $\nup =\nup (s,I)$ is given in~\eqref{eq:nu}.

The trajectory starting at $\bK_+\in\SecPlus$ and ending at $\bK_-\in\SecMinus$ intersects $\SecLocal\setminus\SecPlusMinus$ exactly once if and only if the following \emph{one-cut condition} holds
\begin{equation}\label{eq:OneCut0}
    t_+  <  \Tint  \leq t_+ + \frac{2\pi}{\nu_p},
\end{equation}
where $\Tint$ is the internal transit time~\eqref{eq:Tint}.

Indeed, the inequality $t_+  <  \Tint$ ensures that the first intersection occurs before $\Tint$, whereas $\displaystyle \Tint  \leq t_+ + \frac{2\pi}{\nu_p}$ guarantees that the second intersection takes place after $\Tint$. Under this condition, it also follows that $\bk_+=\bk=\bk_-$ and $t_+ + t_-=\Tint$.

Substituting the expression for $\Tint$ given in~\eqref{eq:Tint}, the one-cut condition~\eqref{eq:OneCut0} is equivalent to
\begin{equation}
\label{eq:OneCut}
   \ee^{\displaystyle -2\pi\lambda/\nup} <\frac{s}{\delta^2}\ee^{\displaystyle \lambda t_+}\leq 1,
\end{equation}   
with $t_+$ given in~\eqref{eq:tplusminus} and $\lambda=\lambda(s,I)$, $\nup=\nup(s,I)$ in~(\ref{eq:lambda}-\ref{eq:nu}).

\begin{rem}
Although condition~\eqref{eq:OneCut} for a single intersection with $\SecLocal\setminus\SecPlusMinus$ might seem restrictive, it is easily satisfied in practice. For instance, when $\lambda$ and $\nup $ are of the same order of magnitude ($\lambda/\nup =O(1)$), we find that $\exp\left(-2\pi\lambda/\nup \right)=O(10^{-3})$. Thus, across almost the entire square $\abs{x},\abs{y}\leq \delta$, trajectories inside $U$ cross $\SecLocal\setminus\SecPlusMinus$ exactly once. Trajectories experiencing multiple intersections are confined to a tiny region extremely close to the invariant axes $x=0$ and $y=0$.
\end{rem}

For trajectories with multiple intersections--specifically, $\ell\geq 1$ additional crossings  (yielding $\ell+1\geq 2$ intersections in total)--we shift $t_+$ in~\eqref{eq:OneCut} to $t_+ +\dfrac{2\pi\ell}{\nup }$. This adjustment accounts for the fact that the transit time of the primary Poincaré map $\SecLocal\to\SecLocal$ is $\displaystyle \frac{2\pi}{\nup }$. Consequently, the  \emph{$(\ell+1)$-cuts condition} becomes
\begin{equation*}
   \ee^{\displaystyle -2\pi\lambda/\nup}  \ee^{\displaystyle -2\pi\ell\lambda/\nup} <\frac{s}{\delta^2}\ee^{\displaystyle \lambda t_+}\leq \ee^{\displaystyle -2\pi\ell\lambda/\nup}.
\end{equation*} 

To fully characterize the square $\abs{x},\abs{y}\leq\delta$ away from the invariant manifolds $x=0$ and $y=0$, there remains a 0-cut zone in the $(x,y)$ plane. This zone corresponds to trajectories where $t_+\geq\Tint$ or, equivalently, when the \emph{0-cuts condition} holds:  
\begin{equation*}
\ee^{-\displaystyle \lambda t_+}< \frac{s}{\delta^2}\leq 1. 
\end{equation*}

\section{Transition map in BNF coordinates} \label{sec:transition_map_BNF}

Let $U$ be the BNF domain given by Equation~\eqref{eq:U}, with incoming/outgoing boundary sections $\SecPlusMinus\in\partial U$.
We think of $\textrm{int } U$ as the \emph{interior} domain, containing the NHIM $\NHIM$, and $\textrm{out } U$ as the \emph{exterior} domain, containing the homoclinic channel $\Gamma$.

This section explains how the transition map is computed in a non-perturbative setting, using Equation~\eqref{eq:transition_map}.
In that equation, $\flow_{T_+}\circ \flow_{T_-}$ corresponds to homoclinic orbit segments from the outgoing boundary section $\SecMinus$ to the incoming boundary section $\SecPlus$, while the local wave maps $\Omega_\pm^{\Gamma_\pm}$ become simply standard projections in BNF coordinates.
After that, the dynamics is restricted to a 4-dimensional \emph{reduced phase space} $\redSpace$ by fixing the value of the energy and using a Poincaré section. The dynamics is described by the first return map $\frm$ of the Hamiltonian flow. We derive the \emph{reduced transition map} of $\frm$ from the transition map of the flow.


\subsection{Transition map}
We recall how the homoclinic orbits are computed following~\cite{DGR16}. Let $\bZ_-\in W^u(\NHIM)$ and $\bZ_+\in W^s(\NHIM)$ be two arbitrary points in the outgoing (resp. incoming) boundary section of the BNF domain $U$, with coordinates
\begin{align*}
\bZ_- &= (\delta, 0, \bJ_-, \bphi_-) \in \SecMinus \subset \partial U \\
\bZ_+ &= (0, \delta, \bJ_+, \bphi_+) \in \SecPlus \subset \partial U. \end{align*}
(For simplicity, we consider the case $x=y=\delta>0$, but $x$ or $y=-\delta<0$ is analogous). The variables $(\bJ_-,\bphi_-)$ and $(\bJ_+,\bphi_+)$ determine $\bZ_-$ and $\bZ_+$, and they are parameters to be found during the procedure.

Numerically integrate the trajectory starting at $\bZ_-$ forward in time until it intersects some exterior section $\SecGlobal$ at a point $\bzeta_-=\flow_{T_-}(\bZ_-)$. Similarly, integrate the trajectory starting at $\bZ_+$ backward in time until it intersects $\SecGlobal$ at a point $\bzeta_+=\flow_{T_+}(\bZ_+)$.
To find a homoclinic, iteratively correct $(\bJ_-,\bphi_-)$ and $(\bJ_+,\bphi_+)$ using a Newton-like method, until $\bzeta_+=\bzeta_-=\bzeta$. 

Notice that $[\bZ_-, \bZ_+]$ is a finite time homoclinic orbit segment of length $\Text = T_- + T_+$. The exterior Poincaré map $\Pext=\flow_{T_+}\circ \flow_{T_-}$ takes $Z_-$ to $Z_+$.
The times $T_\pm=T_\pm(\bzeta)$ are continuous functions on $\Gamma$. (See Remark~\ref{rem:transition_alternative}).

\begin{rem}
	These type of homoclinic intersections at the ``first cut'' of the asymptotic manifolds with $\SecGlobal$ are usually called \emph{primary} or \emph{one-round homoclinics}.
\end{rem}

\begin{rem}
	In practice, to determine the homoclinic channels, one can compute all the (primary) homoclinics from $\NHIM$ to $\NHIM$, and then identify continuous families of homoclinics with the property that the restrictions $\Omega_\pm^\Gamma$ are  diffeomorphisms onto their images.
	In~\cite{DGR16}, we illustrate this computation in the case of the SR3BP around the equilibrium point $L_1$, with small enough energy. We show that there are two different homoclinic channels $\Gamma_1, \Gamma_2$. By suitably restricting $\Gamma_1, \Gamma_2$, the associated scattering maps $\SM_1, \SM_2$ are defined in a common domain.
\end{rem}

Assume that $\bzeta$ belongs to a homoclinic channel $\Gamma$ as defined in Section~\ref{sec:scattering_map_transition_map}. 
By construction, $\bZ_+ \in W^{\st}(\bX_{+})$, where $\bX_+ = (0, 0, \bJ_+, \bphi_+) \in \NHIM$. Indeed, using the BNF equations~\eqref{eq:flow_NF}, where $\lambda = \lambda(0, \bJ_{+}) > 0$ and $\bnu = \nu(0, \bJ_{+})$ are given by equations~(\ref{eq:lambda}--\ref{eq:nu}), the forward flows read
\begin{align*}
\Phi_{t}(\bZ_{+}) &= \bigl(0, \ee^{-t\lambda}\delta, \bJ_{+}, \bphi_+ + t\bnu\bigr), \\
\Phi_{t}(\bX_{+}) &= \bigl(0, 0, \bJ_{+}, \bphi_+ + t\bnu\bigr).
\end{align*}
Consequently, their difference satisfies
\[
\Phi_{t}(\bZ_{+}) - \Phi_{t}(\bX_+) = \bigl(0, \ee^{-t\lambda}\delta, \mathbf{0}, \mathbf{0}\bigr) \longrightarrow \mathbf{0} \quad \text{as } t \rightarrow \infty.
\]
Thus, in BNF coordinates, the wave map $\Omega_+^{\Gamma_+}$ is simply the standard projection onto $\NHIM$:
\[
\bZ_{+} = (0, \delta, \bJ_{+}, \bphi_{+}) \in W_{\mathrm{loc}}^\st(\NHIM) \mapsto \bX_+ = (0, 0, \bJ_{+}, \bphi_{+}) \in \NHIM.
\]

\begin{rem}
	$\bZ_{+}$ ``synchronizes in the far future'' with $\bX_+$. The wave map $\bZ_{+} \mapsto \bX_+$ maps a point to its future synchronized state on the NHIM. The stable manifold $W^{\st}(\bX_{+})$ is also referred to as the \emph{stable isochrone manifold} of $\bX_+$. 
\end{rem}


Analogously, $\bZ_- \in W^{\un}(\bX_{-})$ where $\bX_- = (0, 0, \bJ_{-}, \bphi_{-}) \in \NHIM$. 
Thus, in BNF coordinates, the wave map $\Omega_-^{\Gamma_-}$ is simply the standard projection onto $\NHIM$:
\[
\bZ_{-} = (\delta, 0, \bJ_{-}, \bphi_{-}) \in W_{\mathrm{loc}}^\un(\NHIM) \mapsto \bX_- = (0, 0, \bJ_{-}, \bphi_{-}) \in \NHIM.
\]

Finally, the transition map $\TM\colon\ \NHIM \to \NHIM$, $\bX_+ = \TM(\bX_-)$ is given by the composition~\eqref{eq:transition_map} of the exterior Poincaré map $\Pext = \flow_{T_+} \circ \flow_{T_-}$ with the local wave maps $\Omega_\pm^{\Gamma_\pm}$:
\[
\bX_{-} = (0, 0, \bJ_{-}, \bphi_{-}) \mapsto 
\bZ_{-} = (\delta, 0, \bJ_{-}, \bphi_{-}) \mapsto
\bZ_{+} = (0, \delta, \bJ_{+}, \bphi_{+}) \mapsto 
\bX_{+} = (0, 0, \bJ_{+}, \bphi_{+}),
\]
or simply
\[
\bX_{-} = (\bJ_{-}, \bphi_{-}) \mapsto 
\bX_{+} = (\bJ_{+}, \bphi_{+}).
\]

\begin{rem} 
In Definition \ref{defn:scattering_map_flow} the scattering map depends only on the choice of the homoclinic channel $\Gamma$, while in Definition~\ref{defn:transition_map_flow} the transition map
additionally depends on the times $T_\pm$, which in turn depend on the sections $\SecPlusMinus$.  
	On the one hand, the definition of scattering map is conceptually simpler than the transition map.
	On the other hand, the numerical computation of the transition map is more straightforward. 
	To obtain the scattering map, one would first compute the transition map $\bX_-\mapsto \bX_+$, and then derive the scattering map using Equation \eqref{eq:SMAndTM}. 
\end{rem}

\begin{rem}
	Note again that in formula~\eqref{eq:SMAndTM} relating the scattering map and the transition map, the flow $\Phi_{T_{\pm}}$ is only needed on the NHIM $\NHIM$ and therefore only involves translations for the angles $\bphi=(\phi_\pr,\phi_\vr)$ in the BNF variables, that is, the \emph{inner flow} inside the NHIM $\NHIM$ where the actions  $I_\hr=0$, $I=(I_\pr,I_\vr)$ remain constant.
\end{rem}

\subsection{Reduced transition map}
Inside the BNF domain $U$, we can pass from the flow to a discrete-time dynamical system $\frm$, given by the first return map to a Poincaré section $\SecLocal$. The map $\frm$ inherits an associated NHIM $\NHIMhred$ that is easily obtained from $\NHIM$. Moreover, we derive the transition map $\tM\colon \NHIMcred\to \NHIMcred$ of $\frm$ from the transition map $\TM\colon \NHIM\to \NHIM$ of $\flow_t$. 
Since $\tM$ is lower-dimensional, we refer to it as the \emph{reduced transition map}. In what follows, we derive the expressions of these transition maps in BNF coordinates.

The original Hamiltonian system has a conserved quantity, the energy $H$, which is the energy $\mathcal{Z}$ of the BNF when the remainder is neglected. For a fixed energy $H\circ\nftrans=h>\hL$, the action $I_\pr$ can be given as an implicit function of $I_\hr,I_\vr, h$: $I_\pr=(h-\hL-I_\hr-\nuvL I_\vr)/\nupL+O_2(I_\hr,I_\vr)=:\Fh(I_\hr,I_\vr)$.
The corresponding 3-dimensional NHIM on the level of energy $H=h$ is given by 
$\NHIMh=\NHIM\cap \{H\circ\nftrans=h\}$ and is parameterized by the variables $\phi_\pr,\phi_\vr,I_\vr$. Taking it one step further, we can consider the Poincar\'e section 
\[
\SecLocal  \coloneq  
\left\{ (x,y,\bJ,\bphi) \in U ~:~ \phi_\pr = 0 \right\},
\]
which is a global section within $U$, and denote 
\[
\NHIMh^\SecLocal =\NHIMh\cap\SecLocal =\left\{\bx=(0,0,\bJ,\bphi): \phi_\pr=0,I_\pr=\Fh(0,I_\vr), I_\vr\in[0,r], \phi_\vr\in\T\right\}
\]
the associated 2-dimensional NHIM of the return map, which is parameterized by the variables $\phi_\vr,I_\vr$.

Next, we compute the 2-dimensional reduced transition map $\tM\colon \NHIMhred \to \NHIMhred$ associated to $\frm$, which is exact symplectic (see Section~\ref{sec:scattering_map_transition_map}).
Recall from the above construction that $\bZ_\pm \in \SecPlusMinus$  satisfy $\bZ_\pm \in W^{\st/\un}(\bX_{\pm})$, where $\bX_\pm = (0, 0, \bJ_{\pm}, \bphi_{\pm}) \in \NHIM$.
Let $\tT_\pm = T_\pm + t_\pm \geq T_\pm$ be the first times such that the forward and backward trajectories of $\bzeta\in\Gamma$ reach the  section $\SecLocal = \{\phi_\pr = 0\}$ inside $U$.
In other words, $t_+, t_- > 0$ are the transit times of $\bZ_\pm\in\SecPlusMinus$ to the Poincaré section $\SecLocal$. (See Definition~\ref{defn:tplus_tminus}).

Define new points 
\begin{equation*}
	\begin{aligned}
		\bz_+ &= \flow_{t_+}(\bZ_+) 
		= \bigl(0, \rho_+, \bJ_{+}, \widetilde{\bphi}_{+}\bigr)
		:= \bigl(0, \ee^{-t_+\lambda}\delta, \bJ_{+}, \bphi_+ + t_+\bnu\bigr) \in \SecLocal, \\
		\bz_- &= \flow_{-t_-}(\bZ_-)
		= \bigl(\rho_-, 0, \bJ_{-}, \widetilde{\bphi}_{-}\bigr) 
		:= \bigl(\ee^{t_-\lambda}\delta, 0, \bJ_{-}, \bphi_- - t_-\bnu\bigr) \in \SecLocal.
	\end{aligned}
\end{equation*}
Notice that $\bz_\pm$ belong to the reduced phase space $\redSpace$, that is, their planar angle and action satisfy $\widetilde{\phi}_\pr=0$, $I_\pr=\Fh(I_\hr,I_\vr)$.
See Figure~\ref{fig:NearbyTrajectories}.

Due to the conjugacy of invariant fibers $W^{\st/\un}(x)$ by the flow, we have
$ \bz_+ \in W^{\st}(\bx_{+})$ and $\bz_- \in W^{\un}(\bx_{-})$, where
\begin{equation*}
	\begin{aligned}
		\bx_+ &= \flow_{t_+}(\bX_+) 
		= \bigl(0, 0, \bJ_{+}, \widetilde{\bphi}_{+}\bigr)
		= \bigl(0, 0, \bJ_{+}, \bphi_+ + t_+\bnu\bigr) \in \NHIMhred, \\
		\bx_- &= \flow_{-t_-}(\bX_-) 
		= \bigl(0, 0, \bJ_{-}, \widetilde{\bphi}_{-}\bigr) 
		= \bigl(0, 0, \bJ_{-}, \bphi_- - t_-\bnu\bigr) \in \NHIMhred,
	\end{aligned}
\end{equation*}
that is, their planar angle and action satisfy $\phi_\pr=0,I_\pr=\Fh(0,I_\vr)$.

Thus, for the discrete system, in BNF coordinates, the wave maps $\Omega_\pm$ are also the standard projection onto $\NHIMhred$, and the reduced transition map $\tM\colon \NHIMhred\to\NHIMhred$, $\bx_{+}=\tM(\bx_{-})$ is given by
\[ \bx_{-} = (I_{\vr,-},\widetilde{\phi}_{\vr,-}) \mapsto  \bx_+ = (I_{\vr,+},\widetilde{\phi}_{\vr,+}). \]

\begin{rem}
\label{rem:time_dependency_tM}
    The reduced transition map depends on the choice of $\SecLocal$, a global section within $U$, which in turn uniquely determines  
    the times $t_\pm$. These times must satisfy $\phi_{\pr_\pm}\equiv 0\pmod {2\pi}$, so there is an infinite but discrete set of them. However, any two reduced transition maps $\tM=\tM(\tT_\pm)$ and $\tM'=\tM'\left(\tT'_\pm\right)$ are related in the following way
    \[
    \tM = \Phi_{-\left(\tT'_+ - \tT_+\right)} \circ \tM' \circ \Phi_{-\left(\tT'_- - \tT_-\right)}. 
    \]
    If $\tM$ maps the point $\bx_-=\left(\bJ_{\vr,-},\phi_{\vr,-}\right)$ to $\bx_+=\left(\bJ_{\vr,+},\phi_{\vr,+}\right)$, where, for simplicity of notation, we only write the components $\left(\bJ_{\vr},\phi_{\vr}\right)$ of $\bx_\pm$, which are the signicant ones, then $\tM'$ maps the point
    $\bx'_-=\left(\bJ'_{\vr,-},\phi'_{\vr,-}\right) =\left(\bJ_{\vr,-},\phi_{\vr,-} - \left(\tT'_- - \tT_-\right)\nuv\left(\bJ_{\vr,-}\right)\right)$ to
    $\bx'_+=\left(\bJ'_{\vr,+},\phi'_{\vr,+}\right) =\left(\bJ_{\vr,+},\phi_{\vr,+} + \left(\tT'_+ - \tT_+\right)\nuv\left(\bJ_{\vr,+}\right)\right)$.
   Thus, both pre-images $\bx_-$ and $\bx'_-$ have the same action $\bJ_{\vr,-}$ but their angle differs by a shift.
    Similarly, both images $\bx_+$ and $\bx'_+$ have the same action $\bJ_{\vr,+}$ but their angle differs by a shift.
\end{rem}


\noindent\textbf{Notation warning:} When clear from the context, the `reduced transition map' and `reduced scattering map' are abbreviated just as `transition map' and `scattering map'. 

\section{Pseudo-trajectories and shadowing trajectories} \label{sec:pseudo_orbits_shadowing_trajectories}

This section defines in a precise way several key concepts that we will later need to construct drift trajectories for the original Hamiltonian system $\flow_t$. In particular, we introduce 
\begin{itemize}
    \item a \emph{guiding orbit} \MModify{as an orbit of an IFS consisting of the inner dynamics and (reduced) transition maps;}
   
    \item a \emph{pseudo-trajectory} as a sequence of trajectory segments of the flow $\flow_t$ that form an almost continuous trajectory, and
    
    \item a \emph{shadowing trajectory} as a true trajectory of the flow $\flow_t$ that is ``close'' to a pseudo-trajectory.
\end{itemize}
The reader is referred to Section~\ref{sec:semianalytical_results} for an illustration of these concepts in the SR3BP system. In particular, see an example of guiding orbits in Figure~\ref{fig:diffusion_shortest_path}, a pseudo-trajectory in Figure~\ref{fig:pseudo_orbit_refined}, and the corresponding shadowing trajectory in Figure~\ref{fig:shadowing_trajectory}.

\begin{defn}
\MModify{Let $\tM_1, \tM_2, \dotsc, \tM_{L}$ be reduced transition maps defined on $\NHIMhred$.
Any finite orbit of the IFS $\{f,\tM_1,\ldots, \tM_L\}$ consisting of the inner map and transition maps will be referred to as a \emph{guiding orbit}.} 
\end{defn}

\MModify{In particular, we will consider guiding orbits of a special type, generated solely by transition maps, i.e., $(\bx_i)_{i=1}^{N}$ in $\NHIMhred$ such that $\tM_{\alpha(i)}(\bx_i)=\bx_{i+1}$, where $\alpha(i)\in\{1,2,\dotsc,L\}$. 
These guiding orbits are important since they minimize the drift  time  (when $\delta=0.1$); see Section \ref{subsec:Shortest-time_guiding}. 
Therefore,  the main theoretical results 
Theorems \ref{th:shadowing} and \ref{th:shadowing_topological}, as well as the main numerical result Theorem \ref{thm:main_numerical}, 
are concerned with shadowing of guiding orbits of this type. 
}   


Recall that transition maps encode the homoclinic dynamics to $\NHIMhred$, so the points $\bx_i$ in a guiding orbit do not belong to a true trajectory of the original system. To each guiding orbit, we associate a \textit{sequence of $\NHIMh$-homoclinic trajectories} such that each homoclinic synchronizes with the next homoclinic in the sequence in the following way.
As explained in Section~\ref{sec:transition_map_BNF}, for any homoclinic point $\bzeta\in W^\un(\NHIM)\cap W^\st(\NHIM)$, there exist 
\begin{itemize}
    \item unique times $T_\pm>0$ such that $\bZ_\pm=\Phi_{\pm T_\pm}(\bzeta)\in\SecPlusMinus\subset \partial U$, and 
    \item unique times $\tT_\pm=T_\pm+t_\pm$ with $t_\pm>0$, such that $\bz_\pm=\Phi_{\pm \tT_\pm}(\bzeta) \in U_h^\Sigma$ has zero planar angle: $\phi_{\pr,\pm}=0$ in $\bz_\pm$.
\end{itemize}

\begin{rem}\label{rem:Zpm_to_tZpm}
$\bZ_\pm$ and $\bz_\pm$ are simply related by the flow: $\bz_\pm = \Phi_{\pm t_\pm}(\bZ_\pm)$.    
\end{rem}

Let the guiding orbit be fixed. To each iterate
\begin{equation}\label{eq:Xi}
\bx_{i}=(0,0,\bJ^i,\bphi^i)\in\NHIMhred, 
\text{ with } I^i_\pr=\Fh(0,I^i_\vr), \phi^i_\pr=0,
\end{equation}
in the guiding orbit, associate the points
\begin{equation} \label{eq:tZi}
\begin{aligned}
    \bz_i^+ &= \left(0,y= \rho_+^i := \exp \left( -t_+^i\,\lambda\right)\delta, I^i, \phi^i\right)\in W^\st( \bx_i)\subset W^\st(\NHIMhred ),\\
    \bz_i^- &= \left(x= \rho_-^i := \exp \left( -t_-^i\,\lambda \right)\delta,0, I^i, \phi^i\right)\in W^\un( \bx_i)\subset W^\un(\NHIMhred ),
\end{aligned}
\end{equation}
and
\begin{equation} \label{eq:Zi}
\begin{aligned}
    \bZ_i^+ &= \Phi_{-t_+^i}(\bz_i^+) =\left(0,y=\delta,I^i, \phi^i_+=\phi^i- t_+^i\, \bnu \right)\in W^\st(\Phi_{-t_+^i}(\bx_i))\cap \SecPlus, \\
    \bZ_i^- &=\,\, \Phi_{t_-^i}(\bz_i^-)\: =\left(x=\delta,0,I^i,\phi^i_-=\phi^i+ t_-^i\, \bnu \right)\in W^\un(\Phi_{t_-^i}(\bx_i))\cap \SecMinus, 
\end{aligned}
\end{equation}
where $\lambda$, $\bnu$, defined in equations~(\ref{eq:lambda}--\ref{eq:nu}),  are evaluated on $\bI=(0,I^i)=(0,I_{\pr}^i,I_{\vr}^i)$.

\begin{figure}[!ht]
\centering
\includegraphics[width=0.9\textwidth]{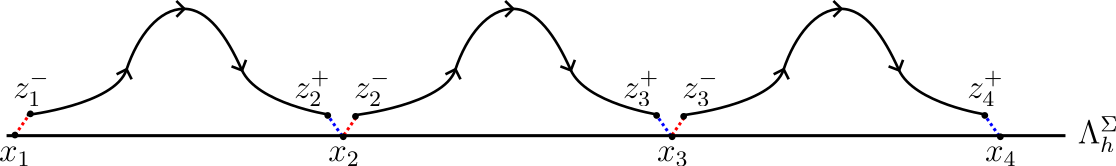}
\caption{Schematic figure showing a guiding orbit  $(\bx_i)_{i=1}^4\in\NHIMhred$ and the corresponding reduced pseudo-trajectory $(\bz^-_i, \bz^+_{i+1})_{i=1}^3$. Finite homoclinic trajectory segments corresponding to the reduced pseudo-trajectory are shown in black; unstable fibers in red, stable fibers in blue. See also Figure~\ref{fig:pseudo_orbit_refined} for a numerical example.}
\label{fig:pseudo_orbit}
\end{figure}

\begin{defn}
The sequence of pairs $(\bZ^-_i, \bZ^+_{i+1})_{i=1}^{N-1}$ will be called a \emph{pseudo-trajectory} of the system, while the sequence of pairs $(\bz^-_i, \bz^+_{i+1})_{i=1}^{N-1}$ will be called a \emph{reduced pseudo-trajectory}. 
\end{defn}
They satisfy the following conditions:
\begin{enumerate}
    \item For every $i=1,\dotsc,N-1$, the pair $(\bZ^-_i, \bZ^+_{i+1})$ are the endpoints in $\partial U$ of a (finite-time $\Text^i=T_-^i+T_+^{i+1}$) homoclinic trajectory, i.e. there exists a homoclinic point $\bzeta_i\in W^\un(\NHIM)\cap W^\st(\NHIM)$ and times $T_-^i, T_+^{i+1}$ such that
    \[
    \bZ^-_{i} = \Phi_{-T_-^i}(\bzeta_i) \qquad \textrm{and} \qquad \bZ^+_{i+1} = \Phi_{T_+^{i+1}}(\bzeta_i).
    \] 

    \item For every $i=2,3,\dotsc,N-1$, the forward trajectory of $\bz^+_i$ and the backward trajectory of $\bz^-_{i}$ synchronize to the same orbit: As $t\to\infty$, we have
    \[
    \norm{\Phi_t(\bz^+_i) - \Phi_t(\bx_i)}\to 0 \qquad \textrm{ and } \qquad
    \norm{\Phi_{-t}(\bz^-_i) - \Phi_{-t}(\bx_i)}\to 0.
    \]
    By Remark~\ref{rem:Zpm_to_tZpm}, the forward trajectory of $\bZ^+_i$ and the backward trajectory of $\bZ^-_{i}$ also tend to the same orbit: As $t\to\infty$, we have
    \[
    \norm{\Phi_{t+t_+}(\bZ^+_i) - \Phi_t(\bx_i)}\to 0 \qquad \textrm{ and } \qquad
    \norm{\Phi_{-(t+t_-)}(\bZ^-_i) - \Phi_{-t}(\bx_i)}\to 0.
    \] 
\end{enumerate}

\begin{rem} \label{rem:reduced_pseudo_orbit}
The sequence of points $(\bz^-_i)_{i=1}^{N-1}$ determines a finite \emph{pseudo-orbit} in the classical sense of dynamical systems~\cite{Bowen75,KPC07}: there is an associated sequence $(h_i)_{i=1}^{N-2}$ of positive times, $h_i = \Text^i+t_-^i+t_+^{i+1}$, such that
\[ \norm{\bz^-_{i+1} - \flow^{h_i}(\bz^-_i)} \leq \delta \quad \text{for }i=1,\dotsc,N-2,\]
where we are using the maximum norm.
Notice that 
$\flow^{h_i}(\bz^-_i) = \bz^+_{i+1}$. Indeed, for all $i$, both points $\bz^-_i$ and $\bz^+_i$ are $\delta$-close to $\bx_i$:
\begin{align*}
    \norm{\bz^+_i - \bx_i} = \abs{y} &= \rho_+:=\ee^{-t_+\,\lambda}\delta\leq \delta, \\
    \norm{\bz^+_i - \bx_i} = \abs{x} &= \rho_-:= \ee^{ -t_-\,\lambda}\delta\leq \delta.
\end{align*}
In our particular scenario, pseudo-orbits specifically consist of points on homoclinic orbits to $\NHIMhred$, not general orbits.
See Figure~\ref{fig:pseudo_orbit}. 
\end{rem}

\begin{defn}
A \emph{shadowing trajectory} is a trajectory of the system that passes close to the pseudo-trajectory, i.e., close to all the points $(\bZ^-_i, \bZ^+_{i+1})_{i=1}^{N-1}$.
\end{defn}

Given a pseudo-trajectory, our goal is to find a trajectory that shadows it.
In contrast to the pseudo-trajectory, the shadowing trajectory cannot lie exactly on the unstable/stable manifolds $W^{\un/\st}(\NHIM)$, otherwise it would approach a neighborhood of $\NHIM$ asymptotically and never leave it. However, the shadowing trajectory needs to be ``close'' to the unstable/stable manifolds if it has to pass close to the pseudo-trajectory. 
On the other hand, the shadowing trajectory should not be too close to $W^{\un/\st}(\NHIM)$ because then the transit time near $\NHIM$ will be too large and therefore inefficient. The choice of a suitable proximity to unstable/stable manifolds $W^{\un/\st}(\NHIM)$ is a key point to achieve adequate shadowing trajectories.
In practice, we want to use as much information from the pseudo-trajectory as possible to guide the search of the shadowing trajectory.

\section{Analytical shadowing}
\label{sec:analytical_shadowing}

In this section, we provide the proof of Theorem \ref{th:shadowing}.
We proceed in a constructive and direct way.

\subsection{Proof of Theorem \ref{th:shadowing}}

\begin{figure}[!ht]
\centering
\includegraphics[height=0.6\textheight]{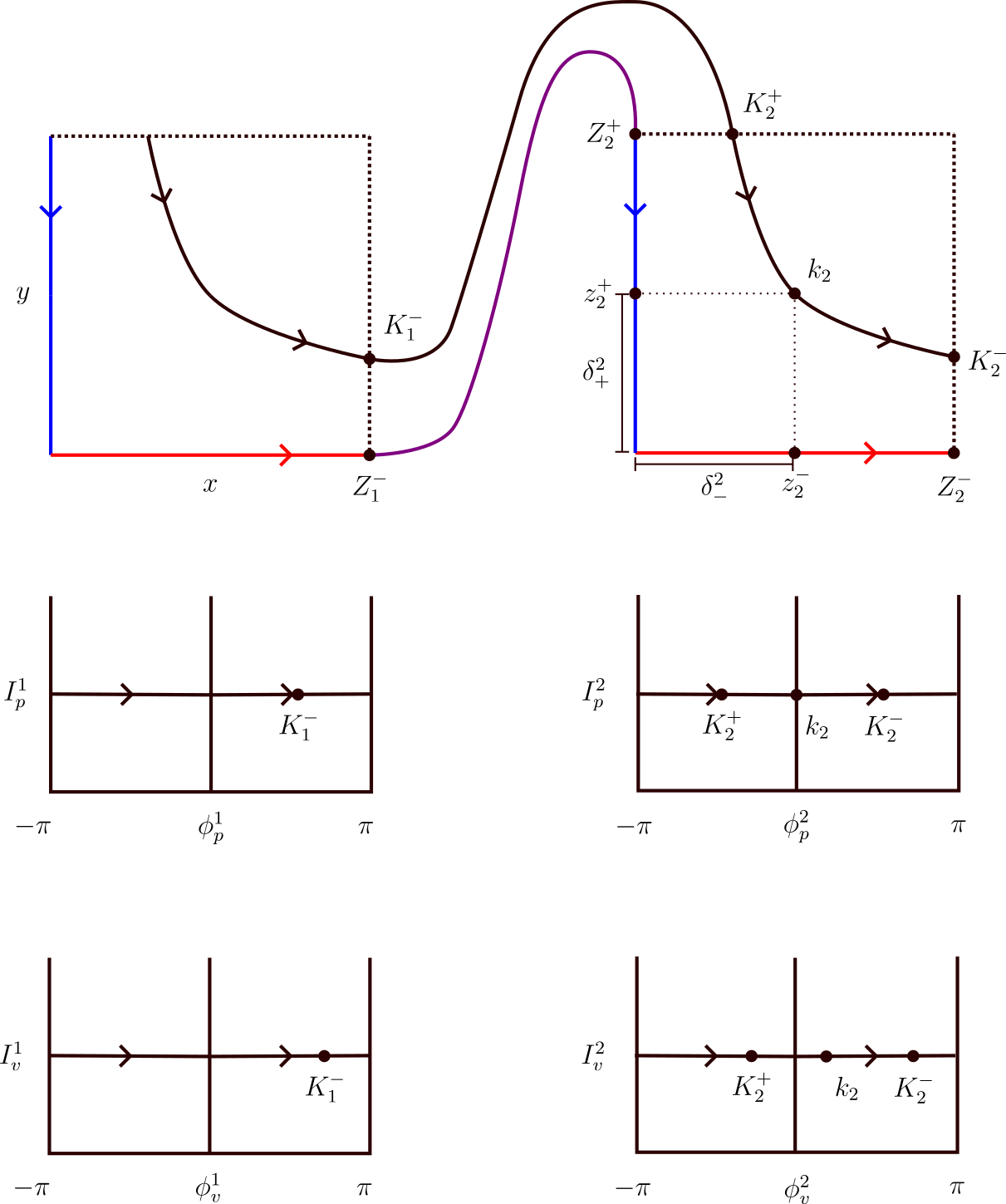}
\caption{First step of shadowing trajectory. \PauModify{The trajectory through $\bK^{+}_2$ crosses the section $\SecLocal$ only once at a point 
$\bk_2$.}}
\end{figure}

Let $\left(\bx_i\right)_{i=1}^N$ be a guiding orbit of the iterated function system (IFS) defined on the reduced NHIM $\NHIMhred$~\eqref{eq:NHIMred} generated by $L\geq 1$ transition maps $\{\tM_1,\ldots,\tM_L\}$, expressed in BNF coordinates~\eqref{eq:BNFcoordinates}. 
Let $(\bZ^-_i, \bZ^+_{i+1})_{i=1}^{N-1}$ be the pseudo-trajectory~\eqref{eq:Zi}, and $(\bz^-_i, \bz^+_{i+1})_{i=1}^{N-1}$ be the reduced pseudo-trajectory~\eqref{eq:tZi} associated with the guiding orbit.

Consider the first homoclinic segment $[\bZ^-_1,\bZ^+_2]$ with $\bZ^-_1\in\SecMinus$, $\bZ^+_2\in\SecPlus$, and  $\bZ^+_2=\Phi_{\Text^1}\left(\bZ^-_1\right)$.
Consider now the trajectory of any point
\[
\bK^{-}_1=\left(\delta,s_1/\delta,\bJ^1+\Delta \bJ^1,\bphi_{-}^1+\Delta\bphi_{-}^1\right)=\bZ^{-}_1+\Delta \bZ^{-}_1\in\SecMinus
\]
that is $\varepsilon_-^1$-close to $\bZ^{-}_1=\left(\delta, 0, \bJ^1,\bphi_{-}^1\right)$, i.e.,
$\Delta \bZ^{-}_1=\left(0,s_1/\delta,\Delta \bJ^1,\Delta\bphi_{-}^1\right)$ is small:
\[
\norm{\Delta \bZ^{-}_1}=\max\left\{s_1/\delta,\norm{\Delta \bJ^1},\norm{\Delta\bphi_{-}^1}\right\}=\varepsilon_-^1.
\]
Outside $U$, if $\varepsilon_-^1$ is small enough, this trajectory will cut the incoming boundary section $\SecPlus $ at a point 
\[
\bK^{+}_2=\Pext(\bK^{-}_1)=\left(s_2/\delta,\delta,\hat{\bJ}^2=\bJ^2+\Delta \bJ^2,\hat{\bphi}_{+}^2=\bphi_{+}^2+\Delta\bphi_{+}^2\right)=\bZ^{+}_2+\Delta \bZ^{+}_2\in \SecPlus,
\]
where $\Pext: \SecMinus\to\SecPlus$ is the Poincar\'e map of the Hamiltonian flow associated to the homoclinic segment
$[\bZ^-_1,\bZ_+^2]$.
$\bK^{+}_2$ is $\varepsilon_+^2$-close to $\bZ_+^2=\left(0, \delta, \bJ^2, \bphi_{+}^2\right)$,  i.e., 
$\Delta \bZ^{+}_2=\left(s_2/\delta,0,\Delta \bJ^2,\Delta\bphi_{+}^2\right)$ satisfies
\[
\norm{\Delta \bZ^{+}_2}=\max\left\{s_2/\delta,\norm{\Delta \bJ^2},\norm{\Delta\bphi_{+}^2}\right\}=\varepsilon_+^2,
\]
see Figure~\ref{sec:analytical_shadowing}. 

Denoting the uniform log-scale bound $\LE$ as
\begin{equation*}
\ee^{\LE} \coloneq \max\left\{ \max_{\bZ_- \in \SecMinus} \norm{D\Pext(\bZ_-)},\max_{\bZ_+ \in \SecPlus} \norm{D\Pext^{-1}(\bZ_+)}\right\},
\end{equation*} 
it follows from the Mean Value Inequality that the variations $\Delta\bZ^-_1$, $\Delta\bZ^+_2$ satisfy the uniform bi-Lipschitz bound
\begin{equation}
\label{eq:LyapunovBound}
\ee^{-\LE} \norm{\Delta\bZ^-_1} \leq \norm{\Delta\bZ^+_2} \leq \ee^{\LE} \norm{\Delta\bZ^-_1},
\end{equation}
so $\varepsilon_+ \leq \ee^{\LE}\varepsilon_-$.

Continuing inside $U$, the trajectory of $\bK^{+}_2$, with $0<s_2<\delta^2$, is given by 
\[
\Phi_{t}\left(\bK^{+}_2\right)
=\left(\ee^{t\lambda}s_2/\delta, \ee^{-t\lambda}\delta, \hat{\bJ}^2,\hat{\bphi}_{+}^2+t\,\bnu\right).
\]
where $\lambda$ and $\bnu$ are evaluated at $\left(s_2, \hat{\bJ^2}\right)$ and are given in~(\ref{eq:lambda}--\ref{eq:nu}),   

We now extend $\Phi_{t}\left(\bK^{+}_2\right)$ until it reaches a point 
\[
\bK^{-}_2=\left(\delta,s_2/\delta,\hat{\bJ}^2,
\hat{\bphi}_{+}^2+\Tint^2\,\bnu\right)
=\bZ^{-}_2+\Delta \bZ^{-}_2\in\SecMinus,
\]
with $\bZ^{-}_2=\left(\delta,0,\bJ^2,\bphi_-^2\right)\in\SecMinus$ and 
$\Delta \bZ^{-}_2=\left(0,s_2/\delta,\Delta\bJ^2,\left(\bphi_{+}^2 +\Tint^2\,\bnu - \bphi_-^2\right)+\Delta\bphi_{+}^2\right)$,
where we have used the dynamical invariance of $I_\hr=xy=s_2$ and $\bJ$, and
the \emph{internal transit time}~\eqref{eq:Tint}
\[
\Tint^2=\Tint^2(\delta,s_2,\hat{\bJ}^2)=\frac{\ln(\delta^2/s_2)}{\lambda\left(s_2,\hat{\bJ}^2\right)}
=\dfrac{1}{\lambdaL}\ln\dfrac{\delta^2}{s_2} \left(1+O(s_2,\hat{\bJ}^2)\right).
\]
We can exploit the dependence of $\Tint^2$ with respect to the parameter $s_2$ 
to ensure that $\bphi_{+}^2 +\Tint^2\,\bnu - \bphi_-^2$ is small. Equivalently, this ensures that $\bphi_+^2 +\Tint^2\,\bnu$ remains close to the angular component $\bphi_-^2 \in \T^2$ of the subsequent point 
$\bZ^{-}_2 = \left(\delta,0,\bJ^2,\bphi_-^2\right) \in \SecMinus$, 
which belongs to the next homoclinic segment $\left[\bZ^{-}_2,\bZ^{+}_3\right]$. 

Notice that $\bZ_2^+$ and $\bZ_2^-$ are linked because their saddle components permute from $(0,\delta)$ to $(\delta,0)$. Meanwhile, their center 
components follow a BNF trajectory from $(\bJ^2, \bphi^2_+)$ to 
$\left(\bJ^2, \bphi^2_-\right)$, satisfying $\bphi^2_- = \bphi^2_+ +\hTint^2 \bnu(0,\bJ^2)$.
When the trajectory through $\bK^{+}_2$ crosses the section $\SecLocal$ \emph{only once} at a point 
$\bk_2=\left(\delta^2_+,\delta^2_-,I^2,\phi^2\right)$, with $\delta^2_+\,\delta^2_-=s_2$ and $\phi^2_\pr=0$, that is, it is contained in the one-cut zone~\eqref{eq:OneCut}, it spends a time
$t^2_+$ during the transition $\SecPlus\to\SecLocal$ (the primary Poincar\'e map), and a time $t^2_-$ during $\SecLocal\to\SecMinus$.
Consequently, the total internal transition time $\hTint^2$
is given by 
\[
\hTint^2 = t^2_+ + t^2_- = \dfrac{1}{\lambda(0,\bJ^2)}\ln \dfrac{\delta^2}{s_2}=\dfrac{1}{\lambdaL}\ln\dfrac{\delta^2}{s_2} \left(1+O(I^2)\right).
\]
Inside the one-cut zone~\eqref{eq:OneCut}, the angular difference simplifies to:
\begin{align*}
\bphi_{+}^2 +\Tint^2\,\bnu - \bphi_-^2&=\bphi_{+}^2 +\hTint^2\,\bnu(0,I^2) - \bphi_-^2+\Tint^2\,\bnu-\hTint^2\,\bnu(0,I^2)\\
&=\Tint^2\,\bnu\left(s_2,\hat{\bJ}^2\right)-\hTint^2\,\bnu(0,I^2)\\
&=\ln\frac{\delta^2}{s_2}\left(\frac{\bnu\left(s_2,\hat{\bJ}^2\right)}{\lambda\left(s_2,\hat{\bJ}^2\right)}-
\frac{\bnu(0,I^2)}{\lambda(0,\bJ^2)}\right)=\ln\frac{\delta^2}{s_2}\, O\left(s_2,\Delta I^2\right),
\end{align*}
which yields
\[
\Delta \bZ^{-}_2=\left(0,s_2/\delta,\Delta\bJ^2,\Delta\bphi_{+}^2+\ln\frac{\delta^2}{s_2}\, O\left(s_2,\Delta I^2\right)\right).
\]
Inside the one-cut zone~\eqref{eq:OneCut} we have that $ \displaystyle\ln \frac{\delta^2}{s}<\frac{4\pi\lambda}{\nup }$. Thus, there exists a constant $c>0$ such that
\[
\norm{\Delta \bZ^{-}_2}\leq c \norm{\Delta \bZ^{+}_2},
\]
which finally leads to
\[
\norm{\Delta \bZ^{-}_2}\leq C \norm{\Delta \bZ^{-}_1}
\]
where $C=\ee^{\LE} c$.

Iterating this process along the guiding orbit, we finally obtain
\[
\norm{\bK^{-}_i-\bZ^{-}_i}=\norm{\Delta \bZ^{-}_i}\leq C^{i-1} \norm{\Delta \bZ^{-}_1},\quad i=2,\dots, N.
\]

Notice that we only need for the recursion that successive $\bK_i^-\in\SecMinus$, so
choosing $\bK^{-}_1-\bZ^{-}_1=\Delta \bZ^{-}_1=\left(0,s_1/\delta,\Delta \bJ^1,\Delta\bphi_{-}^1\right)$ small enough such that 
$C^{N-1} \norm{\Delta \bZ^{-}_1}< \delta$, we obtain a true trajectory $\flow(t)$ of the system that satisfies 
$\norm{\flow(t_i)-\bZ_{i}^-}<\delta$ for some times $0 < t_1 < \dots < t_N$.

It remains to check that the trajectories through successive $\bK^{-}_i$ cross $\SecLocal$ exactly once, so that they are in the one-cut zone~\eqref{eq:OneCut}, that is,
\begin{equation}
\label{eq:OneCuti}
   \ee^{\displaystyle -2\pi\lambda/\nup} <\frac{s_i}{\delta^2}\ee^{\displaystyle \lambda t_{+,i}}\leq 1,
\end{equation}
with $0 < t_{+,i} =\displaystyle \displaystyle  \frac{2\pi-\phi^i_{+,\pr}}{\nup }<\frac{2\pi}{\nup}$ as in~\eqref{eq:tplusminus} and $\lambda=\lambda(s_i,I^i)$, $\nup=\nup(s_i,I^i)$ as in~(\ref{eq:lambda}-\ref{eq:nu}).

By the bounds~\eqref{eq:LyapunovBound} and~\eqref{eq:tplusminus} , we just need to have 

\[
   \ee^{\displaystyle -2\pi\lambda/\nup}
  <\ee^{-(i-1)\LE}\frac{s_1}{\delta^2}\ee^{\displaystyle -2\pi\lambda/\nup} 
   <\ee^{(i-1)\LE}\frac{s_1}{\delta^2}\ee^{\displaystyle 2\pi\lambda/\nup}\leq 1,
    \quad i=1,\dots,N,
\]
which can be written as  
\[
\ee^{(i-1)\LE}<
  \frac{s_1}{\delta^2},\qquad  \frac{s_1}{\delta^2}\ee^{\displaystyle 2\pi\lambda/\nup}\leq \ee^{-(i-1)\LE},
\]
and is implied by the inequalities that take place for $i=N$
\begin{equation}
\label{eq:s1admissible}
  \ee^{(N-1)\LE}<
  \frac{s_1}{\delta^2} <\ee^{-(N-1)\LE}\ee^{-\displaystyle 2\pi\lambda/\nup}.
\end{equation}
For the bounds~\eqref{eq:s1admissible} to be satisfied, we need a left hand side smaller than the right hand side, that is
\begin{equation}
\label{eq:HomecomingCondition}
(N-1)\LE<\pi\lambda/\nup,
\end{equation}
which will be called the  \emph{homecoming condition}. 

Under this condition, choosing any $s_1$ satisfying~\eqref{eq:s1admissible}, for instance $s_1=\delta^2\exp\left(-\pi\lambda/\nup \right)$ or close enough, the bounds~\eqref{eq:OneCuti} are satisfied,
and the sequence $\bK_i^\pm $ is well defined.

Furthermore, since each point $\bk_{i}=\left(\delta^i_-,\delta^i_+,I^ i,\phi^i\right)\in\SecLocal$, where $\delta^i_+\,\delta^i_-=s_i$ and $\phi^i_\pr=0$, belongs to the trajectory that joins $\bK^{+}_i$ to $\bK^{-}_i$, it is also clear that $\norm{\flow(t_i)-\bk_{i}}<\delta$ holds for \emph{any} time of the trajectory joining $\bK^{+}_i$ to $\bK^{-}_i$, in particular for the time associated to the point $\bk_{i}$. For the total drift trajectory we can choose $t_1=\tT_1$, and at the step $i$ the accumulated time can be thus chosen as $t_i=\tT_1+\cdots+\tT_i$.
\qed

\begin{rem}
Condition~\eqref{eq:HomecomingCondition} is sufficient but not necessary to guarantee that the sequence of points $\bK_i^\pm$ remains within $U$ with trajectories intersecting $\SecLocal$ exactly once. For example, in Section~\ref{sec:semianalytical_results} we obtain a numerical trajectory of the SR3BP system that $\delta$-shadows the guiding orbit, whose points $\bK_i^\pm$ remain within $U$ with trajectories intersecting $\SecLocal$ exactly once, even though~\eqref{eq:HomecomingCondition} is not satisfied. See Remark~\ref{rem:shadowing_thm}.
\end{rem}

\section{Topological shadowing}
\label{sec:shadowing_topological}
In this section we provide the proofs of Lemma \ref{lem:shadowing_topological} and Theorem \ref{th:shadowing_topological}.

\subsection{Proof of Lemma \ref{lem:shadowing_topological}}

\emph{Construction of a sequence of balls in $\Lambda$.}
Let $\tau_i$ be a transition map on $\Lambda$, for some corresponding $m_i,n_i>0$.
If   $B^c_r(\bx)$  is a ball of radius $r$ centered at   $ \bx\in\Lambda$, then we define the inner  and outer radii
\[\begin{split}
r^i_{\min}(\bx,r)=\sup\{\rho>0\,|\, B^c_\rho(\tau_i( \bx))\subseteq \tau_i(B^c_r( \bx))\}  \\
r^i_{\max}( \bx,r)=\inf\{\rho>0\,|\, \tau_i(B^c_r( \bx))\subseteq B^c_\rho(\tau_i( \bx))\}.
\end{split}\]
Due to the expansion $\tau_i( \bx+h)=\tau_i( \bx)+D\tau_i( \bx)h+O(h^2)$ we have that
\[ \begin{split}
r^i_{\min}( \bx,r)=s^i_{\min}(D\tau_i( \bx))r+O(r^2)\\
r^i_{\max}( \bx,r)=s^i_{\max}(D\tau_i( \bx))r+O(r^2)
\end{split}\]
where $s^i_{\min}(\tau_i( \bx))$ and $s^i_{\max}(\tau_i( \bx))$ are the smallest and the largest singular values of $D\tau_i( \bx)$.
Since $\tau_i$ is symplectic we have that  $s^i_{\min}\cdot s^i_{\max} =1$ and $r^i_{\min}\cdot r^i_{\max}=r^2+O(r^3)$.
Note that $r^i_{\min}$ and $r^i_{\max}$ depend on $ \bx$, $r$ and $\tau_i$, and implicitly on the number of iterates $M_i=n_i+m_i$.

\begin{figure}
\centering
\includegraphics[width=0.6\textwidth]{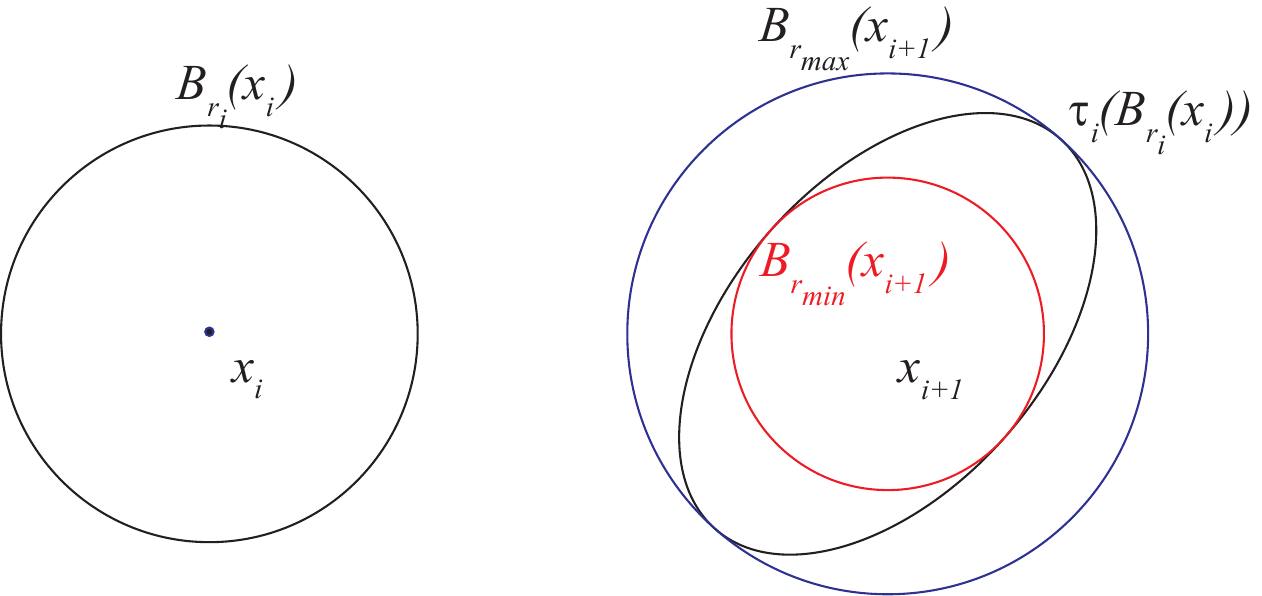}
\caption{Sequence of balls in $\Lambda$.}
\label{fig:balls_chain}.
\end{figure}

As $\Lambda$ is compact and the collection of transition maps is finite, we can find  two continuous functions $r_{\min}(r)$ and $r_{\max}(r)$, independent of $ \bx$ and $i$, but depending on $M_{\max}=\max\{M_1,\ldots,M_{N-1}\}$, such that, for any $i\in\{1,\ldots, N-1\}$, $ \bx\in \Lambda$ and $r>0$ we have
\begin{equation*}
B^c_{r_{\min}(r)}(\tau_i( \bx))\subseteq \tau_i(B^c_r( \bx))\subseteq B^c_{r_{\max}(r)}(\tau_i( \bx)).
\end{equation*}
See Fig.~\ref{fig:balls_chain}.

We now construct inductively a finite sequence of balls  in $\Lambda$ as follows:
\begin{itemize}
\item $B^c_{r_1}(\bx_1)$ for some $r_1>0$;
\item  $B^c_{r_{i+1}}( \bx_{i+1})$ where $r_{i+1}=r_{\max}(r_i)$, for $i=1,\ldots, N-1$.
\end{itemize}
By the construction, we have that
\begin{equation*}
\tau_i(B^c_{r_i}(\bx_i))\subseteq B^c_{r_{i+1}}(\bx_{i+1}), \textrm { for } i=1,\ldots, N-1.
\end{equation*}

By the continuity of the function $r_{\max}(r)$ and finiteness of the number of steps, there exists $r_1$ sufficiently small such that $r_1,\ldots, r_N$ are all less than $\varepsilon>0$. Note that $r_1$ depends on $M_{\max}$.

We remark here  that this construction does not involve the function $r_{\min}(r)$.

\medskip

\emph{Construction of stable and unstable disks at the homoclinic intersection.} 
Let  $\Gamma_{\alpha(i)}$ be a homoclinic channel associated to $\sigma_{\alpha(i)}$.

Choose a homoclinic point $\bk_i\in\Gamma_{\alpha(i)}$.
Then choose a pair of transverse, $1$-dimensional  disks $\mathscr{D}_{i+1}^s$ around $\bk_i$ in $W^s( \bx_{i+1})$ and $\mathscr{D}_{i}^u$ around $\bk_i$ in $W^u(\bx_{i})$, respectively.

By the Lambda Lemma, $\frm^{-n_i}(\mathscr{D}_{i+1}^s)$ approaches a disk $B^s_{\rho_L}(x_i)$ in $W^s(\bx_i)$, relative to the $C^1$ topology, as  $n_i\to\infty$.

Also by the Lambda Lemma, $\frm^{m_i}(\mathscr{D}_{i}^u)$ approaches a disk $B^u_{\rho_L}(x_{i+1})$ in $W^u( x_{i+1})$, relative to the $C^1$ topology, as  $m_i\to\infty$.

For  $m_i, n_i$ sufficiently large (to be made precise later), we denote
 $\bz_i^-$ the $x$-coordinate of $\frm^{-n_i}(\bk_i)\in  \frm^{-n_i}(\mathscr{D}_{i+1}^s)\cap W^u(\Lambda)$, and 
  $\bz_{i+1}^+$ the $y$-coordinate of $\frm^{m_i}(\bk_i)\in \frm^{m_i}(\mathscr{D}_{i}^u)\cap W^s(\Lambda)$.

\medskip

\emph{Construction of correctly aligned windows in the normal form coordinates.}

For the definition of correctly aligned windows and their properties see \cite{zgliczynski2004covering}. 

Define:
\begin{equation*}
  \begin{split}
     W_i^+ =&B^u_{\rho_L}(0)\times B^s_{\rho_S} (\bz_i^+)\times B^c_{r_i}(\bx_i),\\
    W_i^- =&B^u_{\rho_S}(\bz_i^-)\times B^s_{\rho_L} (0)\times  B^c_{r_i}(\bx_i)
  \end{split}
\end{equation*}
where the superscripts $u,s,c$ indicate that the respective  balls are in the unstable, stable and center coordinates $(x,y,\theta)$ described in \eqref{eqn:coordinates}.

We require $\rho_L>\rho_S>0$ to be so that
\begin{itemize}
\item $B^u_{\rho_L}(0)\subset \textrm{Pr}_x [\frm^{m_{i-1}}(\mathscr{D}_{i-1}^u)]$;
\item $\frm^{m_{i-1}}(\mathscr{D}_{i-1}^u)$ is a graph over $x$ that crosses transversally each disk $\{x\}\times B^s_{\rho_S} (\bz_i^+)\times B^c_{r_i}(\bx_i)$ in $W_i^+$, for every $x\in B^u_{\rho_L}(0)$; 
\item $B^s_{\rho_S}(0)\subset \textrm{Pr}_y
[\frm^{-n_i}(\mathscr{D}_{i+1}^s)]$
\item $\frm^{-n_i}(\mathscr{D}_{i+1}^s)$  is a graph over $y$ that crosses transversally each disk $B^u_{\rho_S}(\bz_i^-)\times \{y\}\times  B^c_{r_i}(\bx_i)$ in $W_{i}^-$, for every $y\in B^s_{\rho_L} (0)$.
\end{itemize}

In other words, $W_i^+$ is constructed by `thickening' $\frm^{m_{i-1}}(\mathscr{D}_{i-1}^u)$ in the center-stable directions, and $W_{i}^-$ is constructed by `thickening' $\frm^{-n_i}(\mathscr{D}_{i+1}^s)$ in the center-unstable directions.

We define the exit and entry sets of these windows as follows:
\begin{equation*}
  \begin{split}
  \left( W_i^+\right )^{\rm ex} =&\partial\left[ B^u_{\rho_L}(0) \right ] \times B^s_{\rho_S} (\bz_i^+)\times B^c_{r_i}(\bx_i),\\
     \left( W_i^+\right )^{\rm en} =&B^u_{\rho_L}(0) \times \partial\left[ B^s_{\rho_S} (\bz_i^+)\times B^c_{r_i}(\bx_i) \right ] \\
     \left( W_i^-\right )^{\rm ex}  =&\partial\left[ B^u_{\rho_S}(\bz_i^-) \right]\times B^s_{\rho_L} (0)\times  B^c_{r_i}(\bx_i) \\
     \left( W_i^-\right )^{\rm en}  =&B^u_{\rho_S}(\bz_i^-) \times \partial \left[ B^s_{\rho_L} (0)\times  B^c_{r_i}(\bx_i) \right].
  \end{split}
\end{equation*}
Note that the topologically stable direction of the windows comprises both the hyperbolic stable direction and the center directions.

There exist  $M_{\max}>M_{\min}$ and $\rho_L>\rho_S$ with  $\rho_L<\varepsilon$, such that,  if $M_i=m_i+n_i$ satisfies $M_{\min}\le M_i=m_i+n_i\le M_{max}$ we have:
\begin{itemize}
\item $B^u_{\rho_L}(0)\times B^s_{\rho_S} (\bz_i^+)$ is correctly aligned under the identity $\textrm{Id}_{x,y}$ with $B^u_{\rho_S}(\bz_i^-)\times B^s_{\rho_L} (0)$;
\item $B^u_{\rho_S}(\bz_i^-)\times B^s_{\rho_L} (0)$ is correctly aligned under $\textrm{Pr}_{x,y}\circ \frm^{m_i+n_i}(\cdot,\cdot,\bx_i)$ with  $B^u_{\rho_L}(0)\times B^s_{\rho_S} (\bz_{i+1}^+)$.
\end{itemize}

For  $m_i,n_i$ fixed as above, by the construction of the sequence of balls at the beginning of the proof, there exist $r_1,\ldots, r_N$ such that    $B_{r_i}(\bx_i)$ is correctly aligned with $B_{r_{i+1}}(\bx_{i+1})$ under the transition map $\tau_i$, in the sense that
\[\tau_i(B_{r_i}(\bx_i))\subset  B_{r_{i+1}}(\bx_{i+1}).\]

The correct alignment of the hyperbolic components and that of the center components of the windows do not automatically imply that correct alignment of the full windows. Therefore, we choose  $r_1$ (and implicitly $r_2, \ldots, r_N$) so that
\begin{itemize}\item $B^u_{\rho_S}(\bz_i^-)\times B^s_{\rho_L} (0)$ is correctly aligned under $\textrm{Pr}_{x,y}\circ \frm^{m_i+n_i}(\cdot,\cdot,\by)$ with  $B^u_{\rho_L}(0)\times B^s_{\rho_S}(\bz_{i+1}^+)$
for all $\bx\in B^c_{r_i}(\bx_i)$.
\end{itemize}

\begin{figure}
\centering
\includegraphics[width=0.45\textwidth]{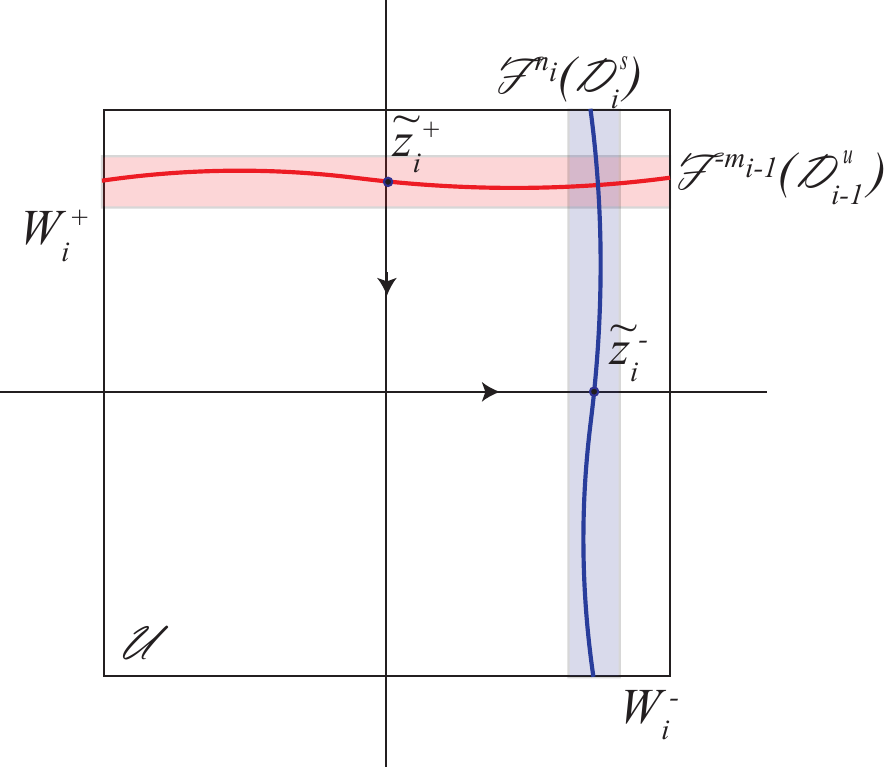}
\caption{Correctly aligned windows near $\Lambda$.}
\label{fig:windows_chain}
\end{figure}

Therefore, we obtain, for all $i=1,\ldots,N-1$,
 \begin{itemize}
\item $W_i^+$ is correctly aligned under the identity $\textrm{Id}$ with $W_i^-$;
\item $W_i^-$ is correctly aligned under $\frm^{n_i+m_i} $ with  $W_{i+1}^+$;
\end{itemize}
See Fig.~\ref{fig:windows_chain}.

In addition, we choose $r_1$ so that  all $r_1,\ldots, r_N$ are less than $\varepsilon>0$.

Then, by the property of correctly aligned windows\footnote{If $\{W_i\}_{i\in \mathbb{Z}}$ is a sequence of windows in $\mathbb{R}^n$ with unstable-like dimension $n_u$ and  stable-like dimension $n_s$, such that $n_u+n_s=n$, and $W_i$ is correctly aligned with $W_{i+1}$ under $\frm$ for all $i$, then there exists an orbit $\{\by_i\}_{i\in\mathbb{Z}}$ of $\frm$ such that $\by_i\in W_i$ and $\by_{i+1}=\frm(\by_i)$ for all $i$ \cite{zgliczynski2004covering}.}
there exists $\by_1$ whose orbit $\by_1,\ldots, \by_N$ satisfies
\begin{equation*}\begin{split}
\by_{i+1}= & \frm^{M_i}(\by_i), \textrm{ for } i=1,\ldots,N-1, \textrm { and }\\
\by_i\in & W^+_i\cap W^-_i,\, \textrm{ for } i=1,\ldots,N.
\end{split}
\end{equation*}

In particular, the orbit $\by_1,\ldots, \by_N$  $\varepsilon$-shadows $\bx_1,\ldots, \bx_N$, as desired.
\qed


\subsection{Proof of Theorem \ref{th:shadowing_topological}}
We consider the BNF system given by $\mathcal{Z}$. Let $\tilde{\Phi}_t$ denote the flow of $\mathcal{Z}$,  $\frm$ be the first return map to $\Sigma$, and $f=\frm_{\mid\Lambda}$. 

Let $\left(\bx_i\right)_{i=1}^N$ be a guiding orbit in $\Sigma$ of the iterated function system (IFS) generated by the transition maps 
$\{\tM_1,\ldots,\tM_{N-1}\}$, such that
\[\bx_{i+1}=\tau_i(\bx_i)\]
where 
\[\tau_i=f^{m_*}\circ \sigma_{\alpha(i)}\circ f^{n_*},\]
for $i=1,\ldots,N-1$. 
Let $\varepsilon>0$ be arbitrarily small. By Lemma \ref{lem:shadowing_topological} there exist $m_*,n_*$ sufficiently large such that, for the  guided orbit $\left(\bx_i\right)_{i=1}^N$ of the corresponding transition maps, there is   an orbit of the first return map $\frm$ to $\Sigma$ such that 
\[ d(\bx_i,\frm^{M_*}(\bx_i))<\varepsilon,\]
where $M_*=m_*+n_*$.

Consequently, there exists a trajectory  $\tilde{\Phi}_t(\bx)$ of the system given by $\mathcal{Z}$  that $\varepsilon$-shadows $\left(\bx_i\right)_{i=1}^N$, i.e., there exists an sequence of times $0<t_1<\ldots<t_N$ such that
\begin{equation}
\label{eqn:error_1}
d(\bx_i,\tilde{\Phi}_{t_i}(\bx))<\varepsilon,\textrm{ for } i=1,\ldots,N.
\end{equation}

Since the  system $\mathcal{Z}$ is $\errorBNF$-close to the original system $H$, by  Gronwall's inequality, there exists a constant $C_0$ such that
\begin{equation}
\label{eqn:error_2}d(\Phi_t(\bx),\tilde{\Phi}_t(\bx))<\errorBNF e^{C_0 t}, \textrm{ for } t>0.\end{equation}
 
Then, for $C=1+e^{C_0 t_N}$, from \eqref{eqn:error_1} and \eqref{eqn:error_2}, we obtain
\begin{equation*}
d(\bx_i,\Phi_{t_i}(\bx))<\varepsilon+\errorBNF C, \textrm{ for } i=1,\ldots,N .\end{equation*}

\section{Multiple shooting} \label{sec:multiple_shooting}

Given the pseudo-trajectory $(\bZ^-_i, \bZ^+_{i+1})_{i=1}^{N-1}$, we set up a \emph{multiple shooting} scheme to find a shadowing trajectory. 
The method divides the interval over which the trajectory is sought into several smaller intervals, solves an initial value problem in each of the smaller intervals, and imposes additional matching conditions to form a solution on the whole interval.

Concretely, we look for points $(\bK^-_i, \bK^+_{i+1})_{i=1}^{N-1}$ located \emph{off} the unstable/stable manifolds,
\begin{subequations} \label{eq:zeta}
\begin{align}
    \bK_i^- &= \bZ_i^- + \Delta \bZ_i^- = \left(\delta,y,I,\phi\right) \in \SecMinus, \label{eq:zeta_minus}\\
    \bK_i^+ &= \bZ_i^+ + \Delta \bZ_i^+ = \left(x,\delta,I,\phi\right) \in \SecPlus, \label{eq:zeta_plus}
\end{align}
\end{subequations}
that satisfy the \emph{matching conditions}
\begin{enumerate}
    \item[C.1] 
    For every $i=1,\dotsc,N-1$, there exist an intermediate point $\bzeta$ in the global section $\SecGlobal =\{Y=0\}$ and an exterior time $\Text = T_- + T_+$ such that
    \begin{equation} \label{eq:outer} 
    \bK^-_{i} = \Phi_{-T_-}(\bzeta) \qquad \textrm{and} \qquad \bK^+_{i+1} = \Phi_{T_+}(\bzeta). 
    \end{equation}

    \item[C.2] For every $i=2,3,\dotsc,N-1$, the forward trajectory of $\bK^+_i\in\SecPlus$ matches the backward trajectory of $\bK^-_{i}\in\SecMinus$, i.e. there exists an interior time $\Tint$ within $U$ such that 
    \begin{equation} \label{eq:inner}
    \Phi_{\Tint }(\bK^+_i) = \bK^-_i. 
    \end{equation}
\end{enumerate}
Notice that the first integral $I_\hr=xy$ of the BNF is $x\delta$ in $\bK^+_i$ and $\delta y$ in $\bK^-_i$. Thus, if Condition 2 is satisfied, then we must have $x=y$ in \eqref{eq:zeta}.

\begin{figure}[!ht]
\centering
\includegraphics[width=\textwidth]{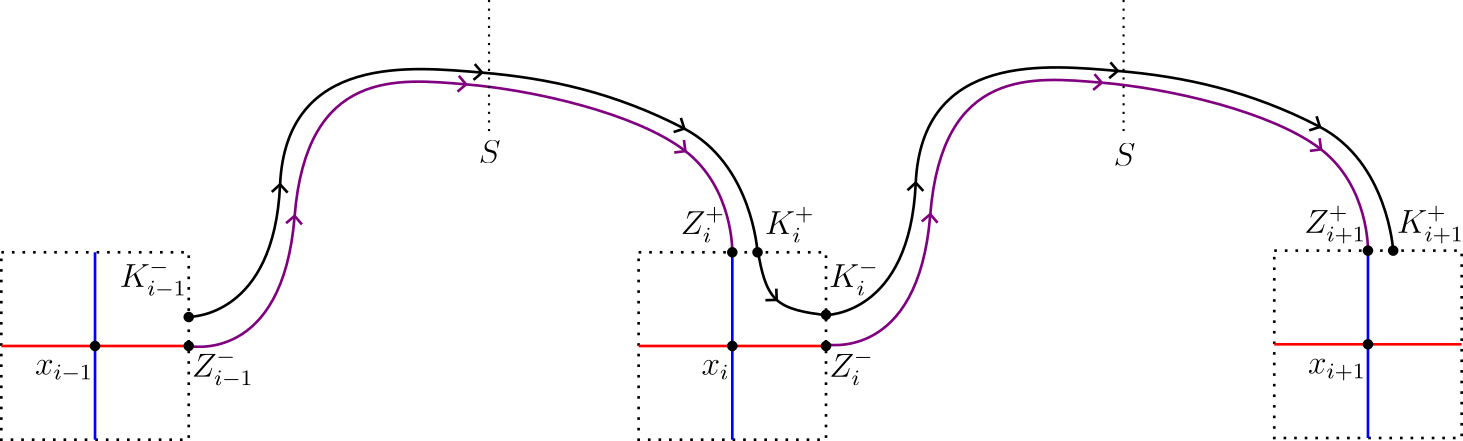}
\caption{Schematic figure showing the pseudo-trajectory $(\bZ^-_i, \bZ^+_{i+1})_{i}$, the corresponding chain of homoclinic trajectories (in purple), and a shadowing trajectory (in black). Condition 1 is verified on the Poincar\'e global section $\SecGlobal =\{Y=0\}$. 
See also Figure~\ref{fig:shadowing_trajectory} for a numerical example of shadowing trajectory.}
\label{fig:shadowing}
\end{figure}

Together, these two conditions warrant $\Phi_{\Text}(\bK^-_i)=\bK^+_{i+1}$ and $\Phi_{\Tint }(\bK^+_{i+1}) = \bK^-_{i+1}$, so every point in $(\bK^-_i, \bK^+_{i+1})_{i=1}^{N-1}$ belongs to one and the same trajectory. If the corrections $\Delta \bZ^\pm_i$ turn out to be small, then the points $\bK_\pm^i$ in the shadowing trajectory will be close to $\bZ^\pm_i$ in the pseudo-trajectory.
See Figure~\ref{fig:shadowing}.

\begin{rem}
    Figure~\ref{fig:shadowing} only shows a projection of the trajectories onto the saddle variables $x,y$. In fact, the action and angle variables $I,\phi$ must also match at the points $(\bK^-_i, \bK^+_{i+1})_{i=1}^{N-1}$.
\end{rem}

Typically, a modification of Newton's method is used to solve the system of equations~\eqref{eq:outer} and~\eqref{eq:inner}. 

    We look for the unknown points $(\bK^-_i, \bK^+_{i+1})_{i=1}^{N-1}$ in the boundary sections $\SecMinus , \SecPlus  $, respectively. Thus, they are at the boundary $\partial U$ of the normal form's domain (See Section~\ref{sec:BNF}). Numerically, this is very natural: Inside $U$, we use the integrable normal form to approximate the flow, thus avoiding the exponential growth of numerical errors due to hyperbolicity. (As shown in Section~\ref{sec:NumErrorSEBNF}, within $U$ the normal form is a very good approximation to the SR3BP flow.) 
    Outside $U$, we integrate the flow numerically.

To improve numerical stability, close to $\NHIM$ we always integrate the flow forward along the unstable manifold and backward along the stable manifold. Thus, it is natural to reformulate condition C.1 in the following way:
\begin{enumerate}    
\item[C.1'] For every $i=1,\dotsc,N-1$, there exists an intermediate point $\bzeta$ in the global section $\SecGlobal =\{Y=0\}$ and times $T_\pm>0$ such that
\begin{equation} \label{eq:outer2}
\Phi_{T_-}(\bK^-_{i}) =  \bzeta = \Phi_{-T_+}(\bK^+_{i+1}). 
\end{equation}
\end{enumerate}
This condition is checked in the intermediate Poincar\'e global section $\SecGlobal$, different from the boundary sections $\SecPlusMinus$.

We need to solve equations~\eqref{eq:inner} and~\eqref{eq:outer2} for $(\bK^-_i, \bK^+_{i+1})$. 
Notice that \eqref{eq:inner} consists of 5 equations in the boundary section $\SecMinus $, and \eqref{eq:outer2} consists of 5 equations in the global section $\SecGlobal$. Condition C.1' amounts to $5(N-1)$ equations, while Condition C.2 amounts to $5(N-2)$ more. Thus, the total number of equations is $10N-15$.
On the other hand, points $\bK^-_i$ consist of 5 variables in the boundary section $\SecMinus $, and $\bK^+_{i+1}$ consist of $5$ variables in the boundary section $\SecPlus  $. Thus, the total number of variables is $10N-10$.
Hence, we may choose up to 5 extra conditions that the shadowing orbit must satisfy. In our computations, we fix some of the variables of the first and last point in the shadowing trajectory, namely the following (marked in boldface):
\begin{equation}
 \begin{aligned} 
    \bK^-_1 &= \left(\delta,\boldsymbol{y=0},\boldsymbol{I_\pr=I_\pr^1}, \boldsymbol{I_\vr=I_\vr^1}, \bphi\right) &\in W^\un(\NHIMh) \cap \SecMinus, \\
    \bK^+_N &= \left(\boldsymbol{x=0},\delta,I_{\pr},\boldsymbol{I_{\vr}=I_{\vr}^N},\bphi\right) &\in W^\st(\NHIMh) \cap  \SecPlus.
\end{aligned}
\label{eq:extra_conditions}
\end{equation}
In this way, we impose that the shadowing trajectory tends to the first and last torus of the guiding orbit, with actions $I_\pr^1, I_\vr^1$ and $I_{\vr}^N$, as $t\to\pm\infty$ (see Equation~\eqref{eq:Xi}).

\subsection{Initial guess} \label{sec:initial_guess}

To solve the system of $10N-15$ equations given by~\eqref{eq:inner} and~\eqref{eq:outer2}, we need an initial guess of the unknown points $(\bK^-_i, \bK^+_{i+1})_{i=1}^{N-1}$ in the boundary sections $\SecMinus , \SecPlus  $. 
We cannot simply take $\bZ^\pm_i$ as the initial guess for $\bK^\pm_i$, because $\bZ^\pm_i$ is exactly on the invariant manifold. Although $\bZ^\pm_i$ satisfies condition C.1', it does not satisfy condition C.2.
As explained in Section~\ref{sec:pseudo_orbits_shadowing_trajectories}, the shadowing trajectory must be close to the invariant manifolds $W^{\un/\st}(\NHIMh)$, and the choice of an appropriate closeness is key. 
This means choosing the distance $y$ from the point $\bK^-_i$ to the unstable manifold $W^u(\NHIMh)$ in Equation~\eqref{eq:zeta_minus}, and the distance $x$ from the point $\bK^+_i$ to the stable manifold $W^s(\NHIMh)$ in Equation~\eqref{eq:zeta_plus}.

\begin{figure}
\centering
\includegraphics{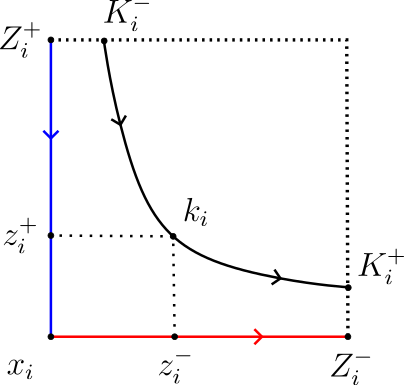}
\caption{Construction of an initial guess for the multiple shooting procedure. $(\bk_i)_i$ is an initial guess for shadowing the reduced pseudo-trajectory $(\bz^-_i, \bz_+^{i+1})_i$, whereas $(\bK^-_i, \bK^+_{i+1})_i$ is an initial guess for shadowing the non-reduced pseudo-trajectory $(\bZ^-_i, \bZ^+_{i+1})_i$.}
\label{fig:initial_guess}
\end{figure}

To obtain an appropriate initial guess, we will use the reduced pseudo-trajectory $(\bz^-_i, \bz^+_{i+1})_i$. Recall from Remark~\ref{rem:reduced_pseudo_orbit} that $\norm{\bz^+_i - \bz^-_i} \leq \delta$. For every $i=2,3,\dotsc,N-1$, consider the point $\bk_i\in\SecLocal$ defined as
\[ 
\bk_i = \left(x=\delta^i_-,y=\delta^i_+,\bJ^i,\bphi^i\right), \text{ with } I^i_\pr=\Fh(s_i,I^i_\vr), \,\phi^i_\pr=0, 
\]
where $\delta^i_-,\delta^i_+,I^i,\phi^i$ are given in Equation~\eqref{eq:tZi}, and $s_i=\delta^i_-\,\delta^i_+$.
Notice that $\bk_i$ and $\bz^-_i$ differ only in the variable $y$, while $\bk_i$ and $\bz^+_i$ differ only in the variable $x$.
See Figure~\ref{fig:initial_guess}.

Notice that $\bk_i$ is $\delta$-close to both $\bz^+_i$ and $\bz^-_i$, and is located off the unstable/stable manifolds. Moreover, the trajectory of $\bk_i$ approaches the trajectory of $\bz^+_i$ in backward time and that of $\bz^-_i$ in forward time,
\begin{align*}
    \norm{\Phi_{-t}(\bk_i) - \Phi_{-t}(\bz^+_i)} \leq \exp \left( -t \lambda  \right)\delta^i_+ \quad \text{for } t\leq t^i_-, \\ 
    \norm{\Phi_{t}(\bk_i) - \Phi_{t}(\bz^-_i)} \leq \exp \left( -t \lambda  \right)\delta^i_- \quad \text{for } t\leq t^i_+.
\end{align*}
Thus, $(\bk_i)_i$ is a good initial guess for shadowing the reduced pseudo-trajectory.

Finally, we transport $\bk_i$ to the boundary sections $\SecMinus $ and $\SecPlus  $ using the flow, in order to obtain an initial guess for shadowing the non-reduced pseudo-trajectory.
\begin{alignat*}{2}
    \bK_i^- &= \Phi_{t_-^i}(\bk_i) &&= \left(x=\delta, y=s_i/\delta,I^i,\phi^i + t_-^i \bnu\right)\in \SecMinus , \\
    \bK_i^+ &= \Phi_{-t_+^i}(\bk_i) &&= \left(x=s_i/\delta,y=\delta,I^i,\phi^i - t_+^i \bnu\right)\in \SecPlus  ,
\end{alignat*}
where $s_i=\,\delta^i_-\,\delta^i_+$, and $\bnu$ is evaluated at $\bI=(s_i,I^i)$.

\begin{rem}
Since $s_i=\delta^i_-\delta^i_+\neq 0$, the time $t_\pm^i=t_\pm^i(s_i,I^i)$ and frequency $\bnu=\bnu(s_i,I^i)$ used here is different from the time  $t_\pm^i=t_\pm^i(0,I^i)$ and frequency $\bnu=\bnu(0,I^i)$ used for $s=0$ in \eqref{eq:Zi}: $t_\pm^i(s_i,I^i)=(1+O(s_i))t_\pm^i(0,I^i)$, so $\bK^-_i$ and $\bZ^-_i$ differ not only in coordinate $y$, but also slightly in $\bphi$.
\end{rem}

Therefore, the sequence $(\bK^-_i, \bK^+_{i+1})_i$ is an initial guess to shadow the  pseudo-trajectory $(\bZ^-_i, \bZ^+_{i+1})_i$. 
The initial guess constructed $(\bK^-_i, \bK^+_{i+1})_i$ \textbf{satisfies condition C.2} but does not exactly satisfy C.1'. However, in practice, this initial guess is close enough to the solution of the multiple shooting method, and we will see that the method converges (Section~\ref{sec:Shadowing_trajectory}).

\section{General error estimates for the  BNF}
\label{sec:GeneralErrorBNF}
After a translation and a symplectic linear change of variables
\begin{equation*}
(X,Y,Z,P_X,P_Y,P_Z) = L+ C(Q,P)= L+C(Q_1, Q_2,Q_3, P_1, P_2, P_3),
\end{equation*}
the Hamiltonian of a general saddle-center-center takes the form $H=H_2+O_3$, 
where the quadratic part $H_2$ is already in its real normal form, which is formula (14.8) in~\cite{DelshamsGR16}
\[
H_{2}=\lambdaL Q_{1} P_{1}+\frac{\nupL}{2}\left(Q_{2}^{2}+P_{2}^{2}\right)+\frac{\nuvL}{2}\left(Q_{3}^{2}+P_{3}^{2}\right).
\]

We can construct the BNF of order $M\geq 2$ provided that the so called \emph{small divisors} associated with the equilibrium point $L$ remain strictly positive for all $s=2,\dots,M$. These small divisors are defined as follows: 
\begin{equation}
\label{eq:alphas}
\alpha_{s}=\alpha_{s}(\mu)=
\min \abs{k_{\hr} \lambdaL+i\left(k_{\pr} \nupL+k_{\vr} \nuvL\right)}
=\min\sqrt{k_{\hr}^{2} \left(\lambdaL\right)^{2}+
\left(k_{\pr} \nupL+k_{\vr} \nuvL\right)^{2}}>0,
\end{equation}
where the minimum is taken over all $k_{\hr}, k_{\pr}, k_{\vr} \in \Z$ so that
$0<\left|k_{\hr}\right|+\left|k_{\pr}\right|+\left|k_{\vr}\right| \leq s$.
By definition, these small divisors form a non-increasing sequence $\alpha_{2} \geq \alpha_{3} \geq \cdots\geq \alpha_{M}$.

Then, there exists a symplectic analytic change $(Q,P)=\nftrans^{(M)}(q,p)$ such that
\[
\mathcal{H}^{(M)}\coloneq  H\circ \nftrans^{(M)}=H_{2}+\mathcal{Z}^{(M)}(I)+\mathcal{R}^{(M)}(q,p),
\]
with
\[
\mathcal{Z}^{(M)}(\bI)=\sum_{\substack{s \text { even } \\
4 \leq s \leq M}} \mathcal{Z}_{s}(\bI),\quad  \mathcal{R}^{(M)}(q, p)=\sum_{s \geq M+1} \mathcal{R}_{s}^{(M)}(q, p),
\]
where $\mathcal{Z}_{s} (\bI)$ is a homogeneous polynomial of degree $s/ 2$ in the ``action" variables $\bI=(I_\hr,I_\pr,I_\vr)$, where $I_{\hr}=x y$, with $x=q_{1}$, $y=p_{1}$, $I_{\pr}=\dfrac{q_{2}^{2}+p_{2}^{2}}{2}$, $I_{\vr}=\dfrac{q_{3}^{2}+p_{3}^{2}}{2}$, and $\mathcal{R}_{s}^{(M)} (q,p)$ is a homogeneous polynomial of degree $s$ in variables $(q,p)$.

In dealing with normal forms near a fixed point of a Hamiltonian
system, it is usual to consider complex variables $(X,Y)$ for the variables $(Q,P)$, and to introduce the complex  normal form of the quadratic part of $H$. The complex canonical coordinates $(X, Y)$ are 
defined by the linear change
\begin{equation*}
X_1=q_1,\quad Y_1=p_1,\qquad 
X_j=\frac{1}{\sqrt{2}}\left(q_j-\iu p_j\right), \quad
Y_j=-\frac{1}{\sqrt{2}}\left(q_j+\iu p_j\right),\qquad
j=2,3,    
\end{equation*}
so that the complex normal form of $H_2$ becomes
\begin{equation}
\label{eq:H2}
H_{2}=\lambdaL X_{1} Y_{1}+\iu\nupL X_{2} Y_{2}+\iu\nuvL X_{3}Y_{3}.
\end{equation}

The following proposition is adapted from \cite[Proposition~1]{DelshamsG96}, which exclusively addresses the totally elliptic case; nevertheless, the formal computations remain completely analogous. Alternative bounds for the normal form terms and the remainder can be found in \cite{CaraccioloL20, EfthyGC04, GuzzoCP26}. Furthermore, the sizes of all variables are uniform; consequently, the results of Proposition~\ref{prop:bounds} apply when $\delta=r$ within the neighborhood $U=U(\delta,r)$ defined in \eqref{eq:U} where the BNF is established.

\begin{prop}\label{prop:bounds}
   Let $a=20$ and $b=548\cdot 3$. Let $H(X, Y)=\sum_{s \geqslant 2} H_{s}(X, Y)$ with $H_{2}$ given in~\eqref{eq:H2} and $\left\|H_{s}\right\| \leq c^{s-2}d$  for $s \geq 2$. Assume that the lower bounds for the resonances $\alpha_s$ given in~\eqref{eq:alphas} satisfy $\alpha_s>0$ for $2\leq s \leq M$.  Then
   
\noindent a) $\displaystyle \left\|\mathcal{Z}_{s}\right\|\leq\frac{\frac{1}{6}(6cd)^{s-2}(s-2)!}{\alpha_{3} \cdots \alpha_{s-1}}$
for $s$ even, $4\leq  s \leq M$,

\noindent b) $\displaystyle \left\|\mathcal{R}_{s}^{(M)}\right\|\leq\frac{a d^2(a c d)^{s-2}(M-3)!(M-2)^{s-M+2}}{\,\alpha_{3} \cdots \alpha_{M-1} \alpha_{M}^{s-M+2}}$ for $s\geq M+1$,

\noindent c) Defining $\displaystyle r_{M}^{*}=\frac{\alpha_{M}}{b c d M}$,
$\nftrans^{(M)}$ is analytic on $\hat{\mathcal{B}}_{r_{R}^{*}}$ and, for any $r \leq r_{R}^{*}$, one has the inclusion
$\nftrans^{(M)}\left(\mathcal{B}_{r}\right)\supset \mathcal{B}_{r / 2}$.

\end{prop}

It remains to specify the norms and neighborhoods used, which are:
\[
\begin{aligned}
\left\|H_{s}\right\|&\coloneq  \sum_{|\beta+\beta^{\prime}|_{1}=s}\left| H_{\beta, \beta^{\prime}}\right| \text{ for }
H_{s}(X, Y)=\sum_{\left|\beta+\beta^{\prime}\right|_{1}=s} H_{\beta, \beta^{\prime}} X^{\beta} Y^{\beta'},\\
 \mathcal{B}_{r}&\coloneq  \left\{(q, p) \in \mathbb{R}^{6}:\left|(q, p)\right|\leq r\right\},
 \text{ with } |(q, p)|\coloneq  \max _{\ell=1,2,3} \sqrt{q_{\ell}^{2}+p_{\ell}^{2}},\\
 \mathcal{\hat{B}}_{r}&\coloneq=\left\{(X, Y) \in \mathbb{C}^{6}:\left|(X, Y)\right|\leq r,(X_1,Y_1)\in\R^2,\bar{Y}_j=\iu X_j, j=2,3\right\},\\
 &\qquad\text{ with }  \left|(X, Y)\right|=\max _{\ell=1,2,3}
  \sqrt{\left|X_{\ell}\right|^{2}+\left| Y_{\ell}\right|^{2}}, \\
  |f|_{r}&\coloneq  \max_{(X,Y)\in\hat{\mathcal{B}}_{r} } |f(X, Y)| \text{ for a function } f(X,Y).
\end{aligned}
\]

Notice that given $f=\sum_{s} f_{s}$, then $|f|_{r} \leq \sum_{s}\left\|f_{s}\right\|r^{s}$, and that if there exists $A, B$ such that $\left\|f_{s}\right\| \leq A^s B$ for all $s$, then $f$ is analytic on $\hat{\mathcal{B}}_{r}$ for $r<1/A$. Therefore
\[
\left|\mathcal{R}^{(M)}\right|_r\leq\sum_{s\geq M+1} \left\|\mathcal{R}^{(M)}_s\right\| r^s=
\frac{(M-3)!}{a c^2 (M-2)^{M-2}} \frac{\alpha_M^{M-2}}{\alpha_{3} \cdots \alpha_{M-1}}
\sum_{s\geq M+1} \left(\frac{a c d (M-2)}{\alpha_M}r\right)^s,
\]
which is a convergent series for $\displaystyle r<\frac{\alpha_{M}}{a c d(M-2)}$, in which  case we have
\begin{equation}
\label{eq:errorBNF}
\left|\mathcal{R}^{(M)}\right|_r\leq\errorBNF:= \frac{(M-3)!}{a c^2 (M-2)^{M-2}}\cdot \frac{\alpha_M^{M-2}}{\alpha_{3} \cdots \alpha_{M-1}}
\cdot\frac{\left(\dfrac{a c d (M-2)}{\alpha_M}r\right)^{M+1}}{1-\dfrac{a c d (M-2)}{\alpha_M}r}.
\end{equation}

\noindent\textbf{Remark 1} One can use Stirling's formula
$\displaystyle n!=\sqrt{2 \pi n}\left(\frac{n}{e}\right)^{n}\left(1+O\left(\frac{1}{n}\right)\right)$
 to get 
 \[
 \frac{(M-3)!}{(M-2)^{M-2}}=\frac{(M-3)}{\mathrm{e}^{M-2}}\left(1+O\left(\frac{1}{M-2}\right)\right).
 \]

\noindent\textbf{Remark 2} $\displaystyle\frac{\alpha_{M}^{M-2}}{\alpha_{3} \cdots \alpha_{M-1}}=
\frac{\alpha_{M}}{\alpha_{3}} \cdots \frac{\alpha_{M}}{\alpha_{M-1}}\leq 1$.

\section{Analytical error estimates for the Sun–Earth BNF}
\label{sec:AnalErrorSEBNF}

The Hamiltonian of the SR3BP (Spatial Restricted 3 Body Problem) is
\[
H=\frac{1}{2}\left(P_{X}^{2}+P_{Y}^{2}+P_{Z}^{2}\right)+Y P_{X}-X P_{Y}-\frac{1-\mu}{r_{1}}-\frac{\mu}{r_{2}},
\]
Where $r_{1}^{2}=r_{\mathrm{S}}^{2}=(X-\mu)^{2}+Y^{2}+Z^{2}, r_{2}^{2}=r_{\mathrm{E}}^{2}=(X-\mu+1)^{2}+Y^{2}+Z^{2}$

From~\cite{DelshamsGR16}, the distance $\gamma$ from $L_{1}$ at $(X,Y,Z)=(\mu-1+\gamma,0,0)$ to the Earth at $(X, Y, Z)=(\mu-1, 0,0)$
satisfies
\[
\gamma^{5}-(3-\mu) \gamma^{4}+(3-2 \mu) \gamma^{3}-\mu \gamma^{2}+2 \mu \gamma-\mu=0,
\]
or equivalently
\[
\gamma^{3}=\frac{\mu(1-\gamma)^{2}}{3-2 \mu-\gamma(3-\mu-\gamma)}=:f(\mu, \gamma).
\]
Iterating $\gamma_{0}^{3}=\dfrac{\mu}{3(1-\mu)}$, $\gamma_{n+1}^{3}=f\left(\mu, \gamma_{n}\right)$, or Taylor expanding in powers of $\left(\dfrac{\mu}{3}\right)^{1 / 3}$, one arrives at \cite[formula (28), p. 137]{Szebehely67}
\[
\gamma=\left(\frac{\mu}{3}\right)^{1 / 3}\left(1-\frac{1}{3}\left(\frac{\mu}{3}\right)^{1 / 3}-\frac{1}{9}\left(\frac{\mu}{3}\right)^{2 / 3}+\cdots\right),\text{ or }\gamma^{3}=\frac{\mu}{3}\left(1+O\left(\frac{\mu}{3}\right)^{1 / 3}\right).
\]
Scaling coordinates~\cite[formula (4.1)]{DelshamsGR16},
\[
X=-\gamma x+\mu-1+\gamma, Y=-\gamma y, Z=\gamma z,
\]
and expanding in Legendre polynomials, we get the Hamiltonian (14.5) of~\cite{DelshamsGR16}
\begin{equation}\label{eq:H_original_coords}
H=\frac{1}{2}\left(p_{x}^{2}+p_{y}^{2}+p_{z}^{2}\right) +y p_{x}-x p_{y}-\sum_{n\geq 2} c_{m}(\mu) \varepsilon^{n} P_{n}\left(\frac{x}{\varepsilon}\right),
\end{equation}
where $\displaystyle c_{n}(\mu)=\frac{1}{\gamma^{3}}\left(\mu+(-1)^{n} \frac{(1-\mu) \gamma^{n+1}}{(1-\gamma)^{n+1}}\right)$, $\varepsilon^{2}=x^{2}+y^{2}+z^{2}$,
and the expansion above is valid for $\varepsilon<1$, where $1$ is the distance from $L_{1}$ at $(x, y, z)=(0,0,0)$ to the Earth at $(x, y, z)=(1,0,0)$.

In particular, the quadratic part of $H$ at the origin is~\cite[formula (14.6)]{DelshamsGR16}
\[
H_{2}=\frac{1}{2} p_{x}^{2}+\frac{1}{2} p_{y}^{2}+y _{x}-x p_{y}-c_{2} x^{2}+\frac{c_{2}}{2} y^{2}+\frac{1}{2} p_{z}^{2}+\frac{c_{2}}{2} z^{2}, \]
where
\[
c_{2}=\frac{1}{\gamma^{3}}\left(\mu+\frac{(1-\mu)^{2}}{(1-\gamma)^{3}}\right)=\frac{\mu}{\gamma^{3}}+\frac{1-\mu}{(1-\gamma)^{3}}
=4\left(1+O\left(\frac{\mu}{3}\right)^{1/3}\right).
\]
A symplectic linear change of variables
\begin{equation}\label{eq:linear_change_variables}
(x,y,z,p_x,p_y,p_z) = C(Q,P)= C(Q_1, Q_2,Q_3, P_1, P_2, P_3),
\end{equation}
casts $H_2$ into its real normal form
\[
H_{2}=\lambdaL Q_{1} P_{1}+\frac{\nupL}{2}\left(Q_{2}^{2}+P_{2}^{2}\right)+\frac{\nuvL}{2}\left(Q_{3}^{2}+P_{3}^{2}\right),
\]
where (see, for instance,~\cite{Jorba99})
\[
\lambdaL=\frac{\sqrt{c_{2}-2+\sqrt{9 c_{2}^{2}-8 c_{2}}}}{2}, \quad \nupL=\sqrt{\frac{9 c_{2}^{2}-8 c_{2}+2-c_{2}}{2}}, \quad \nuvL=\sqrt{c_{2}}.
\]

As the Hamiltonian $H$~\eqref{eq:H_original_coords} is analytic for $\varepsilon<1$, one can expect  $\norm{H_{s}}\leq c^{s-2} d$, $s \geq 2$ only for $c \geq 1$ in the original coordinates. 
Taking into account the linear symplectic change~\eqref{eq:linear_change_variables}, in the coordinates $(x,y)$ we obtain
$\norm{H_s}\leq c^{s-2} d$, with $c=2.09$ and $d=\norm{H_2}=\lambdaL+\nupL+\nuvL=6.63$. 
Indeed, we have numerically computed the norm of all homogeneous polynomials $H_s$ up to order 32, see Figure~\ref{fig:Hs_real}.
The value $c=2.09$ is obtained by fitting the function $c^{s-2}d$ to these numerical data.
Other theoretical estimates can be found in~\cite{GiorgilliDFGS89}.

For our numerical computations, we have implemented a BNF of order \(M=16\). This choice requires computing the small divisors $\alpha_{3} \geq \alpha_{4} \geq \cdots \geq\alpha_{16}$ given in equation~\eqref{eq:alphas} for the Sun-Earth mass ratio $\mu_{\text{Sun-Earth}} = 3.040423398444176 \cdot 10^{-6}$. In this system, $c_2$ takes the value $c_2 = 4.061074016255354$,  and $\lambdaL =2.532659174052968$, $\nupL =2.086453564223108$, $\nuvL =2.015210662996639$.

The frequency ratio $\frL \coloneq \nuvL / \nupL = 0.9658545474253123$ can be expressed 
by the continued fraction expansion:
\[
\frL = 1 + \frac{1}{-29 + \frac{1}{-3 + \frac{1}{-2}}}
\]
Consequently, the resonance order of $\frL$ is approximately $\operatorname{round}\left(\frac{1}{1-\frL}\right) = 29$. 
This implies that $28/29$ is a remarkably close rational approximation of $\frL$, 
with a residual of $\frL - 28/29 = 0.00033730604$, and there is no better rational approximation with a denominator smaller than $29$.

Evaluating the                                                                                                                            small divisors~\eqref{eq:alphas} reveals that:
\begin{align*}
\alpha(3) = \dots = \alpha(26) &= 0.21372870367940491, \\
\alpha(27) &= 0.16289523110846341, \\
\alpha(28) &= 0.091652329881995109, \\
\alpha(29) &= 0.020409428655526807.
\end{align*}
Thus, the divisors $\alpha(s)$ remain stable at approximately $0.21$ up to the order $26$, 
where the Birkhoff normal form (BNF) remains highly reliable. However, at order $29$, 
they drop significantly to $0.02$. This drop indicates that a less restrictive 
normal form, such as the Birkhoff--Gustavson Normal Form (BGNF), would provide greater 
reliability and better capture resonant effects for larger action values, including 
those corresponding to halo orbits.

Finally, plugging the radius $r=r_{16} = 2.8303 \cdot 10^{-5}$ 
into Equation~\eqref{eq:errorBNF} yields an error 
$\abs{\mathcal{R}^{(16)}}_{r}\leq \epsilon_{16} = 3.57 \cdot 10^{-7}$ 
for the BNF remainder. Although this error is reasonably small, the underlying theoretical bounds 
demand an extremely restrictive radius of applicability. In the next section, we 
compare this analytical prediction with direct numerical estimates of the actual 
validity radius of the BNF.

\section{Numerical error estimates for the Sun–Earth BNF}
\label{sec:NumErrorSEBNF}

The Proposition~\ref{prop:bounds} gives bounds on the norm of the normal form and its reminder. In this section, we consider the SR3BP and numerically compute its BNF up to order $M=16$. In this way, we can explicitly compute the norms $\norm{\mathcal{Z}_s}, \norm{\mathcal{R}_s}$ and compare them with their theoretical bounds.

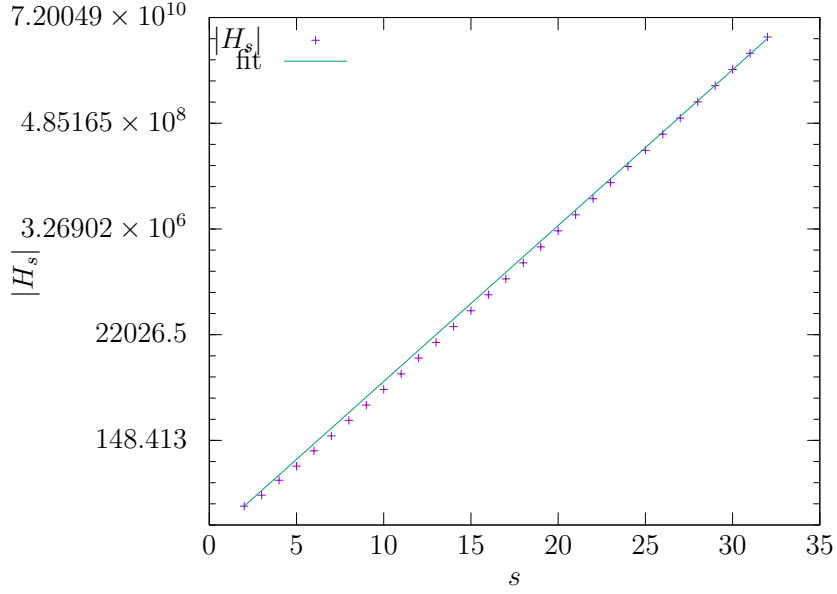
\begin{figure}[!ht]
\centering
	\scalebox{0.9}{\begin{tikzpicture}[gnuplot]
\path (0.000,0.000) rectangle (12.500,8.750);
\gpcolor{color=gp lt color border}
\gpsetlinetype{gp lt border}
\gpsetdashtype{gp dt solid}
\gpsetlinewidth{1.00}
\draw[gp path] (2.976,0.985)--(3.066,0.985);
\draw[gp path] (11.947,0.985)--(11.857,0.985);
\draw[gp path] (2.976,1.296)--(3.066,1.296);
\draw[gp path] (11.947,1.296)--(11.857,1.296);
\draw[gp path] (2.976,1.606)--(3.066,1.606);
\draw[gp path] (11.947,1.606)--(11.857,1.606);
\draw[gp path] (2.976,1.917)--(3.066,1.917);
\draw[gp path] (11.947,1.917)--(11.857,1.917);
\draw[gp path] (2.976,2.228)--(3.156,2.228);
\draw[gp path] (11.947,2.228)--(11.767,2.228);
\node[gp node right] at (2.792,2.228) {$148.413$};
\draw[gp path] (2.976,2.228)--(3.066,2.228);
\draw[gp path] (11.947,2.228)--(11.857,2.228);
\draw[gp path] (2.976,2.538)--(3.066,2.538);
\draw[gp path] (11.947,2.538)--(11.857,2.538);
\draw[gp path] (2.976,2.849)--(3.066,2.849);
\draw[gp path] (11.947,2.849)--(11.857,2.849);
\draw[gp path] (2.976,3.160)--(3.066,3.160);
\draw[gp path] (11.947,3.160)--(11.857,3.160);
\draw[gp path] (2.976,3.470)--(3.066,3.470);
\draw[gp path] (11.947,3.470)--(11.857,3.470);
\draw[gp path] (2.976,3.781)--(3.156,3.781);
\draw[gp path] (11.947,3.781)--(11.767,3.781);
\node[gp node right] at (2.792,3.781) {$22026.5$};
\draw[gp path] (2.976,3.781)--(3.066,3.781);
\draw[gp path] (11.947,3.781)--(11.857,3.781);
\draw[gp path] (2.976,4.092)--(3.066,4.092);
\draw[gp path] (11.947,4.092)--(11.857,4.092);
\draw[gp path] (2.976,4.402)--(3.066,4.402);
\draw[gp path] (11.947,4.402)--(11.857,4.402);
\draw[gp path] (2.976,4.713)--(3.066,4.713);
\draw[gp path] (11.947,4.713)--(11.857,4.713);
\draw[gp path] (2.976,5.024)--(3.066,5.024);
\draw[gp path] (11.947,5.024)--(11.857,5.024);
\draw[gp path] (2.976,5.334)--(3.156,5.334);
\draw[gp path] (11.947,5.334)--(11.767,5.334);
\node[gp node right] at (2.792,5.334) {$3.26902\times10^{6}$};
\draw[gp path] (2.976,5.334)--(3.066,5.334);
\draw[gp path] (11.947,5.334)--(11.857,5.334);
\draw[gp path] (2.976,5.645)--(3.066,5.645);
\draw[gp path] (11.947,5.645)--(11.857,5.645);
\draw[gp path] (2.976,5.956)--(3.066,5.956);
\draw[gp path] (11.947,5.956)--(11.857,5.956);
\draw[gp path] (2.976,6.266)--(3.066,6.266);
\draw[gp path] (11.947,6.266)--(11.857,6.266);
\draw[gp path] (2.976,6.577)--(3.066,6.577);
\draw[gp path] (11.947,6.577)--(11.857,6.577);
\draw[gp path] (2.976,6.888)--(3.156,6.888);
\draw[gp path] (11.947,6.888)--(11.767,6.888);
\node[gp node right] at (2.792,6.888) {$4.85165\times10^{8}$};
\draw[gp path] (2.976,6.888)--(3.066,6.888);
\draw[gp path] (11.947,6.888)--(11.857,6.888);
\draw[gp path] (2.976,7.198)--(3.066,7.198);
\draw[gp path] (11.947,7.198)--(11.857,7.198);
\draw[gp path] (2.976,7.509)--(3.066,7.509);
\draw[gp path] (11.947,7.509)--(11.857,7.509);
\draw[gp path] (2.976,7.820)--(3.066,7.820);
\draw[gp path] (11.947,7.820)--(11.857,7.820);
\draw[gp path] (2.976,8.130)--(3.066,8.130);
\draw[gp path] (11.947,8.130)--(11.857,8.130);
\draw[gp path] (2.976,8.441)--(3.156,8.441);
\draw[gp path] (11.947,8.441)--(11.767,8.441);
\node[gp node right] at (2.792,8.441) {$7.20049\times10^{10}$};
\draw[gp path] (2.976,8.441)--(3.066,8.441);
\draw[gp path] (11.947,8.441)--(11.857,8.441);
\draw[gp path] (2.976,0.985)--(2.976,1.165);
\draw[gp path] (2.976,8.441)--(2.976,8.261);
\node[gp node center] at (2.976,0.677) {$0$};
\draw[gp path] (4.258,0.985)--(4.258,1.165);
\draw[gp path] (4.258,8.441)--(4.258,8.261);
\node[gp node center] at (4.258,0.677) {$5$};
\draw[gp path] (5.539,0.985)--(5.539,1.165);
\draw[gp path] (5.539,8.441)--(5.539,8.261);
\node[gp node center] at (5.539,0.677) {$10$};
\draw[gp path] (6.821,0.985)--(6.821,1.165);
\draw[gp path] (6.821,8.441)--(6.821,8.261);
\node[gp node center] at (6.821,0.677) {$15$};
\draw[gp path] (8.102,0.985)--(8.102,1.165);
\draw[gp path] (8.102,8.441)--(8.102,8.261);
\node[gp node center] at (8.102,0.677) {$20$};
\draw[gp path] (9.384,0.985)--(9.384,1.165);
\draw[gp path] (9.384,8.441)--(9.384,8.261);
\node[gp node center] at (9.384,0.677) {$25$};
\draw[gp path] (10.665,0.985)--(10.665,1.165);
\draw[gp path] (10.665,8.441)--(10.665,8.261);
\node[gp node center] at (10.665,0.677) {$30$};
\draw[gp path] (11.947,0.985)--(11.947,1.165);
\draw[gp path] (11.947,8.441)--(11.947,8.261);
\node[gp node center] at (11.947,0.677) {$35$};
\draw[gp path] (2.976,8.441)--(2.976,0.985)--(11.947,0.985)--(11.947,8.441)--cycle;
\node[gp node right] at (3.896,8.107) {$|H_s|$};
\gpcolor{rgb color={0.580,0.000,0.827}}
\gpsetpointsize{4.00}
\gp3point{gp mark 1}{}{(3.489,1.262)}
\gp3point{gp mark 1}{}{(3.745,1.423)}
\gp3point{gp mark 1}{}{(4.001,1.642)}
\gp3point{gp mark 1}{}{(4.258,1.850)}
\gp3point{gp mark 1}{}{(4.514,2.075)}
\gp3point{gp mark 1}{}{(4.770,2.294)}
\gp3point{gp mark 1}{}{(5.027,2.522)}
\gp3point{gp mark 1}{}{(5.283,2.747)}
\gp3point{gp mark 1}{}{(5.539,2.977)}
\gp3point{gp mark 1}{}{(5.795,3.205)}
\gp3point{gp mark 1}{}{(6.052,3.437)}
\gp3point{gp mark 1}{}{(6.308,3.668)}
\gp3point{gp mark 1}{}{(6.564,3.900)}
\gp3point{gp mark 1}{}{(6.821,4.133)}
\gp3point{gp mark 1}{}{(7.077,4.367)}
\gp3point{gp mark 1}{}{(7.333,4.601)}
\gp3point{gp mark 1}{}{(7.590,4.836)}
\gp3point{gp mark 1}{}{(7.846,5.071)}
\gp3point{gp mark 1}{}{(8.102,5.307)}
\gp3point{gp mark 1}{}{(8.359,5.542)}
\gp3point{gp mark 1}{}{(8.615,5.779)}
\gp3point{gp mark 1}{}{(8.871,6.015)}
\gp3point{gp mark 1}{}{(9.128,6.252)}
\gp3point{gp mark 1}{}{(9.384,6.489)}
\gp3point{gp mark 1}{}{(9.640,6.727)}
\gp3point{gp mark 1}{}{(9.896,6.964)}
\gp3point{gp mark 1}{}{(10.153,7.202)}
\gp3point{gp mark 1}{}{(10.409,7.440)}
\gp3point{gp mark 1}{}{(10.665,7.678)}
\gp3point{gp mark 1}{}{(10.922,7.916)}
\gp3point{gp mark 1}{}{(11.178,8.155)}
\gp3point{gp mark 1}{}{(4.538,8.107)}
\gpcolor{color=gp lt color border}
\node[gp node right] at (3.896,7.799) {fit};
\gpcolor{rgb color={0.000,0.620,0.451}}
\draw[gp path] (4.080,7.799)--(4.996,7.799);
\draw[gp path] (3.489,1.262)--(3.566,1.331)--(3.644,1.401)--(3.722,1.470)--(3.799,1.540)%
  --(3.877,1.609)--(3.955,1.678)--(4.032,1.748)--(4.110,1.817)--(4.188,1.887)--(4.265,1.956)%
  --(4.343,2.025)--(4.421,2.095)--(4.498,2.164)--(4.576,2.234)--(4.654,2.303)--(4.731,2.372)%
  --(4.809,2.442)--(4.887,2.511)--(4.964,2.581)--(5.042,2.650)--(5.120,2.719)--(5.197,2.789)%
  --(5.275,2.858)--(5.353,2.928)--(5.430,2.997)--(5.508,3.066)--(5.586,3.136)--(5.663,3.205)%
  --(5.741,3.275)--(5.819,3.344)--(5.896,3.413)--(5.974,3.483)--(6.052,3.552)--(6.129,3.622)%
  --(6.207,3.691)--(6.285,3.760)--(6.362,3.830)--(6.440,3.899)--(6.518,3.969)--(6.595,4.038)%
  --(6.673,4.107)--(6.751,4.177)--(6.828,4.246)--(6.906,4.315)--(6.984,4.385)--(7.061,4.454)%
  --(7.139,4.524)--(7.217,4.593)--(7.295,4.662)--(7.372,4.732)--(7.450,4.801)--(7.528,4.871)%
  --(7.605,4.940)--(7.683,5.009)--(7.761,5.079)--(7.838,5.148)--(7.916,5.218)--(7.994,5.287)%
  --(8.071,5.356)--(8.149,5.426)--(8.227,5.495)--(8.304,5.565)--(8.382,5.634)--(8.460,5.703)%
  --(8.537,5.773)--(8.615,5.842)--(8.693,5.912)--(8.770,5.981)--(8.848,6.050)--(8.926,6.120)%
  --(9.003,6.189)--(9.081,6.259)--(9.159,6.328)--(9.236,6.397)--(9.314,6.467)--(9.392,6.536)%
  --(9.469,6.606)--(9.547,6.675)--(9.625,6.744)--(9.702,6.814)--(9.780,6.883)--(9.858,6.953)%
  --(9.935,7.022)--(10.013,7.091)--(10.091,7.161)--(10.168,7.230)--(10.246,7.300)--(10.324,7.369)%
  --(10.401,7.438)--(10.479,7.508)--(10.557,7.577)--(10.634,7.647)--(10.712,7.716)--(10.790,7.785)%
  --(10.867,7.855)--(10.945,7.924)--(11.023,7.994)--(11.100,8.063)--(11.178,8.132);
\gpcolor{color=gp lt color border}
\draw[gp path] (2.976,8.441)--(2.976,0.985)--(11.947,0.985)--(11.947,8.441)--cycle;
\node[gp node center,rotate=-270.0] at (0.292,4.713) {$|H_s|$};
\node[gp node center] at (7.461,0.215) {$s$};
\gpdefrectangularnode{gp plot 1}{\pgfpoint{2.976cm}{0.985cm}}{\pgfpoint{11.947cm}{8.441cm}}
\end{tikzpicture}
}
    \caption{Numerically computed values of $\norm{H_s}$, and their fit $c^{s-2}d$. The vertical axis is plotted in logarithmic scale.}
    \label{fig:Hs_real}
\end{figure}

First, we expand the Hamiltonian~\eqref{eq:H_original_coords} up to order $32$ and apply the diagonalizing change of coordinates~\eqref{eq:linear_change_variables}, obtaining 
\[ H(Q, P)=\sum_{s = 2}^{32} H_{s}(Q, P). \] 
The size of each term $\norm{H_s}$ is shown in Figure~\ref{fig:Hs_real}. Fitting the bound $c^{s-2}d$ to these data allows us to obtain the constants $c=2.09$ and $d=6.63$.

\begin{figure}[!ht]
\centering
	\scalebox{0.9}{\begin{tikzpicture}[gnuplot]
\path (0.000,0.000) rectangle (12.500,8.750);
\gpcolor{color=gp lt color border}
\gpsetlinetype{gp lt border}
\gpsetdashtype{gp dt solid}
\gpsetlinewidth{1.00}
\draw[gp path] (2.668,0.985)--(2.848,0.985);
\draw[gp path] (11.947,0.985)--(11.767,0.985);
\node[gp node right] at (2.484,0.985) {$1$};
\draw[gp path] (2.668,0.985)--(2.758,0.985);
\draw[gp path] (11.947,0.985)--(11.857,0.985);
\draw[gp path] (2.668,1.198)--(2.758,1.198);
\draw[gp path] (11.947,1.198)--(11.857,1.198);
\draw[gp path] (2.668,1.411)--(2.758,1.411);
\draw[gp path] (11.947,1.411)--(11.857,1.411);
\draw[gp path] (2.668,1.624)--(2.758,1.624);
\draw[gp path] (11.947,1.624)--(11.857,1.624);
\draw[gp path] (2.668,1.837)--(2.758,1.837);
\draw[gp path] (11.947,1.837)--(11.857,1.837);
\draw[gp path] (2.668,2.050)--(2.848,2.050);
\draw[gp path] (11.947,2.050)--(11.767,2.050);
\node[gp node right] at (2.484,2.050) {$148.413$};
\draw[gp path] (2.668,2.050)--(2.758,2.050);
\draw[gp path] (11.947,2.050)--(11.857,2.050);
\draw[gp path] (2.668,2.263)--(2.758,2.263);
\draw[gp path] (11.947,2.263)--(11.857,2.263);
\draw[gp path] (2.668,2.476)--(2.758,2.476);
\draw[gp path] (11.947,2.476)--(11.857,2.476);
\draw[gp path] (2.668,2.689)--(2.758,2.689);
\draw[gp path] (11.947,2.689)--(11.857,2.689);
\draw[gp path] (2.668,2.902)--(2.758,2.902);
\draw[gp path] (11.947,2.902)--(11.857,2.902);
\draw[gp path] (2.668,3.115)--(2.848,3.115);
\draw[gp path] (11.947,3.115)--(11.767,3.115);
\node[gp node right] at (2.484,3.115) {$22026.5$};
\draw[gp path] (2.668,3.115)--(2.758,3.115);
\draw[gp path] (11.947,3.115)--(11.857,3.115);
\draw[gp path] (2.668,3.328)--(2.758,3.328);
\draw[gp path] (11.947,3.328)--(11.857,3.328);
\draw[gp path] (2.668,3.541)--(2.758,3.541);
\draw[gp path] (11.947,3.541)--(11.857,3.541);
\draw[gp path] (2.668,3.754)--(2.758,3.754);
\draw[gp path] (11.947,3.754)--(11.857,3.754);
\draw[gp path] (2.668,3.967)--(2.758,3.967);
\draw[gp path] (11.947,3.967)--(11.857,3.967);
\draw[gp path] (2.668,4.180)--(2.848,4.180);
\draw[gp path] (11.947,4.180)--(11.767,4.180);
\node[gp node right] at (2.484,4.180) {$3.26902\times10^{6}$};
\draw[gp path] (2.668,4.180)--(2.758,4.180);
\draw[gp path] (11.947,4.180)--(11.857,4.180);
\draw[gp path] (2.668,4.393)--(2.758,4.393);
\draw[gp path] (11.947,4.393)--(11.857,4.393);
\draw[gp path] (2.668,4.606)--(2.758,4.606);
\draw[gp path] (11.947,4.606)--(11.857,4.606);
\draw[gp path] (2.668,4.820)--(2.758,4.820);
\draw[gp path] (11.947,4.820)--(11.857,4.820);
\draw[gp path] (2.668,5.033)--(2.758,5.033);
\draw[gp path] (11.947,5.033)--(11.857,5.033);
\draw[gp path] (2.668,5.246)--(2.848,5.246);
\draw[gp path] (11.947,5.246)--(11.767,5.246);
\node[gp node right] at (2.484,5.246) {$4.85165\times10^{8}$};
\draw[gp path] (2.668,5.246)--(2.758,5.246);
\draw[gp path] (11.947,5.246)--(11.857,5.246);
\draw[gp path] (2.668,5.459)--(2.758,5.459);
\draw[gp path] (11.947,5.459)--(11.857,5.459);
\draw[gp path] (2.668,5.672)--(2.758,5.672);
\draw[gp path] (11.947,5.672)--(11.857,5.672);
\draw[gp path] (2.668,5.885)--(2.758,5.885);
\draw[gp path] (11.947,5.885)--(11.857,5.885);
\draw[gp path] (2.668,6.098)--(2.758,6.098);
\draw[gp path] (11.947,6.098)--(11.857,6.098);
\draw[gp path] (2.668,6.311)--(2.848,6.311);
\draw[gp path] (11.947,6.311)--(11.767,6.311);
\node[gp node right] at (2.484,6.311) {$7.20049\times10^{10}$};
\draw[gp path] (2.668,6.311)--(2.758,6.311);
\draw[gp path] (11.947,6.311)--(11.857,6.311);
\draw[gp path] (2.668,6.524)--(2.758,6.524);
\draw[gp path] (11.947,6.524)--(11.857,6.524);
\draw[gp path] (2.668,6.737)--(2.758,6.737);
\draw[gp path] (11.947,6.737)--(11.857,6.737);
\draw[gp path] (2.668,6.950)--(2.758,6.950);
\draw[gp path] (11.947,6.950)--(11.857,6.950);
\draw[gp path] (2.668,7.163)--(2.758,7.163);
\draw[gp path] (11.947,7.163)--(11.857,7.163);
\draw[gp path] (2.668,7.376)--(2.848,7.376);
\draw[gp path] (11.947,7.376)--(11.767,7.376);
\node[gp node right] at (2.484,7.376) {$1.06865\times10^{13}$};
\draw[gp path] (2.668,7.376)--(2.758,7.376);
\draw[gp path] (11.947,7.376)--(11.857,7.376);
\draw[gp path] (2.668,7.589)--(2.758,7.589);
\draw[gp path] (11.947,7.589)--(11.857,7.589);
\draw[gp path] (2.668,7.802)--(2.758,7.802);
\draw[gp path] (11.947,7.802)--(11.857,7.802);
\draw[gp path] (2.668,8.015)--(2.758,8.015);
\draw[gp path] (11.947,8.015)--(11.857,8.015);
\draw[gp path] (2.668,8.228)--(2.758,8.228);
\draw[gp path] (11.947,8.228)--(11.857,8.228);
\draw[gp path] (2.668,8.441)--(2.848,8.441);
\draw[gp path] (11.947,8.441)--(11.767,8.441);
\node[gp node right] at (2.484,8.441) {$1.58601\times10^{15}$};
\draw[gp path] (2.668,8.441)--(2.758,8.441);
\draw[gp path] (11.947,8.441)--(11.857,8.441);
\draw[gp path] (2.668,0.985)--(2.668,1.165);
\draw[gp path] (2.668,8.441)--(2.668,8.261);
\node[gp node center] at (2.668,0.677) {$0$};
\draw[gp path] (3.994,0.985)--(3.994,1.165);
\draw[gp path] (3.994,8.441)--(3.994,8.261);
\node[gp node center] at (3.994,0.677) {$5$};
\draw[gp path] (5.319,0.985)--(5.319,1.165);
\draw[gp path] (5.319,8.441)--(5.319,8.261);
\node[gp node center] at (5.319,0.677) {$10$};
\draw[gp path] (6.645,0.985)--(6.645,1.165);
\draw[gp path] (6.645,8.441)--(6.645,8.261);
\node[gp node center] at (6.645,0.677) {$15$};
\draw[gp path] (7.970,0.985)--(7.970,1.165);
\draw[gp path] (7.970,8.441)--(7.970,8.261);
\node[gp node center] at (7.970,0.677) {$20$};
\draw[gp path] (9.296,0.985)--(9.296,1.165);
\draw[gp path] (9.296,8.441)--(9.296,8.261);
\node[gp node center] at (9.296,0.677) {$25$};
\draw[gp path] (10.621,0.985)--(10.621,1.165);
\draw[gp path] (10.621,8.441)--(10.621,8.261);
\node[gp node center] at (10.621,0.677) {$30$};
\draw[gp path] (11.947,0.985)--(11.947,1.165);
\draw[gp path] (11.947,8.441)--(11.947,8.261);
\node[gp node center] at (11.947,0.677) {$35$};
\draw[gp path] (2.668,8.441)--(2.668,0.985)--(11.947,0.985)--(11.947,8.441)--cycle;
\node[gp node right] at (3.588,8.107) {$|H_s|$};
\gpcolor{rgb color={0.580,0.000,0.827}}
\draw[gp path] (3.772,8.107)--(4.688,8.107);
\draw[gp path] (3.198,1.388)--(3.463,1.498)--(3.728,1.649)--(3.994,1.791)--(4.259,1.945)%
  --(4.524,2.096)--(4.789,2.252)--(5.054,2.406)--(5.319,2.564)--(5.584,2.720)--(5.849,2.879)%
  --(6.114,3.037)--(6.380,3.197)--(6.645,3.357)--(6.910,3.517)--(7.175,3.678)--(7.440,3.839)%
  --(7.705,4.000)--(7.970,4.161)--(8.235,4.323)--(8.501,4.485)--(8.766,4.647)--(9.031,4.810)%
  --(9.296,4.972)--(9.561,5.135)--(9.826,5.298)--(10.091,5.461)--(10.356,5.624)--(10.621,5.788)%
  --(10.887,5.951)--(11.152,6.114);
\gpsetpointsize{4.00}
\gp3point{gp mark 1}{}{(3.198,1.388)}
\gp3point{gp mark 1}{}{(3.463,1.498)}
\gp3point{gp mark 1}{}{(3.728,1.649)}
\gp3point{gp mark 1}{}{(3.994,1.791)}
\gp3point{gp mark 1}{}{(4.259,1.945)}
\gp3point{gp mark 1}{}{(4.524,2.096)}
\gp3point{gp mark 1}{}{(4.789,2.252)}
\gp3point{gp mark 1}{}{(5.054,2.406)}
\gp3point{gp mark 1}{}{(5.319,2.564)}
\gp3point{gp mark 1}{}{(5.584,2.720)}
\gp3point{gp mark 1}{}{(5.849,2.879)}
\gp3point{gp mark 1}{}{(6.114,3.037)}
\gp3point{gp mark 1}{}{(6.380,3.197)}
\gp3point{gp mark 1}{}{(6.645,3.357)}
\gp3point{gp mark 1}{}{(6.910,3.517)}
\gp3point{gp mark 1}{}{(7.175,3.678)}
\gp3point{gp mark 1}{}{(7.440,3.839)}
\gp3point{gp mark 1}{}{(7.705,4.000)}
\gp3point{gp mark 1}{}{(7.970,4.161)}
\gp3point{gp mark 1}{}{(8.235,4.323)}
\gp3point{gp mark 1}{}{(8.501,4.485)}
\gp3point{gp mark 1}{}{(8.766,4.647)}
\gp3point{gp mark 1}{}{(9.031,4.810)}
\gp3point{gp mark 1}{}{(9.296,4.972)}
\gp3point{gp mark 1}{}{(9.561,5.135)}
\gp3point{gp mark 1}{}{(9.826,5.298)}
\gp3point{gp mark 1}{}{(10.091,5.461)}
\gp3point{gp mark 1}{}{(10.356,5.624)}
\gp3point{gp mark 1}{}{(10.621,5.788)}
\gp3point{gp mark 1}{}{(10.887,5.951)}
\gp3point{gp mark 1}{}{(11.152,6.114)}
\gp3point{gp mark 1}{}{(4.230,8.107)}
\gpcolor{color=gp lt color border}
\node[gp node right] at (3.588,7.799) {$|Z_s|$};
\gpcolor{rgb color={0.000,0.620,0.451}}
\draw[gp path] (3.772,7.799)--(4.688,7.799);
\draw[gp path] (3.728,1.050)--(4.259,1.014)--(4.789,1.304)--(5.319,1.792)--(5.849,2.311)%
  --(6.380,2.849)--(6.910,3.401);
\gp3point{gp mark 2}{}{(3.728,1.050)}
\gp3point{gp mark 2}{}{(4.259,1.014)}
\gp3point{gp mark 2}{}{(4.789,1.304)}
\gp3point{gp mark 2}{}{(5.319,1.792)}
\gp3point{gp mark 2}{}{(5.849,2.311)}
\gp3point{gp mark 2}{}{(6.380,2.849)}
\gp3point{gp mark 2}{}{(6.910,3.401)}
\gp3point{gp mark 2}{}{(4.230,7.799)}
\gpcolor{color=gp lt color border}
\node[gp node right] at (3.588,7.491) {$|R_s|$};
\gpcolor{rgb color={0.337,0.706,0.914}}
\draw[gp path] (3.772,7.491)--(4.688,7.491);
\draw[gp path] (7.175,4.229)--(7.440,4.497)--(7.705,4.771)--(7.970,5.003)--(8.235,5.235)%
  --(8.501,5.541)--(8.766,5.845)--(9.031,6.098)--(9.296,6.387)--(9.561,6.662)--(9.826,6.989)%
  --(10.091,7.234)--(10.356,7.536)--(10.621,7.811)--(10.887,8.134)--(11.152,8.393);
\gp3point{gp mark 3}{}{(7.175,4.229)}
\gp3point{gp mark 3}{}{(7.440,4.497)}
\gp3point{gp mark 3}{}{(7.705,4.771)}
\gp3point{gp mark 3}{}{(7.970,5.003)}
\gp3point{gp mark 3}{}{(8.235,5.235)}
\gp3point{gp mark 3}{}{(8.501,5.541)}
\gp3point{gp mark 3}{}{(8.766,5.845)}
\gp3point{gp mark 3}{}{(9.031,6.098)}
\gp3point{gp mark 3}{}{(9.296,6.387)}
\gp3point{gp mark 3}{}{(9.561,6.662)}
\gp3point{gp mark 3}{}{(9.826,6.989)}
\gp3point{gp mark 3}{}{(10.091,7.234)}
\gp3point{gp mark 3}{}{(10.356,7.536)}
\gp3point{gp mark 3}{}{(10.621,7.811)}
\gp3point{gp mark 3}{}{(10.887,8.134)}
\gp3point{gp mark 3}{}{(11.152,8.393)}
\gp3point{gp mark 3}{}{(4.230,7.491)}
\gpcolor{color=gp lt color border}
\draw[gp path] (2.668,8.441)--(2.668,0.985)--(11.947,0.985)--(11.947,8.441)--cycle;
\node[gp node center] at (7.307,0.215) {$s$};
\gpdefrectangularnode{gp plot 1}{\pgfpoint{2.668cm}{0.985cm}}{\pgfpoint{11.947cm}{8.441cm}}
\end{tikzpicture}
}
    \caption{Numerically computed values of $\norm{\mathcal{Z}_s}$ and $\norm{\mathcal{R}_s}$, compared to $\norm{H_s}$. The vertical axis is plotted in logarithmic scale.}
    \label{fig:Zs_Rs_real}
\end{figure}
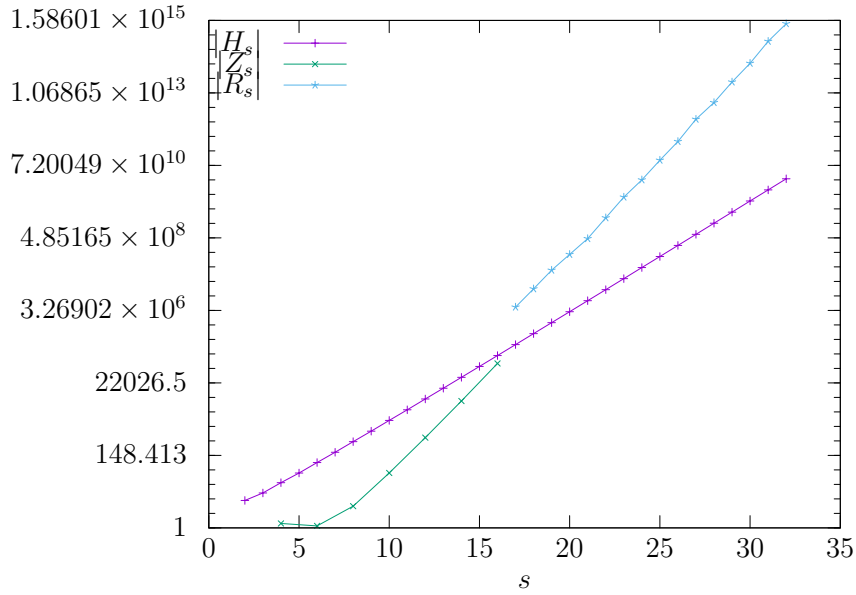

Then, we compute its normal form up to order $M=16$, namely, 
\[ \mathcal{H}^{(16)}=H_{2}+\mathcal{Z}^{(16)}(\bI)+\mathcal{R}^{(16)}(q,p),
\]
with
\[
\mathcal{Z}^{(16)}(\bI)=\sum_{\substack{s \text { even } \\
4 \leq s \leq 16}} \mathcal{Z}_{s}(\bI),\quad  \mathcal{R}^{(16)}(q, p)=\sum_{s = 17}^{32} \mathcal{R}_{s}^{(16)}(q, p).
\]
The size of each homogeneous polynomial $\norm{\mathcal{Z}_s}$ and $\norm{\mathcal{R}_s}$ is shown in Figure~\ref{fig:Zs_Rs_real}.

\begin{figure}[!ht]
\centering
	\scalebox{0.9}{\input{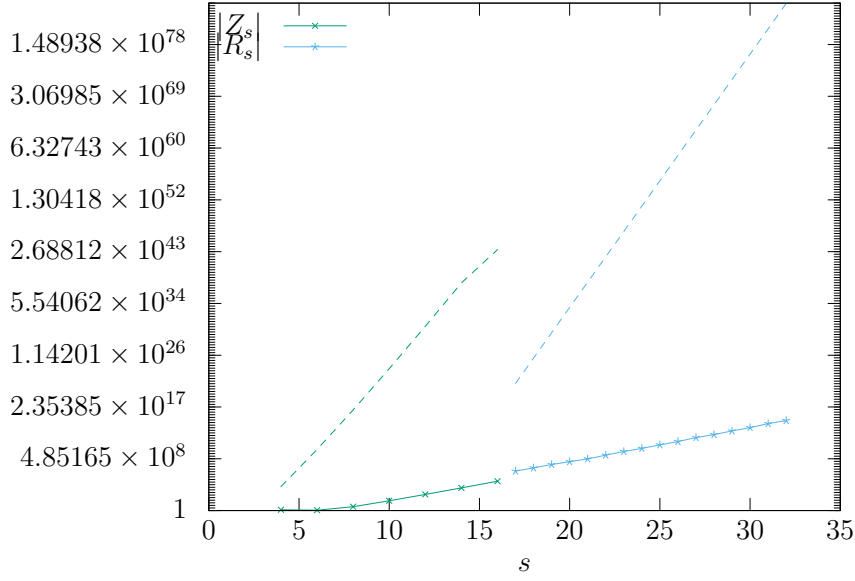}}
    \caption{Numerically computed values of $\norm{\mathcal{Z}_s}$ and $\norm{\mathcal{R}_s}$, compared to their theoretical bounds (dashed lines). The vertical axis is plotted in logarithmic scale.}
    \label{fig:Zs_Rs_real_bounds}
\end{figure}

One can also compare $\norm{\mathcal{Z}_s}$ and $\norm{\mathcal{R}_s}$ with their theoretical bounds, given in Proposition~\ref{prop:bounds}. For the theoretical bounds, we have used the value $a=20$, the numerically fitted coefficients $c=2.09$ and $d=6.63$, as well as the small divisors $\alpha(3)=\cdots=\alpha(16)=0.21372870367940491$.
Figure~\ref{fig:Zs_Rs_real_bounds} compares $\norm{\mathcal{Z}_s}$ and $\norm{\mathcal{R}_s}$ with their theoretical bounds.
Clearly, the theoretical bounds are very conservative.

To complement the theoretical estimates of the error in the normal form, derived in Section~\ref{sec:AnalErrorSEBNF}, we will now measure this error more precisely using semi-analytical methods. Specifically, we will look for the domain of validity $U$ of the computed normal form transformation $\bzeta=\nftrans(\bZ)$, which is used to parameterize all points close to the NHIM $\NHIM$. (See Equation~\eqref{eq:nf}).
Recall from Section~\ref{sec:BNF} that 
\[\NHIM=W_{\mathrm{loc}}^{\cc}=\left\{\bX=(x=0, y=0, I,\phi)\right\}=\{(0,0)\}\times[0, r]^{2}\times\T^2 \]
is the local center manifold of $L_1$ for the normal form system, which is a 4-dimensional NHIM. 

For points $\bZ$ located either in $\NHIM$ or in the stable manifold $W^\st(\NHIM)$, we define the \emph{error in the orbit of $\bZ$} 
\[ \error(\bZ) = \sup_{t\in[0,T=3]} \norm{\nftrans(\flow_t(\bZ)) - \flow_t(\nftrans(\bZ))}. \]
In the first term, $\flow_t(\bZ)$ is computed exactly using Equation~\eqref{eq:flow_NF}. In the second term, $\flow_t(\nftrans(\bZ))$ is computed using numerical integration, more specifically a Runge-Kutta Prince-Dormand (8, 9) method with local error $10^{-16}$. 

\begin{rem}
    The orbit error of $\bZ$ is calculated over a time span of $T=3$ units because this is the approximate period of vertical and horizontal oscillations of solutions around the equilibrium point.
\end{rem}

For ease of visualization, we will summarize the orbit error over all horizontal angles $\phi_\pr$. Let $\bZ=\left(x,y,I,\phi\right)\in U$, and define its \emph{maximal orbit error over $\phi_\pr$}
\[ \bar\error(\bZ) = \sup_{\phi_\pr\in[0,2\pi]} \error(\bZ). \]

\begin{figure}[!ht]
\centering
	\scalebox{0.9}{\input{figs/test-nf-delta_0}}
    \caption{Maximal orbit error $\bar\error(\bX)$ of points $\bX=\left(0,0,I,\phi\right)\in\NHIMc$ as a function of $I_\vr,\phi_\vr$.}
    \label{fig:test-nf-delta_0}
\end{figure}

First, we compute the error on the NHIM $\NHIMc$ of Jacobi constant $C=c:=3.00088$. For each point 
$\bX\in\NHIMc$, 
we compute its maximal orbit error using a normal form of order $M=16$. 
See Figure~\ref{fig:test-nf-delta_0}.

\begin{figure}[!ht]
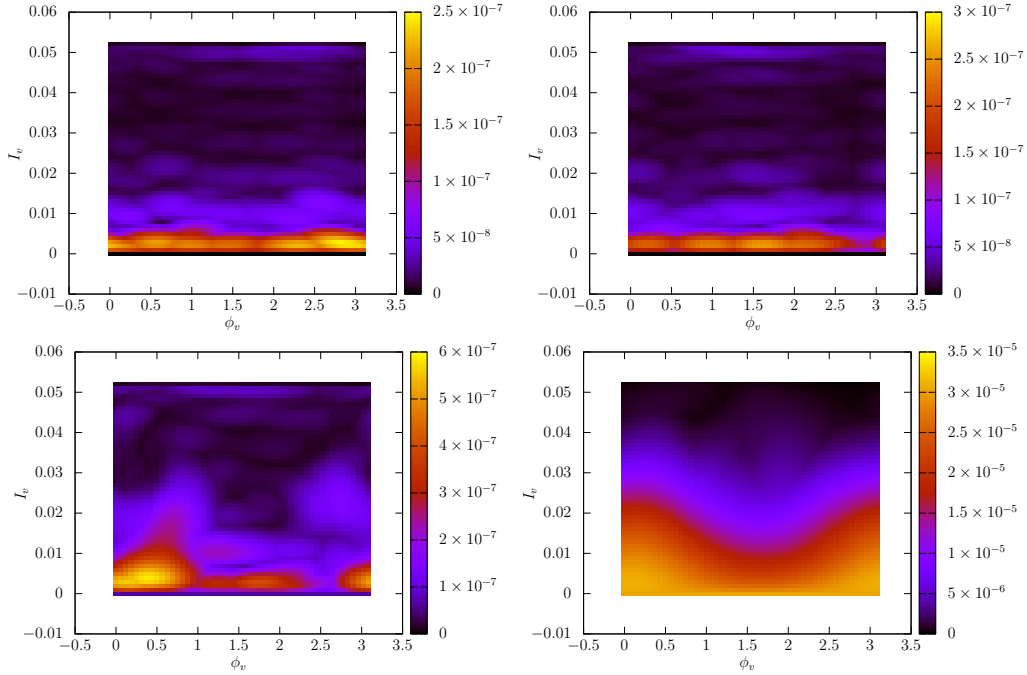

\centering
	\scalebox{0.5}{\input{figs/test-nf-delta_1.e-3}}
	\scalebox{0.5}{\input{figs/test-nf-delta_1.e-1}}
	\scalebox{0.5}{\input{figs/test-nf-delta_0.25}}
	\scalebox{0.5}{\input{figs/test-nf-delta_5.e-1}}
    \caption{Maximal orbit error $\bar\error(\bZ)$ of points $\bZ=(0,\delta,I,\phi)\in W^\st(\Lambda_c)$ for several values $\delta=0.001, 0.1, 0.25, 0.5$. Similar results are obtained for the unstable manifold.}
    \label{fig:test-nf-delta}
\end{figure}

Then, we compute the error on the stable manifold $W^\st(\NHIMc)$.
The maximal orbit error is computed for each point $\bZ=(0,\delta,I,\phi)\in W^\st(\NHIMc)$ located at a distance $\delta$ from the NHIM along its stable fiber. See Figure~\ref{fig:test-nf-delta}.
This experiment is repeated for several values of $\delta$ in order to determine the domain of validity of the normal form of order $16$ to represent the stable manifold. Figure~\ref{fig:test-nf-delta} shows that the stable manifold can be extended to $\delta=0.1$ without significantly affecting the orbit error.
Similar results are obtained for the unstable manifold, integrating backwards over a time span of $T=-3$ time units.


Finally, we compute the error in $U=[-\delta,\delta]^2\times [0,r]^2\times\T^2$, \textit{away} from the invariant manifolds. In this case, it is not guaranteed that $\flow_t(\bZ)$ will stay in the neighborhood $U$ for all $t\in[0,3]$. In particular, the unstable $x$-component grows exponentially, so it can escape $[-\delta,\delta]$. Thus, we slightly modify the definition of the error in the orbit of $\bZ$. Let $\bZ=(x,y,I,\phi)\in U$ with $0<x<\delta$, and define 
\[ \error(\bZ) = \sup_{t\in[0,T]} \norm{\nftrans(\flow_t(\bZ)) - \flow_t(\nftrans(\bZ))}\qquad \text{with }T=\min\{T_\text{exit}, 3\}, \]
where $T_\text{exit}(\bZ)$ is the time to exit $U$, that is, $\flow_{T_\text{exit}}(\bZ)\in\partial U$.
We summarize the orbit error over all actions $I$ and angles $\phi$. Let $\bZ=\left(x,y,I,\phi\right)\in U_c$ be a point with energy value $C=c$, and define
\[ \bar\error(x,y) = \sup_{I\in[0,r]^2, \phi\in[0,2\pi]^2} \error(\left(x,y,I,\phi\right)). \]
\begin{figure}
\centering
    \includegraphics[width=0.7\linewidth]{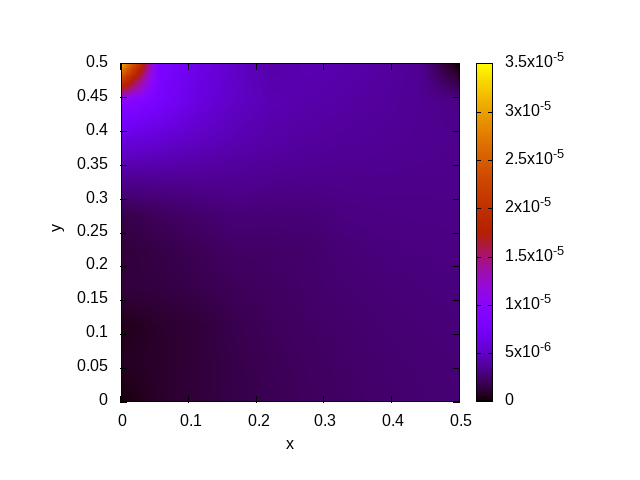}
    \caption{Maximal orbit error $\bar\error(x,y)$ of points $\bZ=(x,y,I,\phi)\in U$ as a function of $x,y$.}
    \label{fig:test-nf-x_sample-y_sample}
\end{figure}
For each point 
$\bZ\in U_c$, 
we compute its maximal orbit error $\bar\error(x,y)$ using a normal form of order $M=16$. 
See Figure~\ref{fig:test-nf-x_sample-y_sample}. This figure shows that the neighborhood $U_c$ can be extended to $\delta=0.1$ without significantly affecting the orbit error, which remains smaller than $3\times 10^{-7}$.

\section{Semi-analytical results}
\label{sec:semianalytical_results}

\subsection{Series representation of the reduced transition maps}
\label{subsec:Series_transition_maps}
Consider the saddle-center-center equilibrium point $L_1$ of the spatial restricted three-body problem (SR3BP). The Birkhoff normal form (BNF) introduced in Section~\ref{sec:BNF} can be used to approximate the flow of the SR3BP in a neighborhood of $L_1$ (see~\cite{GomezMS98,DelshamsGR16, DelshamsGR26}). 
Fix the degree $M=16$ of the normal form and the energy value $C=3.00088$.

Recall that the normal form is computed in a neighborhood $U(\delta,r) = [-\delta, \delta]^2 \times [0,r]^2 \times \T^2$.
As explained in Remark~\ref{rem:time_dependency_tM}, the transition maps depend on the choice of $\delta$. To explore this dependence, we will consider two neighborhoods, a smaller $U(\delta=0.001, r=0.05)$ and a larger $U(\delta=0.1, r=0.05)$. The normal form provides an excellent approximation of the flow in both neighborhoods. As shown in Section~\ref{sec:NumErrorSEBNF}, the maximal orbit error is $2.5\times 10^{-7}$.

\MModify{\noindent\textbf{Notation change:}  
Starting with this section, we rescale the action and write $(J,\phi)$ instead of $(1000\,\bJ_\vr,\phi_\vr)$.}

\begin{figure}
\centering
	\begin{subfigure}{0.49\textwidth}
	    \includegraphics[width=\textwidth]{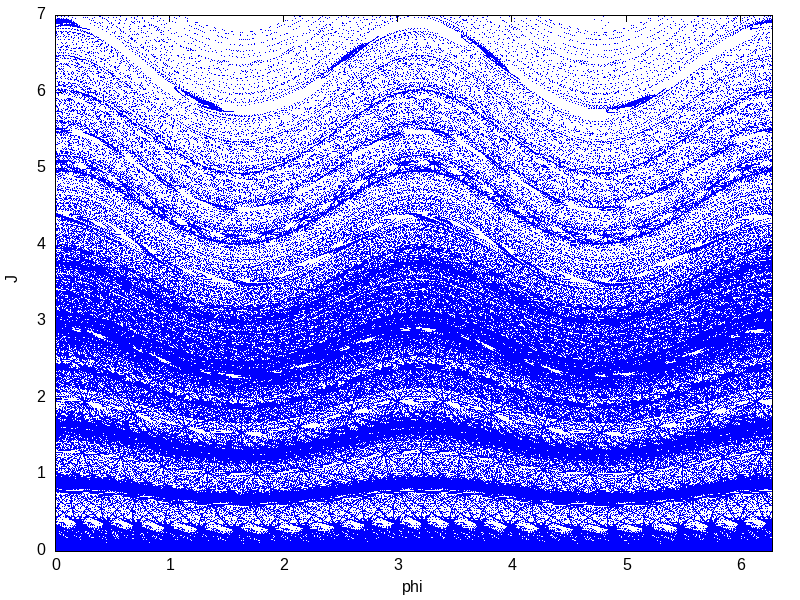}
		\caption{$\tM_1(\delta=0.1)$}
	\end{subfigure}
	\hfill
	\begin{subfigure}{0.49\textwidth}
	    \includegraphics[width=\textwidth]{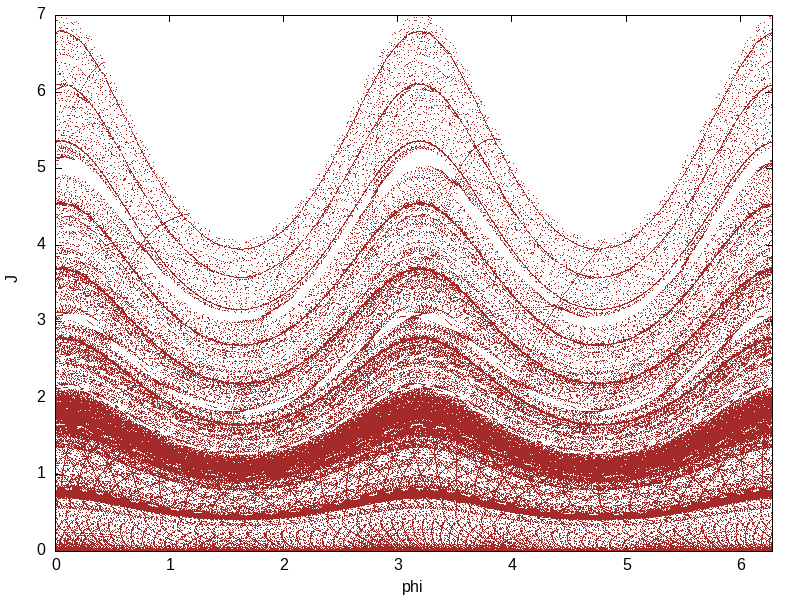}
		\caption{$\tM_2(\delta=0.1)$}
	\end{subfigure}
	~
	\begin{subfigure}{0.49\textwidth}
	    \includegraphics[width=\textwidth]{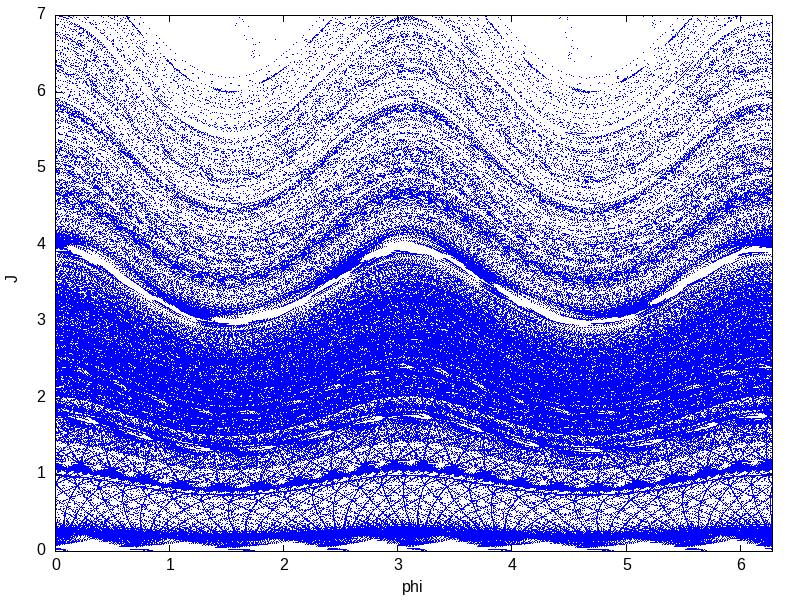}
		\caption{$\tM_1(\delta=0.001)$}
	\end{subfigure}
	\hfill
	\begin{subfigure}{0.49\textwidth}
	    \includegraphics[width=\textwidth]{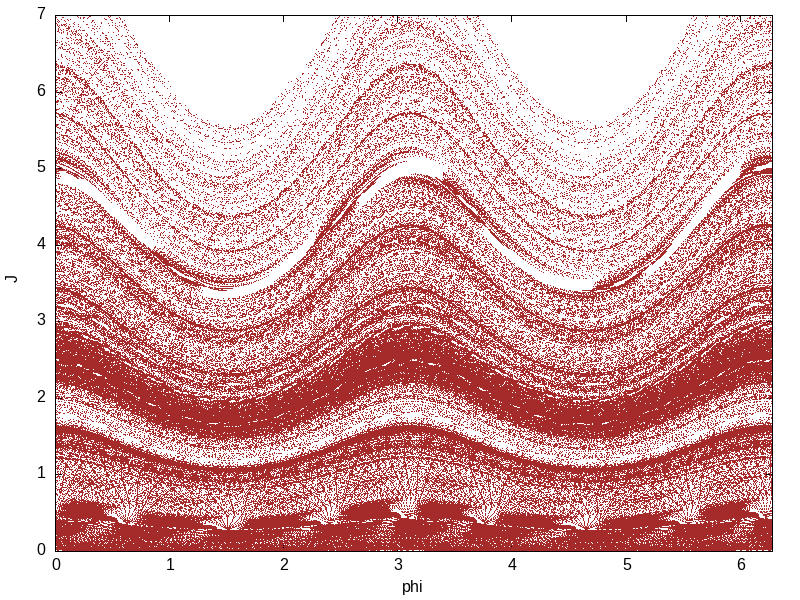}
		\caption{$\tM_2(\delta=0.001)$}
	\end{subfigure}	

	\caption{Phase portrait of the reduced transition maps $\tM_1$ (left) and $\tM_2$ (right) when $\delta=0.1$ (top) and $\delta=0.001$ (bottom).}
	\label{fig:phase_port_SM}
\end{figure}

\MModify{We consider the annulus $\domain \subset \NHIMcred$ given by
\[
\domain = \{ (J, \phi)\colon\ J\in[1,7] \text{ and } \phi\in[0,2\pi) \}.
\]
Using the techniques developed in~\cite{DelshamsGR26}, we compute two transition maps $\tM_{1,2}$ on $\domain$. 
More precisely, 
we start by computing the transition maps on a few invariant tori $J=1,\dotsc,7$.} 
Then we compute a (Fourier-Polynomial) series representation of $\tM_{1,2}$. After that, we can apply the transition maps to any point $(J, \phi)$ in their domain $\domain$, that is, $(J', \phi') = \tM_{1,2}(J, \phi)$. Figure~\ref{fig:phase_port_SM} shows the global phase portraits of $\tM_{1,2}$, obtained by taking a grid of $300\times 300$ initial conditions evenly distributed in the domain $\domain$ and iterating each initial condition $100$ times. See~\cite{DelshamsGR26} for a detailed analysis of the phase portraits.
Some comments are in order:
\begin{itemize}
\PauModify{
    \item The dependence of $\tM_{1,2}$ on the size $\delta$ of the normal form neighborhood $U$ is very small. In this example,  the angle shifts $-(t'_- - t_-)\nuv(J_-)$ and $(t'_+ - t_+)\nuv(J_+)$ mentioned in Remark~\ref{rem:time_dependency_tM} are close to $2\pi$. However, this small difference has dramatic consequences for the phase portraits, such as the size and location of the main resonances (see Figure~\ref{fig:phase_port_SM}). It also has dramatic consequences for the shortest-time guiding orbit, as we shall see in Section~\ref{subsec:Shortest-time_guiding}.
}
    \item $\tM_2$ allows for larger ``jumps'' in terms of action $J$ than $\tM_1$. However, we will show that the shortest-time guiding orbit uses both $\tM_1$ and $\tM_2$.
    
    \item The $J$-component of $\tM_{1,2}$ is a $\pi$-periodic function of the angle, as seen in Figure~\ref{fig:phase_port_SM} and proved in~\cite{DelshamsGR26}. However, the $\phi$-component is not.
\end{itemize}

\subsection{Shortest-time drift guiding orbit}
\label{subsec:Shortest-time_guiding}
Let $\iM\colon \NHIMcred\to\NHIMcred$ be the inner map given by $J'=J,\ \phi'=\phi+\fr(J)$.
In terms of the inner flow of the SR3BP, an application of the inner map corresponds to integrating an initial condition $(J_v, \phi_p=0, \phi_v)\in \NHIMcred$ during the amount of time 
it takes to return to the global section $\SecLocal $ within $U$. Thus, each iterate of the inner map takes time
$\tin = 2\pi/\fr(J)$.
Numerically, we find that $2.0764 < \nu_p(J) < 2.0781$, and therefore the inner time is bounded by
\[ 3.0235 < \tin < 3.0261. \]

In terms of SR3BP flow, one application of the transition map corresponds to a (segment of) a homoclinic trajectory from $\bz_- \in W^u(\NHIMcred)$ to $\bz_+ \in W^s(\NHIMcred)$. 
The transition time $\tout$ from $\bz_-$ to $\bz_+$ depends on the size $\delta$ of the neighborhood $U$ and on the transition map ($\tM_1$ or $\tM_2$). Of course, it also depends on the particular homoclinic, but this dependence is quite mild.
We have computed the transition time $\tout$ for all homoclinics starting on the invariant tori $J=1,2,\dotsc,7$. The results are shown in Table~\ref{tab:tout}. 
\begin{table}[ht!]
	\centering
	\begin{tabular}{lcc}
		\toprule
		& $\tM_1$ & $\tM_2$ \\
		\midrule
		$\delta=0.1$ & $5.45\leq\tout\leq 5.71$ & $4.95\leq\tout\leq 5.21$ \\
		$\delta=0.001$ & $8.32\leq\tout\leq 8.52$ & $8.00\leq\tout\leq 8.21$ \\
		\bottomrule
	\end{tabular}
	\caption{Transition times $\tout$.}
	\label{tab:tout}
\end{table}

Notice that enlarging the neighborhood $U$ by two orders of magnitude (from $\delta=0.001$ to $\delta=0.1$) considerably reduces $\tout$, which then becomes less than twice the inner time $\tin$.

In~\cite{DelshamsGR26}, we introduced different algorithms to produce \textit{drift guiding orbits}, that is, orbits of the iterated function system (IFS) $\left\{ \iM, \tM_1, \tM_2 \right\}$ consisting of the inner map $\iM$ and transition maps $\tM_{1,2}$, whose action variable $J$ increases from $J=1$ to $J>7$. 
Using the previous time estimates, we can estimate the \textit{total drift time} of the guiding orbit as
\[ t = n_0 \tin + n_1 \tout + n_2 \tout, \]
where $n_0$, $n_1$ and $n_2$ denote the number of iterations of $\iM$, $\tM_1$ and $\tM_2$, respectively.
In particular, the Shortest-Time Algorithm (\cite{DelshamsGR26}, Section 12.2) looks for a drift guiding orbit with the shortest total time.

\begin{figure}[!ht]

	\begin{subfigure}{\textwidth}
	\centering		\includegraphics[width=0.7\textwidth]{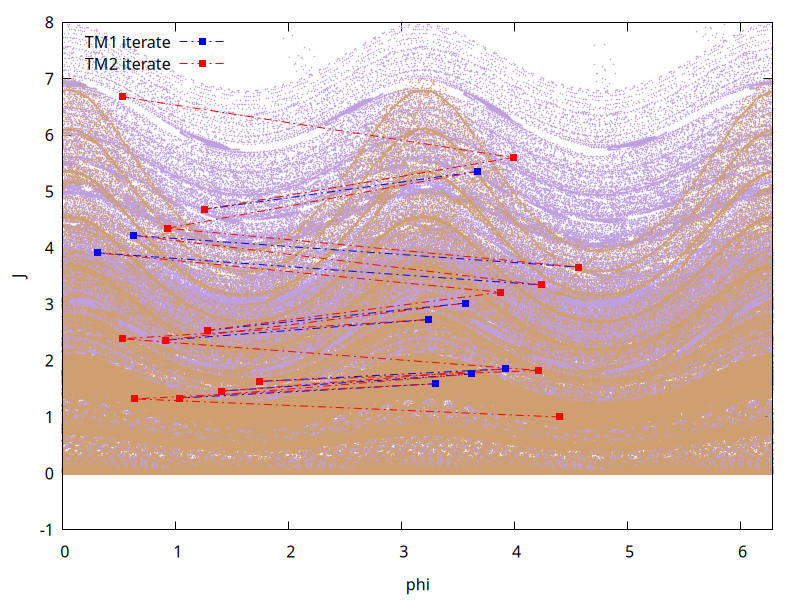}
		\caption{$\delta=0.1$. Total drift time $121.3$.}
		\label{fig:diffusion-shortest-path-delta-0.1}
	\end{subfigure}
	~
	\begin{subfigure}{\textwidth}
			\centering
		\includegraphics[width=0.7\textwidth]{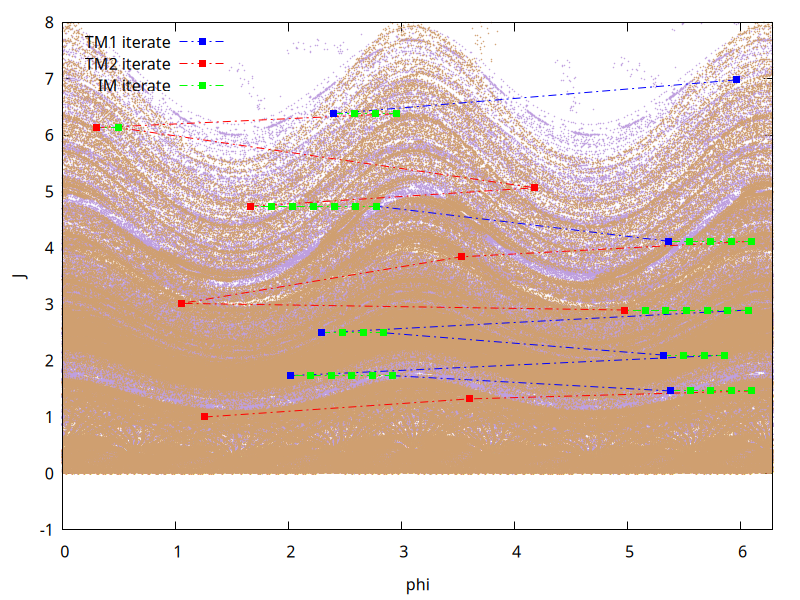}
		\caption{$\delta=0.001$. Total drift time $220.2$.}
		\label{fig:diffusion-shortest-path-delta-0.001}
	\end{subfigure}
	
	\caption{Shortest-time drift orbit found by iterating the triple dynamical system $\{\iM, \tM_1, \tM_2\}$ when $\delta=0.1$ (above) and $\delta=0.001$ (below). A point $x_i$ in the orbit is plotted in blue if $\tM_1(x_i)=x_{i+1}$, red if $\tM_2(x_i)=x_{i+1}$, and green if $\iM(x_i)=x_{i+1}$. Consecutive iterates are joined by dashed segments to help visualize the orbit.}
	\label{fig:diffusion_shortest_path}
\end{figure}

Consider first the case $\delta=0.1$, with an IFS $\left\{ \iM, \tM_1(\delta=0.1), \tM_2(\delta=0.1) \right\}$. In this case, the Shortest-Time Algorithm produces the guiding orbit of Figure~\ref{fig:diffusion-shortest-path-delta-0.1} with a total drift time $121.3$ in normalized time units, equivalent to $19.3$ years.
Consider next the case $\delta=0.001$, with an IFS $\left\{ \iM, \tM_1(\delta=0.001), \tM_2(\delta=0.001) \right\}$. In this case, the Shortest-Time Algorithm produces the guiding orbit of Figure~\ref{fig:diffusion-shortest-path-delta-0.001} with a total drift time $220.2$, or $35$ years.

The resulting guiding orbits are strikingly different in each case. When $\delta=0.001$, it uses all three maps $\iM, \tM_1, \tM_2$. In particular, the guiding orbit uses the inner map $\iM$ repeatedly to decrease the angle $\phi$ until it reaches a location where $\tM_{1,2}$ can be applied to increase action $J$.

In contrast, when $\delta=0.1$, the shortest-time orbit does not use the inner map at all: it only uses $\tM_1$ and $\tM_2$. Recall that $\tM_2$ makes larger ``jumps'' in terms of action than $\tM_1$. The algorithm manages to find a path where $\tM_2$ is applied often, always to increase the action, even at the expense of applying $\tM_1$ too, always to \textit{decrease} the action but not as much.

Another factor that explains the qualitative difference between the shortest-time orbits is the ratio between the inner time $\tin$ and the transition time $\tout$. When $\delta=0.001$, $\tout$ is much larger than $\tin$, so $\tM_{1,2}$ are used sparingly. When $\delta=0.1$, $\tout$ is not much larger than $\tin$, so $\tM_{1,2}$ can be used more often.

In the first case, when $\delta=0.001$, the obtained shortest-time guiding orbit are orbits of the IFS $\{\iM, \tM_1, \tM_2\}$. 
In the second case, when $\delta=0.1$, the obtained shortest-time guiding orbit are orbits of the IFS $\{\tM_1, \tM_2\}$. Lemma~\ref{lem:shadowing_topological} and Theorem~\ref{th:shadowing} give conditions that guarantee that such guiding orbits can be shadowed by true orbits.

While Lemma \ref{lem:shadowing_topological} and Theorem~\ref{th:shadowing} are general theoretical results that guarantee the existence of shadowing orbits under some conditions, the multiple shooting scheme presented in  
Section~\ref{sec:multiple_shooting} provides the explicit numerical construction of such shadowing orbits. 

\begin{rem}
The scattering map cannot be used to design optimal time guiding orbits. It is defined through asymptotic homoclinics, so it does not carry any information on time. In this situation, it is necessary to use the transition map.
\end{rem}

\subsection{Shadowing trajectory}
\label{sec:Shadowing_trajectory}
Consider the case $\delta=0.1$. Let the guiding orbit consist of the last $N=17$ iterates of the shortest-time drift orbit given in Figure~\ref{fig:diffusion-shortest-path-delta-0.1}, starting at $\bx_1=(\phi,J)=(1.743404, 1.630975)$ and ending at $\bx_{17}=(\phi,J)=(6.677670, 0.528935)$. Each iterate in the guiding orbit $\{\bx_i\in\NHIMcred\}_{i=1,\dotsc,N}$ satisfies either $\tM_1(\bx_i)=\bx_{i+1}$ or $\tM_2(\bx_i)=\bx_{i+1}$. Our goal is to shadow this guiding orbit using the method of Section~\ref{sec:multiple_shooting}.

\begin{figure}[!ht]
\centering
    \begin{subfigure}{0.7\textwidth}
	\centering
	    \resizebox{\textwidth}{!}{
        \small
        \input{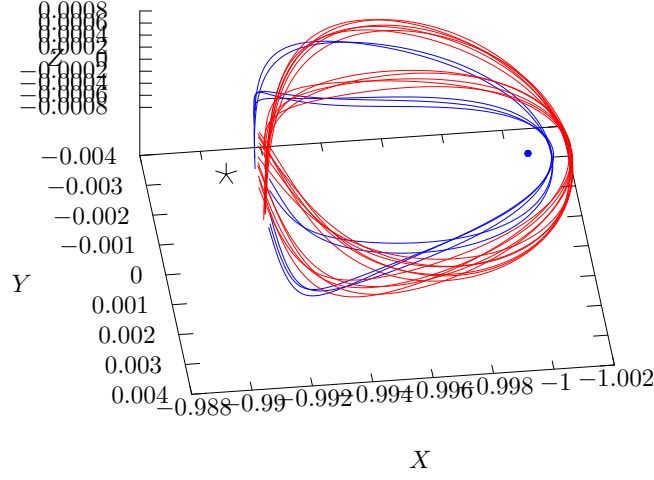}
        }
		\caption{Homoclinic trajectory segments ending on the pseudo-trajectory $(\bZ^-_i, \bZ^+_{i+1})_{i=1}^{N-1}$}
	\end{subfigure}
	\begin{subfigure}{0.7\textwidth}
	    \resizebox{\textwidth}{!}{
        \small
        \input{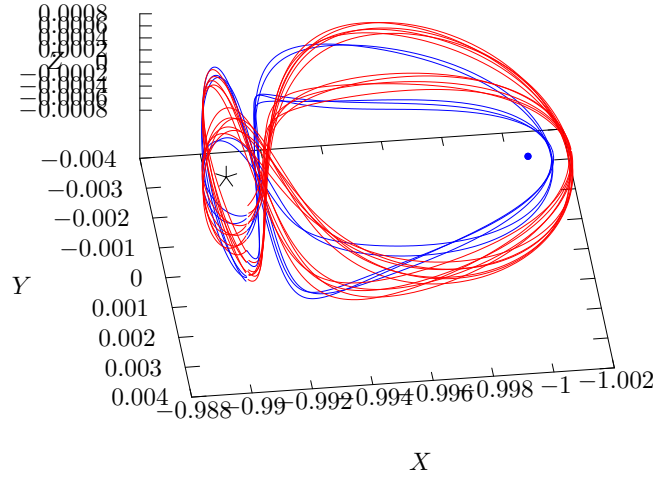}
        }
		\caption{Homoclinic segments ending on the reduced pseudo-trajectory $\{(\bz_-^i, \bz_+^{i+1})\}_{i=1}^{N-1}$}
        \label{fig:pseudo-orbit-refined_2}
	\end{subfigure}
\caption{Sequence of $N-1=16$ homoclinic trajectory segments (in original SR3BP coordinates) associated to the guiding orbit. Blue segments are associated to the transition map $\tM_1$, red segments to $\tM_2$. The Earth location is marked by the blue circle, while the equilibrium point $L_1$ is marked by a black star.}
\label{fig:pseudo_orbit_refined}
\end{figure}

To this guiding orbit, we associate a sequence of $N-1=16$ $\NHIMc$-homoclinic trajectory segments. 
When the segments end at the boundary $\partial U$ of the normal form domain, the endpoints $(\bZ^-_i, \bZ^+_{i+1})_{i=1}^{N-1}$ make up the pseudo-trajectory. 
If the segments are further extended to the global section $\SecLocal =\{\phi_{p,\pm}=0\}$ within $U$, then the endpoints $(\bz^-_i, \bz^+_{i+1})_{i=1}^{N-1}$ make up the reduced pseudo-trajectory. 
See Figure~\ref{fig:pseudo_orbit_refined}. This reduced pseudo-trajectory has a total drift time $83.5896$ adimensional time units, or 13.3 years.

As explained in Remark~\ref{rem:reduced_pseudo_orbit}, the reduced pseudo-trajectory determines a pseudo-orbit in the classical sense~\cite{Bowen75}, that is, the last point of each trajectory segment is close to the first point of the next segment. This is clear in Figure~\ref{fig:pseudo-orbit-refined_2}: there are small (but visible) gaps between consecutive segments, gaps located around $(X=-0.99, Y=0.004)$. See a closeup of the gaps in Figure~\ref{closecasebguess}.

The error of this pseudo-trajectory can be estimated as the total size of these gaps. Let $\bz^+_i=(x^+_i,v^+_i)$ denote the last point of a segment in position-velocity coordinates, and $\bz^-_i=(x^-_i,v^-_i)$ the first point of the next segment. 
Each of these gaps has an average position error $\norm{\Delta x_i}=\norm{x^+_i - x^-_i} \approx 8.2\times 10^{-5}$ adimensional units, and an average velocity error $\norm{\Delta v_i}=\norm{v^+_i - v^-_i} \approx 2.1\times 10^{-4}$ adimensional units.
Since there are 15 gaps, the total position error is $\norm{\Delta x} = 1.23\times 10^{-3}$, and total velocity error $\norm{\Delta v} = 3.15\times 10^{-3}$. In the Sun-Earth system, the adimensional length unit corresponds to $1$ AU (Astronomical Unit) = $1\,496\times 10^{8}$ km, while $2\pi$ adimensional time units correspond to 1 year. Therefore, the total position and velocity errors amount to   $\norm{\Delta x} = 184\,005$ km and $\norm{\Delta v} = 338\,000 \text{ km/h} = 93.89 \text{ km/s}$.

Finally, we apply the multiple shooting method described in Section~\ref{sec:multiple_shooting}. We look for $2(N-1)=32$ points $(\bK^-_i, \bK^+_{i+1})_{i=1}^{N-1}$ in the boundary sections $\SecPlusMinus$ that satisfy the matching conditions C.1' and C.2. This corresponds to solving the system of equations \eqref{eq:outer2}, \eqref{eq:inner} and~\eqref{eq:extra_conditions}, consisting of $10N-10=160$ real equations in $10N-10=160$ variables. We use the multidimensional Newton's method to solve the system. 

Since the shadowing trajectory is hyperbolic and thus unstable, this requires the initial guess to be very close to an actual but unknown trajectory. Indeed, we have used the initial guess proposed in Section~\ref{sec:initial_guess}.

A root of Newton's method is found with absolute error $10^{-12}$. This implies that the total position and velocity error of the shadowing trajectory is within $10^{-12}$; that is, $\norm{\Delta x} < 10^{-12}$, $\norm{\Delta v} < 10^{-12}$. In physical units, this amounts to $\norm{\Delta x} < 1.49598$ m and $\norm{\Delta v} < 2.9806\times 10^{-5}$ m/s.
Compare these errors to those of the differential correction method used in Section~\ref{sec:two_level_diff_corr} ($29$ m and $0.006$ m/s). Clearly, the multiple-shooting method produces a more precise shadowing trajectory.

\begin{figure}[!ht]
\centering
	\input{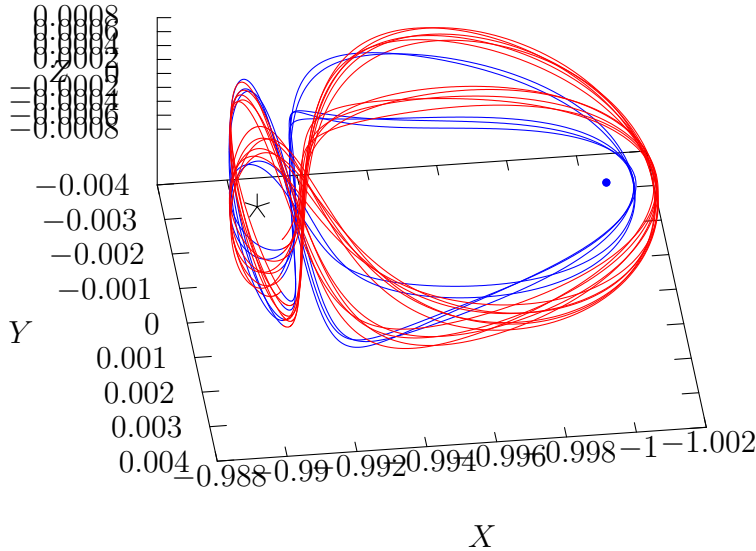}
    \caption{Shadowing trajectory corresponding to the pseudo-trajectory of  Figure~\ref{fig:pseudo_orbit_refined}. The small gaps that were present in the pseudo-trajectory have been closed.}
    \label{fig:shadowing_trajectory}
\end{figure}

The resulting shadowing trajectory is shown in Figure~\ref{fig:shadowing_trajectory}. Notice that the small gaps present in the pseudo-trajectory~\ref{fig:pseudo_orbit_refined} have been closed.
The shadowing trajectory has a total drift time $83.5874$, which is very similar to that of the pseudo-trajectory. 

\begin{table}[ht!]
	\centering
    \resizebox{\textwidth}{!}{%
	\begin{tabular}{llllllllll}
		\toprule
		$x$ & $y$ & $I_\pr$ & $I_\vr$ & $\phi_\pr$ & $\phi_\vr$ & $t_+$ & $\Tint$ & $t_+ + \frac{2\pi}{\nuv}$ \\
		\midrule
1.209337e-04 & -1.000000e-01 & 4.836035e-02 & 2.003625e-03 & 5.065988e+00 & 2.829394e+00 & 1.655701e+00 & 2.652066e+00 & 4.667774e+00 \\
1.202584e-04 & -1.000000e-01 & 4.782276e-02 & 2.557167e-03 & 5.064408e+00 & 5.473926e+00 & 3.879481e-01 & 2.654277e+00 & 3.400022e+00 \\
4.869715e-05 & -1.000000e-01 & 4.765769e-02 & 2.718287e-03 & 5.066777e+00 & 1.912417e+00 & 2.095287e+00 & 3.011174e+00 & 5.107360e+00 \\
4.833341e-05 & -1.000000e-01 & 4.792152e-02 & 2.446521e-03 & 4.313133e+00 & 5.124596e+00 & 5.554119e-01 & 3.014134e+00 & 3.567485e+00 \\
4.881234e-05 & -1.000000e-01 & 4.736196e-02 & 3.022895e-03 & 5.056019e+00 & 2.243910e+00 & 1.936374e+00 & 3.010242e+00 & 4.948447e+00 \\
4.815479e-05 & -1.000000e-01 & 4.782672e-02 & 2.544133e-03 & 4.322214e+00 & 5.506668e+00 & 3.722520e-01 & 3.015596e+00 & 3.384326e+00 \\
1.142642e-04 & -1.000000e-01 & 4.712836e-02 & 3.271653e-03 & 5.041091e+00 & 2.545108e+00 & 1.791984e+00 & 2.674462e+00 & 4.804057e+00 \\
4.897135e-05 & -1.000000e-01 & 4.657093e-02 & 3.837777e-03 & 5.040550e+00 & 5.271619e+00 & 4.849309e-01 & 3.008958e+00 & 3.497004e+00 \\
4.814316e-05 & -1.000000e-01 & 4.707655e-02 & 3.316801e-03 & 4.338323e+00 & 2.231280e+00 & 1.942428e+00 & 3.015691e+00 & 4.954502e+00 \\
4.903114e-05 & -1.000000e-01 & 4.623255e-02 & 4.186421e-03 & 5.021040e+00 & 5.580117e+00 & 3.370415e-01 & 3.008476e+00 & 3.349115e+00 \\
4.798710e-05 & -1.000000e-01 & 4.681267e-02 & 3.588618e-03 & 4.356312e+00 & 2.590640e+00 & 1.770156e+00 & 3.016973e+00 & 4.782229e+00 \\
1.040047e-04 & -1.000000e-01 & 4.604983e-02 & 4.381515e-03 & 5.000651e+00 & 5.869450e+00 & 1.983390e-01 & 2.711603e+00 & 3.210412e+00 \\
4.936078e-05 & -1.000000e-01 & 4.511055e-02 & 5.342652e-03 & 4.993494e+00 & 2.334512e+00 & 1.892940e+00 & 3.005831e+00 & 4.905014e+00 \\
4.783102e-05 & -1.000000e-01 & 4.577267e-02 & 4.660125e-03 & 4.384510e+00 & 5.605969e+00 & 3.246485e-01 & 3.018259e+00 & 3.336722e+00 \\
9.506716e-05 & -1.000000e-01 & 4.481874e-02 & 5.649086e-03 & 4.963041e+00 & 2.627957e+00 & 1.752266e+00 & 2.747075e+00 & 4.764340e+00 \\
		\bottomrule
	\end{tabular}}
	\caption{Points $K_+^i\in\SecPlus$ in the shadowing trajectory, for $i=2,\cdots,16$. We also tabulate the time $t_+$ to first cut, the interior time $\Tint$, and the time $t_+ + \frac{2\pi}{\nuv}$ to second cut.}
	\label{tab:OneCut}
\end{table}

After obtaining the shadowing trajectory, we check that it satisfies condition~\eqref{eq:OneCut0}, which guarantees cutting the section $\SecLocal$ exactly once during each passage through the BNF neighborhood $U$. 
As shown in Table~\ref{tab:OneCut}, the one-cut condition is verified for every point $K_+^i\in\SecPlus$ in the shadowing trajectory.

\begin{rem}\label{rem:shadowing_thm}
Theorem~\ref{th:shadowing} does not apply in this situation. This particular guiding orbit has $N=17$. All other variables ($\LE, \lambda, \nup$) depend on the homoclinic point $\bzeta^{i}$, so they vary slightly for each homoclinic segment. However, in practice, we find that these variables do not change much across homoclinics, and remain very close to the following values: $\LE \approx 2.79$, $\lambda \approx \lambdaL = 2.53$, $\nup \approx \nupL = 2.08$. Thus, the homecoming condition~\eqref{eq:HomecomingCondition} is not satisfied,
\[ (N-1)\LE \approx 44.64 \nless \pi\lambda/\nup \approx 3.82. \]
\end{rem}

\section{Two-level differential correction of pseudo-orbits}
\label{sec:two_level_diff_corr}

In this section, we describe a two-level differential correction method to reduce the discontinuities present in a pseudo-trajectory and apply it to the pseudo-trajectory of Section~\ref{sec:semianalytical_results}. 
After convergence of the differential correction method, the converged trajectory is compared with the shadowing trajectory of Section~\ref{sec:semianalytical_results}.

\subsection{Methodology}

A brief overview of the method is given first. More
detailed descriptions of this approach may be found in
Wilson~\cite{Wilson:2003}, Pernicka~\cite{Pernicka:1986}, or Parker and
Anderson~\cite{Parker:2014}.   Using this approach,
a level-1 step is taken where velocities at patch points along the trajectory are
updated followed by a level-2 step where positions and times are updated.

\subsubsection*{Level 1}

In this step, the variational equations, or the state transition matrix (STM),
are computed between a series of \emph{patch points} that have been obtained
from another method or model and placed in the current model.  The STM is used
to update velocities at each patch point so as to target the position at the
next patch point.

\subsubsection*{Level 2}

In the level-2 step, positions and times are updated using a least squares
approach.  Following the derivation in Anderson and Parker, this is computed using
\begin{equation*}
\begin{bmatrix}
  \delta r_1 \\
  \delta t_1 \\
  \delta r_2 \\
  \delta t_2 \\
  \delta r_3 \\
  \delta t_3 \\
  \vdots \\
  \delta r_n \\
  \delta t_n \\
  \end{bmatrix}
  =
  M^T (M M^T)^{-1}
  \begin{bmatrix}
  \delta \Delta v_2 \\
  \delta \Delta v_3 \\
  \vdots \\
  \delta \Delta v_{n-1} \\
  \end{bmatrix}
  .
\end{equation*}
Here, $\delta r_n$ is the discontinuity in position at patch point $n$, $\delta t_n$
the discontinuity in time, and $\delta \Delta v_n$ the discontinuity in velocity.  $M$ is 
computed from 
\begin{equation*}
M = 
\begin{bmatrix}                                                                                             
    \frac{\partial \Delta v_2}{\partial r_1} & \frac{\partial \Delta v_2}{\partial t_1} & \frac{\partial
  \Delta v_2}{\partial r_2} & \frac{\partial \Delta v_2}{\partial t_2} & \frac{\partial \Delta v_2}{\partial
  r_3} & \frac{\partial \Delta v_2}{\partial t_3} & 0 & 0 & \cdots & 0 & 0\\
    0 & 0 & \frac{\partial \Delta v_3}{\partial r_2} & \frac{\partial \Delta v_3}{\partial t_2} & \frac{\partial \Delta v_3}{\partial r_3} & \frac{\partial \Delta v_3}{\partial t_3} & \frac{\partial \Delta v_3}{\partial r_4} & \frac{\partial \Delta v_3}{\partial t_4} & \cdots & 0 & 0\\
    \vdots & \vdots & \vdots & \ddots & \vdots & \ddots & \vdots & \vdots & \vdots & \vdots & \vdots\\
    0 & 0 & 0 & 0 & \cdots &  \frac{\partial \Delta v_{n-1}}{\partial r_{n-2}} & \frac{\partial \Delta v_{n-1}}{\partial t_{n-2}} & \frac{\partial
  \Delta v_{n-1}}{\partial r_{n-1}} & \frac{\partial \Delta v_{n-1}}{\partial t_{n-1}} & \frac{\partial \Delta v_{n-1}}{\partial
  r_n} & \frac{\partial \Delta v_{n-1}}{\partial t_n}
  \end{bmatrix}                 
\end{equation*}
This adjusts the positions and times, which introduces discontinuities in the
velocities at each patch point.  

\subsubsection*{Overall iteration}

The entire process is repeated converging first with the level-1 step until a
desired tolerance is met and then with the level-2 step.  This is typically
continued until the changes in the positions become very small within some
tolerance, or until the discontinuities satisfy a convergence criterion.   The
number of patch points is often selected heuristically and modified depending
on how well the process is converging.  Constraints may also be added for each
patch point as needed for the problem or to aid in convergence.

\subsection{Application to SR3BP pseudo-orbit}

The two-level differential corrector was applied to the pseudo-orbit of Section~\ref{sec:Shadowing_trajectory} to
obtain a nearly continuous trajectory and remove position and velocity
discontinuities down to a small tolerance.  
The only constraint
was that the initial and final positions were kept constant. 
The initial case was run with a larger number of patch points, and the number of patch points
was reduced as needed.

The pseudo-orbit of Section~\ref{sec:Shadowing_trajectory} was used as the initial trajectory (see Figure~\ref{fig:pseudo-orbit-refined_2}).
This trajectory consists of homoclinic segments that start and end
at a large distance $\delta=0.1$ from the center manifold (see Section~\ref{subsec:Shortest-time_guiding}).  
This reduces the duration of the homoclinic
segments, but introduces larger $\Delta$Vs or discontinuities.  
For this case,
16 heteroclinics were used with 64 total sub-segments, and the time spent near L$_1$ is reduced compared to case $\delta=0.001$.  
The pseudo-trajectory has total position and velocity discontinuities of $184\,005$ km and $338\,000 \text{ km/h}$, as described in Section~\ref{sec:Shadowing_trajectory}.

\begin{figure}[!ht]
	\centering
		\includegraphics[width=0.4\textwidth]{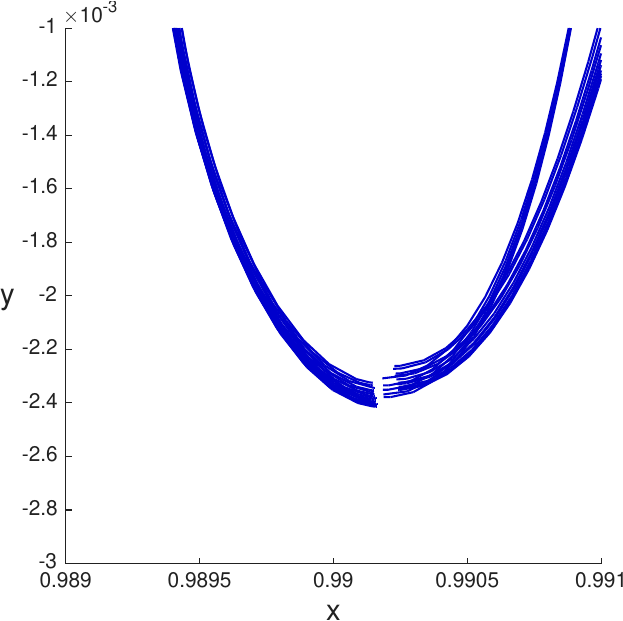}
		\caption{Close xy view of the pseudo-orbit}
		\label{closecasebguess}
\end{figure}

The position discontinuities can be seen more
clearly in the closer view in Figure \ref{closecasebguess}.
Given the large
numbers of segments and discontinuities, this is typically a very challenging
case for differential correction.

\begin{figure}[ht!]
	\centering
	\begin{subfigure}{0.45\textwidth}
		\includegraphics[width=\textwidth]{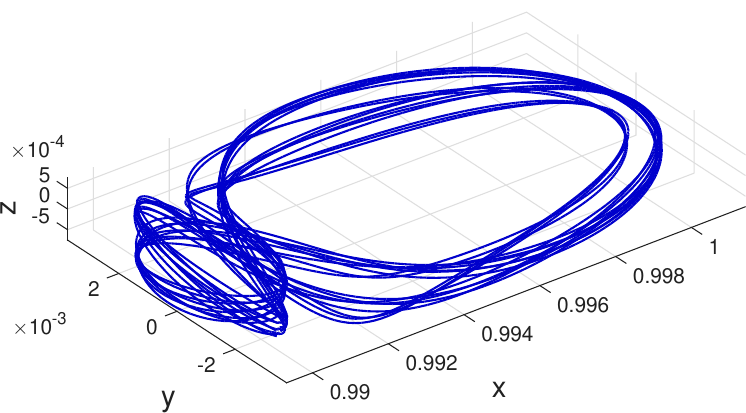}
		\caption{Pseudo-Orbit}
	\end{subfigure}
	\begin{subfigure}{0.45\textwidth}
		\includegraphics[width=\textwidth]{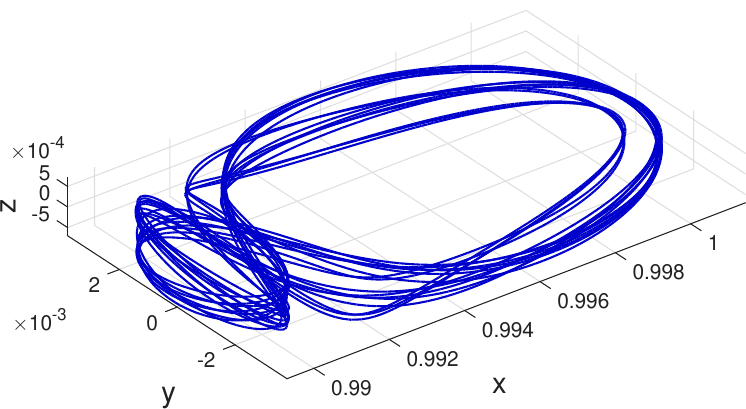}
		\caption{Converged Trajectory}
	\end{subfigure}
	\caption{
		The short trajectory converged using the two-level differential corrector
		compared to the original pseudo-orbit.
	}
	\label{shortcompare}
\end{figure}
\begin{figure}[ht!]
	\centering
	\includegraphics[width=0.55\textwidth]{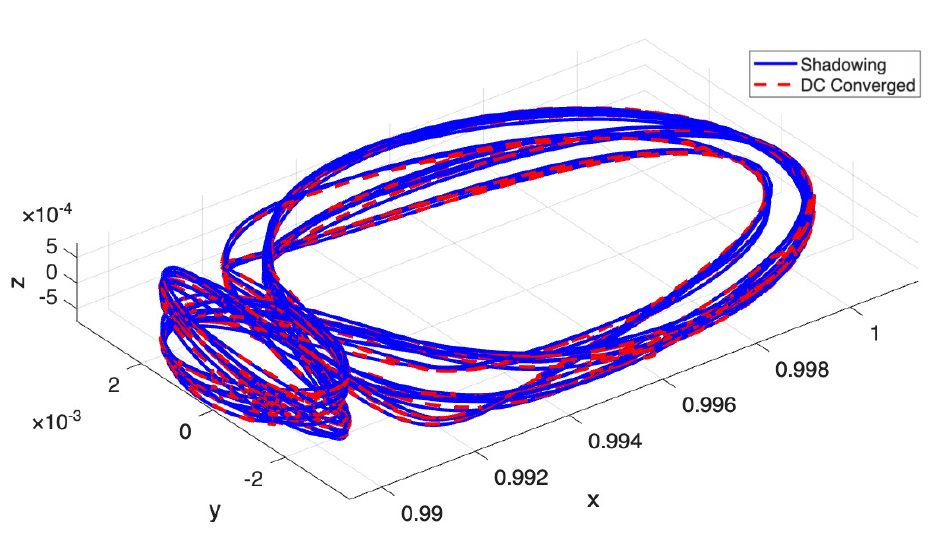}
	\caption{
		Comparison of the trajectory obtained using the two-level differential
		corrector and the shadowing method.
	}
	\label{compareall}
\end{figure}
\begin{figure}[ht!]
	\centering
	\begin{subfigure}{0.45\textwidth}
		\includegraphics[width=\textwidth]{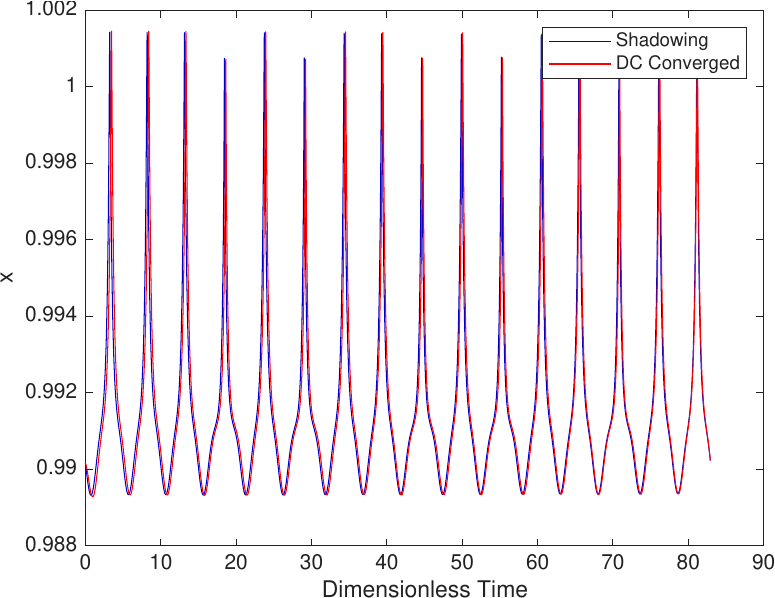}
		\caption{x Difference}
	\end{subfigure}
	\begin{subfigure}{0.45\textwidth}
		\includegraphics[width=\textwidth]{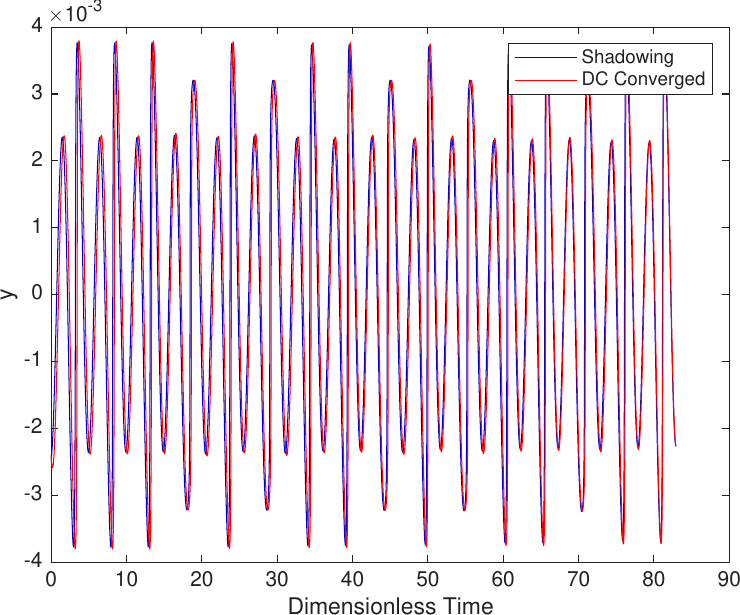}
		\caption{y Difference}
	\end{subfigure}
	\begin{subfigure}{0.45\textwidth}
		\includegraphics[width=\textwidth]{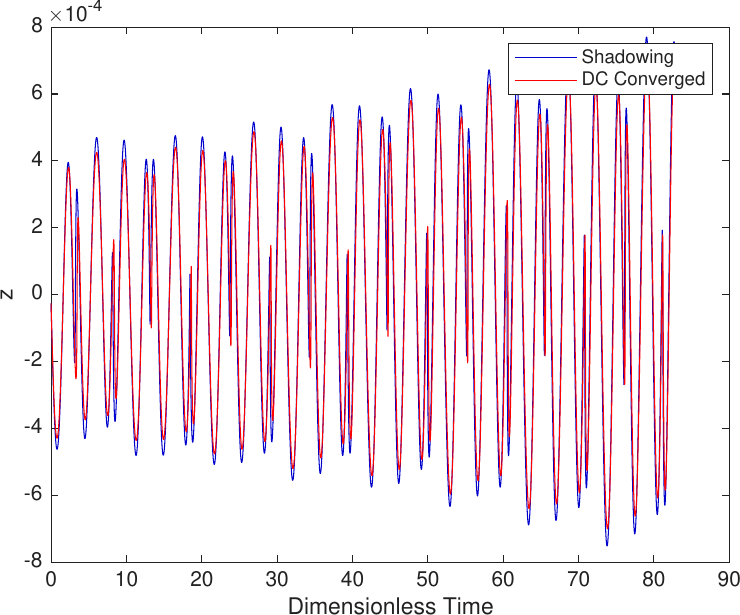}
		\caption{z Difference}
	\end{subfigure}
	\caption{
        Difference in position of shadowing and DC-converged trajectory over time.
	}
	\label{xyzdiff}
\end{figure}

As a first step in the differential correction process, the segments themselves
were used with patch points at each discontinuity.  This failed to converge
using the typical setup, so intermediate patch points were introduced and then
gradually removed as the corrector converged.  This process has been found to
work well in other cases with long trajectories, and is often used in
combination with a reduced correction step size to improve convergence.  For
this particular case, 400 initial patch points were used with an approximate
position update factor of 0.001.  Patch points were then interactively added or
removed in subsequent steps until a converged solution was found with 40 patch points.
The final converged
trajectory from the two-level corrector is 
compared to the original pseudo-orbit in Figure \ref{shortcompare}.
The converged trajectory had total position and velocity
discontinuities of 29 m and 0.006 m/s, respectively.
This final trajectory was compared to the shadowing trajectory obtained with the
multiple-shooting method of Section~\ref{sec:Shadowing_trajectory} in Figures \ref{compareall} and \ref{xyzdiff}.  
Although there is some
difference in the phasing of the trajectory, 
it can be seen that 
both methods
converged to a similar final trajectory.  The largest differences are seen
in small variations in the z component, and near the beginning of the trajectory.
This difference at the beginning might be expected because the differential corrector solution held the
positions constant at the initial and final patch points.

These results demonstrate that the pseudo-orbit process described in this paper can generate accurate initial guesses for higher fidelity studies.

\section*{Acknowledgments}

Part of the research presented in this paper has been carried out at the Jet
Propulsion Laboratory, California Institute of Technology, under a contract
with the National Aeronautics and Space Administration (80NM0018D0004).  
Research of M.G. was partially supported by NSF grant DMS-2154725 and DMS-2307718.
Research of A.D and P.R. was partially supported by Spain's State Research Agency grant PID2024-158570NB-I00.

 \bibliography{mybiblio}{}

@article {FontichM01,
    AUTHOR = {Fontich, Ernest and Mart\'in, Pau},
     TITLE = {Arnold diffusion in perturbations of analytic integrable
              {H}amiltonian systems},
   JOURNAL = {Discrete Contin. Dynam. Systems},
  FJOURNAL = {Discrete and Continuous Dynamical Systems},
    VOLUME = {7},
      YEAR = {2001},
    NUMBER = {1},
     PAGES = {61--84},
      ISSN = {1078-0947,1553-5231},
   MRCLASS = {37J40 (37D10 70H08)},
  MRNUMBER = {1806373},
MRREVIEWER = {Pere\ Guti\'errez},
       DOI = {10.3934/dcds.2001.7.61},
       URL = {https://doi-org.recursos.biblioteca.upc.edu/10.3934/dcds.2001.7.61},
}

@book{KoonLMR08,
  author    = {Koon, W. S. and Lo, M. W. and Marsden, J. E. and Ross, S. D.},
  title     = {Dynamical Systems, the Three-Body Problem and Space Mission Design},
  publisher = {Marsden Books},
  year      = {2008},
  url       = {https://www.researchgate.net/publication/328913173_Dynamical_Systems_the_Three-Body_Problem_and_Space_Mission_Design}
}

@book {GomezLMS01,
    AUTHOR = {G\'omez, G. and Llibre, J. and Mart\'inez, R. and Sim\'o, C.},
     TITLE = {Dynamics and mission design near libration points. {V}ol.
              {I}.},
    SERIES = {World Scientific Monograph Series in Mathematics},
    VOLUME = {2},
      NOTE = {Fundamentals: the case of collinear libration points,
              With a foreword by Walter Flury},
 PUBLISHER = {World Scientific Publishing Co., Inc., River Edge, NJ},
      YEAR = {2001},
     PAGES = {xviii+443},
      ISBN = {981-02-4285-9},
   MRCLASS = {70M20 (37N05 70-08 70F07 70F15)},
  MRNUMBER = {1867240},
MRREVIEWER = {Wang\ Sang\ Koon},
       DOI = {10.1142/9789812810632\_bmatter},
       URL = {https://doi-org.recursos.biblioteca.upc.edu/10.1142/9789812810632_bmatter},
}

@article{LoWBHMCR01,
  author    = {Lo, M. W. and Williams, B. G. and Bollman, W. E. and Han, D. and Muellerschoen, Y. and Cruikshank, D. and Roby, B.},
  title     = {Genesis Mission Design},
  journal   = {The Journal of the Astronautical Sciences},
  volume    = {49},
  number    = {1},
  pages     = {169--184},
  year      = {2001}
}

@article {CanaliasM06,
    AUTHOR = {Canalias, E. and Masdemont, J. J.},
     TITLE = {Homoclinic and heteroclinic transfer trajectories between
              planar {L}yapunov orbits in the {S}un-{E}arth and {E}arth-{M}oon
              systems},
   JOURNAL = {Discrete Contin. Dyn. Syst.},
  FJOURNAL = {Discrete and Continuous Dynamical Systems},
    VOLUME = {14},
      YEAR = {2006},
    NUMBER = {2},
     PAGES = {261--279},
      ISSN = {1078-0947,1553-5231},
   MRCLASS = {70F07 (34C37 37C29 37N05 70F15)},
  MRNUMBER = {2163533},
MRREVIEWER = {Ana\ Nunes},
       DOI = {10.3934/dcds.2006.14.261},
       URL = {https://doi-org.recursos.biblioteca.upc.edu/10.3934/dcds.2006.14.261},
}

@Article{GomezMS98,
author={G{\'o}mez, G.
and Masdemont, J.
and Sim{\'o}, C.},
title={Quasihalo Orbits Associated with Libration Points},
journal={The Journal of the Astronautical Sciences},
year={1998},
month={Jun},
day={01},
volume={46},
number={2},
pages={135-176},
issn={2195-0571},
doi={10.1007/BF03546241},
url={https://doi.org/10.1007/BF03546241}
}

@article {AkingbadeGS23,
    AUTHOR = {Akingbade, Samuel W. and Gidea, Marian and M-Seara, Tere},
     TITLE = {Arnold diffusion in a model of dissipative system},
   JOURNAL = {SIAM J. Appl. Dyn. Syst.},
  FJOURNAL = {SIAM Journal on Applied Dynamical Systems},
    VOLUME = {22},
      YEAR = {2023},
    NUMBER = {3},
     PAGES = {1983--2023},
      ISSN = {1536-0040},
   MRCLASS = {37J40 (37E40 70H08)},
  MRNUMBER = {4623277},
MRREVIEWER = {Jianlu\ Zhang},
       DOI = {10.1137/22M1525508},
       URL = {https://doi.org/10.1137/22M1525508},
}

@article {AlessiBGG25,
    AUTHOR = {Alessi, Elisa Maria and Baldom\'a, Inmaculada and Giralt, Mar
              and Guardia, Marcel},
     TITLE = {On the {A}rnold diffusion mechanism in medium {E}arth orbit},
   JOURNAL = {J. Nonlinear Sci.},
  FJOURNAL = {Journal of Nonlinear Science},
    VOLUME = {35},
      YEAR = {2025},
    NUMBER = {1},
     PAGES = {Paper No. 8, 54},
      ISSN = {0938-8974,1432-1467},
   MRCLASS = {37N05 (37J40 70F07 70F15)},
  MRNUMBER = {4819219},
MRREVIEWER = {Manuele\ Santoprete},
       DOI = {10.1007/s00332-024-10080-0},
       URL = {https://doi.org/10.1007/s00332-024-10080-0},
}

@article {BertiBB03,
    AUTHOR = {Berti, Massimiliano and Biasco, Luca and Bolle, Philippe},
     TITLE = {Drift in phase space: a new variational mechanism with optimal
              diffusion time},
   JOURNAL = {J. Math. Pures Appl. (9)},
  FJOURNAL = {Journal de Math\'ematiques Pures et Appliqu\'ees. Neuvi\`eme
              S\'erie},
    VOLUME = {82},
      YEAR = {2003},
    NUMBER = {6},
     PAGES = {613--664},
      ISSN = {0021-7824},
   MRCLASS = {37J40 (70H09)},
  MRNUMBER = {1996776},
MRREVIEWER = {Mikhail\ B.\ Sevryuk},
       DOI = {10.1016/S0021-7824(03)00032-1},
       URL = {https://doi-org.recursos.biblioteca.upc.edu/10.1016/S0021-7824(03)00032-1},
}

@article {Bounemoura11,
    AUTHOR = {Bounemoura, Abed},
     TITLE = {Generic super-exponential stability of invariant tori in
              {H}amiltonian systems},
   JOURNAL = {Ergodic Theory Dynam. Systems},
  FJOURNAL = {Ergodic Theory and Dynamical Systems},
    VOLUME = {31},
      YEAR = {2011},
    NUMBER = {5},
     PAGES = {1287--1303},
      ISSN = {0143-3857,1469-4417},
   MRCLASS = {37J40 (37J25)},
  MRNUMBER = {2832246},
MRREVIEWER = {Jes\'us\ F.\ Palaci\'an},
       DOI = {10.1017/S0143385710000441},
       URL = {https://doi-org.recursos.biblioteca.upc.edu/10.1017/S0143385710000441},
}

@book {Bowen75,
    AUTHOR = {Bowen, Rufus},
     TITLE = {Equilibrium states and the ergodic theory of {A}nosov
              diffeomorphisms},
    SERIES = {Lecture Notes in Mathematics},
    VOLUME = {Vol. 470},
 PUBLISHER = {Springer-Verlag, Berlin-New York},
      YEAR = {1975},
     PAGES = {i+108},
   MRCLASS = {58F10 (28A65)},
  MRNUMBER = {442989},
MRREVIEWER = {L.\ A.\ Bunimovich},
}

@article {BurgosGS25,
    AUTHOR = {Burgos-Garc\'ia, Jaime and Gidea, Marian and Sierpe, Claudio},
     TITLE = {Arnold diffusion in the elliptic {H}ill four-body problem:
              geometric method and numerical verification},
   JOURNAL = {Nonlinearity},
  FJOURNAL = {Nonlinearity},
    VOLUME = {38},
      YEAR = {2025},
    NUMBER = {9},
     PAGES = {Paper No. 095016, 51},
      ISSN = {0951-7715,1361-6544},
   MRCLASS = {70F10 (37E40 37J40 70H08)},
  MRNUMBER = {4957834},
       DOI = {10.1088/1361-6544/adfff8},
       URL = {https://doi.org/10.1088/1361-6544/adfff8},
}

@article {CabreFL05,
    AUTHOR = {Cabr\'e, Xavier and Fontich, Ernest and de la Llave, Rafael},
     TITLE = {The parameterization method for invariant manifolds. {III}.
              {O}verview and applications},
   JOURNAL = {J. Differential Equations},
  FJOURNAL = {Journal of Differential Equations},
    VOLUME = {218},
      YEAR = {2005},
    NUMBER = {2},
     PAGES = {444--515},
      ISSN = {0022-0396,1090-2732},
   MRCLASS = {37D10 (37M99)},
  MRNUMBER = {2177465},
       DOI = {10.1016/j.jde.2004.12.003},
       URL = {https://doi-org.recursos.biblioteca.upc.edu/10.1016/j.jde.2004.12.003},
}

@article {CaraccioloL20,
    AUTHOR = {Caracciolo, Chiara and Locatelli, Ugo},
     TITLE = {Computer-assisted estimates for {B}irkhoff normal forms},
   JOURNAL = {J. Comput. Dyn.},
  FJOURNAL = {Journal of Computational Dynamics},
    VOLUME = {7},
      YEAR = {2020},
    NUMBER = {2},
     PAGES = {425--460},
      ISSN = {2158-2491,2158-2505},
   MRCLASS = {37J40 (68W40 70F07 70H08)},
  MRNUMBER = {4128817},
       DOI = {10.3934/jcd.2020017},
       URL = {https://doi-org.recursos.biblioteca.upc.edu/10.3934/jcd.2020017},
}

@article {DelshamsGR13,
    AUTHOR = {Delshams, Amadeu and Gidea, Marian and Rold\'an, Pablo},
     TITLE = {Transition map and shadowing lemma for normally hyperbolic
              invariant manifolds},
   JOURNAL = {Discrete Contin. Dyn. Syst.},
  FJOURNAL = {Discrete and Continuous Dynamical Systems},
    VOLUME = {33},
      YEAR = {2013},
    NUMBER = {3},
     PAGES = {1089--1112},
      ISSN = {1078-0947,1553-5231},
   MRCLASS = {37D10 (34C37 34C45 37N05 70F07 70H05)},
  MRNUMBER = {2988936},
MRREVIEWER = {Jaap\ Eldering},
       DOI = {10.3934/dcds.2013.33.1089},
       URL = {https://doi-org.recursos.biblioteca.upc.edu/10.3934/dcds.2013.33.1089},
}

@article {DelshamsGR16,
    AUTHOR = {Delshams, Amadeu and Gidea, Marian and Roldan, Pablo},
     TITLE = {Arnold's mechanism of diffusion in the spatial circular
              restricted three-body problem: a semi-analytical argument},
   JOURNAL = {Phys. D},
  FJOURNAL = {Physica D. Nonlinear Phenomena},
    VOLUME = {334},
      YEAR = {2016},
     PAGES = {29--48},
      ISSN = {0167-2789,1872-8022},
   MRCLASS = {70F07 (37N05 70M20)},
  MRNUMBER = {3545968},
MRREVIEWER = {Lei\ Zhao},
       DOI = {10.1016/j.physd.2016.06.005},
       URL = {https://doi-org.recursos.biblioteca.upc.edu/10.1016/j.physd.2016.06.005},
}

@article {DelshamsGR26,
    AUTHOR = {Delshams, Amadeu and Gidea, Marian and Roldan, Pablo},
     TITLE = {Semi-analytic construction of global transfers between
              quasi-periodic orbits in the spatial {R}3{BP}},
   JOURNAL = {Commun. Nonlinear Sci. Numer. Simul.},
  FJOURNAL = {Communications in Nonlinear Science and Numerical Simulation},
    VOLUME = {152},
      YEAR = {2026},
     PAGES = {Paper No. 109220, 34},
      ISSN = {1007-5704,1878-7274},
   MRCLASS = {70F07 (34C25)},
  MRNUMBER = {4950183},
       DOI = {10.1016/j.cnsns.2025.109220},
       URL = {https://doi-org.recursos.biblioteca.upc.edu/10.1016/j.cnsns.2025.109220},
}

@article {DelshamsG96,
    AUTHOR = {Delshams, Amadeu and Guti\'errez, Pere},
     TITLE = {Estimates on invariant tori near an elliptic equilibrium point
              of a {H}amiltonian system},
   JOURNAL = {J. Differential Equations},
  FJOURNAL = {Journal of Differential Equations},
    VOLUME = {131},
      YEAR = {1996},
    NUMBER = {2},
     PAGES = {277--303},
      ISSN = {0022-0396,1090-2732},
   MRCLASS = {58F27 (34C30 58F30 70H05)},
  MRNUMBER = {1419016},
MRREVIEWER = {Luigi\ Chierchia},
       DOI = {10.1006/jdeq.1996.0165},
       URL = {https://doi-org.recursos.biblioteca.upc.edu/10.1006/jdeq.1996.0165},
}

@article {DelshamsGS24,
    AUTHOR = {Delshams, A. and Granados, A. and Schaefer, R. G.},
     TITLE = {Arnold diffusion for an a priori unstable {H}amiltonian system
              with {$3 + 1/2$} degrees of freedom},
   JOURNAL = {Chaos},
  FJOURNAL = {Chaos. An Interdisciplinary Journal of Nonlinear Science},
    VOLUME = {34},
      YEAR = {2024},
    NUMBER = {6},
     PAGES = {Paper No. 063118, 17},
      ISSN = {1054-1500,1089-7682},
   MRCLASS = {37J40 (70H05 70H06 70H14)},
  MRNUMBER = {4755372},
MRREVIEWER = {Xifeng\ Su},
       DOI = {10.1063/5.0185044},
       URL = {https://doi.org/10.1063/5.0185044},
}

@article {DelshamsLS16,
    AUTHOR = {Delshams, Amadeu and de la Llave, Rafael and Seara, Tere M.},
     TITLE = {Instability of high dimensional {H}amiltonian systems:
              multiple resonances do not impede diffusion},
   JOURNAL = {Adv. Math.},
  FJOURNAL = {Advances in Mathematics},
    VOLUME = {294},
      YEAR = {2016},
     PAGES = {689--755},
      ISSN = {0001-8708,1090-2082},
   MRCLASS = {37J25},
  MRNUMBER = {3479576},
MRREVIEWER = {Camelia\ Arie\c sanu},
       DOI = {10.1016/j.aim.2015.11.010},
       URL = {https://doi.org/10.1016/j.aim.2015.11.010},
}

@article {DelshamsZ25,
    AUTHOR = {Delshams, Amadeu and Zhang, Ke},
     TITLE = {Generic global diffusion for analytic a priori unstable
              systems},
   JOURNAL = {Comm. Math. Phys.},
  FJOURNAL = {Communications in Mathematical Physics},
    VOLUME = {406},
      YEAR = {2025},
    NUMBER = {6},
     PAGES = {Paper No. 146, 22},
      ISSN = {0010-3616,1432-0916},
   MRCLASS = {37J40 (70H05 70H08)},
  MRNUMBER = {4913967},
MRREVIEWER = {Qinbo\ Chen},
       DOI = {10.1007/s00220-025-05342-1},
       URL = {https://doi.org/10.1007/s00220-025-05342-1},
}

@article {EfthyGC04,
    AUTHOR = {Efthymiopoulos, C. and Giorgilli, A. and Contopoulos, G.},
     TITLE = {Nonconvergence of formal integrals. {II}. {I}mproved estimates
              for the optimal order of truncation},
   JOURNAL = {J. Phys. A},
  FJOURNAL = {Journal of Physics. A. Mathematical and General},
    VOLUME = {37},
      YEAR = {2004},
    NUMBER = {45},
     PAGES = {10831--10858},
      ISSN = {0305-4470,1751-8121},
   MRCLASS = {70H05 (70-08 70H08 70K55)},
  MRNUMBER = {2100146},
       DOI = {10.1088/0305-4470/37/45/008},
       URL = {https://doi-org.recursos.biblioteca.upc.edu/10.1088/0305-4470/37/45/008},
}

@article {GiorgilliDFGS89,
    AUTHOR = {Giorgilli, Antonio and Delshams, Amadeo and Fontich, Ernest
              and Galgani, Luigi and Sim\'o, Carles},
     TITLE = {Effective stability for a {H}amiltonian system near an
              elliptic equilibrium point, with an application to the
              restricted three-body problem},
   JOURNAL = {J. Differential Equations},
  FJOURNAL = {Journal of Differential Equations},
    VOLUME = {77},
      YEAR = {1989},
    NUMBER = {1},
     PAGES = {167--198},
      ISSN = {0022-0396,1090-2732},
   MRCLASS = {70K20 (12H20 58F05 58F30 70F07 70H05)},
  MRNUMBER = {980547},
MRREVIEWER = {James\ A.\ Murdock},
       DOI = {10.1016/0022-0396(89)90161-7},
       URL = {https://doi-org.recursos.biblioteca.upc.edu/10.1016/0022-0396(89)90161-7},
}

@article{GuzzoCP26,
author = {Massimiliano Guzzo and Chiara Caracciolo and Gabriella Pinzari},
title = {Improved stability estimates at elliptic equilibria of Hamiltonian systems},
JOURNAL = {Discrete Contin. Dyn. Syst.},
FJOURNAL = {Discrete and Continuous Dynamical Systems. Series A},
VOLUME = {77},
year = {2026},
pages = {1--31},
issn = {1078-0947},
doi = {10.3934/dcds.2026063},
url = {https://www.aimsciences.org/article/id/69bbb40cd1a02b472f36a46e}
}

@book {HaroCFLM16,
    AUTHOR = {Haro, \`Alex and Canadell, Marta and Figueras, Jordi-Llu\'is
              and Luque, Alejandro and Mondelo, Josep-Maria},
     TITLE = {The parameterization method for invariant manifolds},
    SERIES = {Applied Mathematical Sciences},
    VOLUME = {195},
      NOTE = {From rigorous results to effective computations},
 PUBLISHER = {Springer, [Cham]},
      YEAR = {2016},
     PAGES = {xvi+267},
      ISBN = {978-3-319-29660-9; 978-3-319-29662-3},
   MRCLASS = {37-02 (34C45 37C55 37D10 37J40 37M99 65P40)},
  MRNUMBER = {3467671},
MRREVIEWER = {Martin\ Hermann},
       DOI = {10.1007/978-3-319-29662-3},
       URL = {https://doi-org.recursos.biblioteca.upc.edu/10.1007/978-3-319-29662-3},
}

@article {Jorba99,
	AUTHOR = {Jorba, \`{A}ngel},
	TITLE = {A methodology for the numerical computation of normal forms,
	centre manifolds and first integrals of {H}amiltonian systems},
	JOURNAL = {Experiment. Math.},
	FJOURNAL = {Experimental Mathematics},
	VOLUME = {8},
	YEAR = {1999},
	NUMBER = {2},
	PAGES = {155--195},
	ISSN = {1058-6458,1944-950X},
	MRCLASS = {37M99 (34-04 34C20 37-04 37J40 70F15 70H05 70K50)},
	MRNUMBER = {1700577},
	URL = {http://projecteuclid.org/euclid.em/1047477059},
}

@article {JorbaM99,
	AUTHOR = {Jorba, \`{A}ngel and Masdemont, Josep},
	TITLE = {Dynamics in the center manifold of the collinear points of the
	restricted three body problem},
	JOURNAL = {Phys. D},
	FJOURNAL = {Physica D. Nonlinear Phenomena},
	VOLUME = {132},
	YEAR = {1999},
	NUMBER = {1-2},
	PAGES = {189--213},
	ISSN = {0167-2789,1872-8022},
	MRCLASS = {70F07 (37J99 37N05 70H05)},
	MRNUMBER = {1705705},
	DOI = {10.1016/S0167-2789(99)00042-1},
	URL = {https://doi.org/10.1016/S0167-2789(99)00042-1},
}

@article {DLS08,
	AUTHOR = {Delshams, Amadeu and de la Llave, Rafael and Seara, Tere M.},
	TITLE = {Geometric properties of the scattering map of a normally
	hyperbolic invariant manifold},
	JOURNAL = {Adv. Math.},
	FJOURNAL = {Advances in Mathematics},
	VOLUME = {217},
	YEAR = {2008},
	NUMBER = {3},
	PAGES = {1096--1153},
	ISSN = {0001-8708,1090-2082},
	MRCLASS = {37J40 (37D10 37J99 37N20 81U05)},
	MRNUMBER = {2383896},
	MRREVIEWER = {Enrico\ Valdinoci},
	DOI = {10.1016/j.aim.2007.08.014},
	URL = {https://doi.org/10.1016/j.aim.2007.08.014},
}

@article {DGR16,
	AUTHOR = {Delshams, Amadeu and Gidea, Marian and Roldan, Pablo},
	TITLE = {Arnold's mechanism of diffusion in the spatial circular
	restricted three-body problem: a semi-analytical argument},
	JOURNAL = {Phys. D},
	FJOURNAL = {Physica D. Nonlinear Phenomena},
	VOLUME = {334},
	YEAR = {2016},
	PAGES = {29--48},
	ISSN = {0167-2789,1872-8022},
	MRCLASS = {70F07 (37N05 70M20)},
	MRNUMBER = {3545968},
	MRREVIEWER = {Lei\ Zhao},
	DOI = {10.1016/j.physd.2016.06.005},
	URL = {https://doi.org/10.1016/j.physd.2016.06.005},
}

@article {DelshamsS17,
    AUTHOR = {Delshams, Amadeu and Schaefer, Rodrigo G.},
     TITLE = {Arnold diffusion for a complete family of perturbations},
   JOURNAL = {Regul. Chaotic Dyn.},
  FJOURNAL = {Regular and Chaotic Dynamics. International Scientific
              Journal},
    VOLUME = {22},
      YEAR = {2017},
    NUMBER = {1},
     PAGES = {78--108},
      ISSN = {1560-3547,1468-4845},
   MRCLASS = {37J40 (70H33)},
  MRNUMBER = {3607642},
MRREVIEWER = {Qinbo\ Chen},
       DOI = {10.1134/S1560354717010051},
       URL = {https://doi.org/10.1134/S1560354717010051},
}

@article {DelshamsS18,
    AUTHOR = {Delshams, Amadeu and Schaefer, Rodrigo G.},
     TITLE = {Arnold diffusion for a complete family of perturbations with
              two independent harmonics},
   JOURNAL = {Discrete Contin. Dyn. Syst.},
  FJOURNAL = {Discrete and Continuous Dynamical Systems},
    VOLUME = {38},
      YEAR = {2018},
    NUMBER = {12},
     PAGES = {6047--6072},
      ISSN = {1078-0947,1553-5231},
   MRCLASS = {37J40 (70H09)},
  MRNUMBER = {3917800},
MRREVIEWER = {Patricia\ Yanguas},
       DOI = {10.3934/dcds.2018261},
       URL = {https://doi.org/10.3934/dcds.2018261},
}

@article {CapinskiG23,
    AUTHOR = {Capi\'nski, Maciej J. and Gidea, Marian},
     TITLE = {Arnold diffusion, quantitative estimates, and stochastic
              behavior in the three-body problem},
   JOURNAL = {Comm. Pure Appl. Math.},
  FJOURNAL = {Communications on Pure and Applied Mathematics},
    VOLUME = {76},
      YEAR = {2023},
    NUMBER = {3},
     PAGES = {616--681},
      ISSN = {0010-3640,1097-0312},
   MRCLASS = {70F07 (37N05)},
  MRNUMBER = {4544807},
MRREVIEWER = {Patricia\ Yanguas},
       DOI = {10.1002/cpa.22014},
       URL = {https://doi.org/10.1002/cpa.22014},
}

@article {DelshamsLS06a,
    AUTHOR = {Delshams, Amadeu and de la Llave, Rafael and Seara, Tere M.},
     TITLE = {A geometric mechanism for diffusion in {H}amiltonian systems
              overcoming the large gap problem: heuristics and rigorous
              verification on a model},
   JOURNAL = {Mem. Amer. Math. Soc.},
  FJOURNAL = {Memoirs of the American Mathematical Society},
    VOLUME = {179},
      YEAR = {2006},
    NUMBER = {844},
     PAGES = {viii+141},
      ISSN = {0065-9266,1947-6221},
   MRCLASS = {37J40 (70H09)},
  MRNUMBER = {2184276},
MRREVIEWER = {Jean-Pierre\ Marco},
       DOI = {10.1090/memo/0844},
       URL = {https://doi.org/10.1090/memo/0844},
}

@article {DelshamsLS06b,
    AUTHOR = {Delshams, Amadeu and de la Llave, Rafael and Seara, Tere M.},
     TITLE = {Orbits of unbounded energy in quasi-periodic perturbations of
              geodesic flows},
   JOURNAL = {Adv. Math.},
  FJOURNAL = {Advances in Mathematics},
    VOLUME = {202},
      YEAR = {2006},
    NUMBER = {1},
     PAGES = {64--188},
      ISSN = {0001-8708,1090-2082},
   MRCLASS = {37J45 (37D10 37D40 37J40 53D25 70H08 70H09)},
  MRNUMBER = {2218821},
MRREVIEWER = {Mikhail\ B.\ Sevryuk},
       DOI = {10.1016/j.aim.2005.03.005},
       URL = {https://doi.org/10.1016/j.aim.2005.03.005},
}

@article {GideaLS20,
    AUTHOR = {Gidea, Marian and de la Llave, Rafael and Seara, Tere M.},
     TITLE = {A general mechanism of instability in {H}amiltonian systems:
              skipping along a normally hyperbolic invariant manifold},
   JOURNAL = {Discrete Contin. Dyn. Syst.},
  FJOURNAL = {Discrete and Continuous Dynamical Systems},
    VOLUME = {40},
      YEAR = {2020},
    NUMBER = {12},
     PAGES = {6795--6813},
      ISSN = {1078-0947,1553-5231},
   MRCLASS = {37J25 (70H08 70H33)},
  MRNUMBER = {4160091},
MRREVIEWER = {Lennard\ F.\ Bakker},
       DOI = {10.3934/dcds.2020166},
       URL = {https://doi.org/10.3934/dcds.2020166},
}

@article {GuardiaPS25,
    AUTHOR = {Guardia, Marcel and Paradela, Jaime and Seara, Tere M.},
     TITLE = {A degenerate {A}rnold diffusion mechanism in the restricted
              3-body problem},
   JOURNAL = {Ann. Sci. \'Ec. Norm. Sup\'er. (4)},
  FJOURNAL = {Annales Scientifiques de l'\'Ecole Normale Sup\'erieure.
              Quatri\`eme S\'erie},
    VOLUME = {58},
      YEAR = {2025},
    NUMBER = {6},
     PAGES = {1401--1478},
      ISSN = {0012-9593,1873-2151},
   MRCLASS = {37J40 (30B40)},
  MRNUMBER = {5013757},
       DOI = {10.24033/asens.2632},
       URL = {https://doi.org/10.24033/asens.2632},
}

@ARTICLE{Arnold64,
  author = {Arnold, V.I.},
  title = {Instability of dynamical systems with several degrees of freedom},
  year = 1964,
  journal = {Sov. Math. Doklady},
  volume = 5,
  pages = {581--585},
  note = {See also \emph{Vladimir {I}. {A}rnold---collected works. {V}ol. {I}}, Springer-Verlag, Berlin, 2009, pp.~423--427,
  \href{https://doi.org/10.1007/978-3-642-01742-1}{DOI 10.1007/978-3-642-01742-1}, or
  \emph{Hamiltonian Dynamical Systems}, edited by R.S MacKay, J.D Meiss, CRC Press, 2020,
   pp.~633--637, \href{https://doi.org/10.1201/9781003069515}{DOI 10.1201/9781003069515}}
 }

@article {DelshamsKRS19,
    AUTHOR = {Delshams, Amadeu and Kaloshin, Vadim and de la Rosa, Abraham
              and Seara, Tere M.},
     TITLE = {Global instability in the restricted planar elliptic three
              body problem},
   JOURNAL = {Comm. Math. Phys.},
  FJOURNAL = {Communications in Mathematical Physics},
    VOLUME = {366},
      YEAR = {2019},
    NUMBER = {3},
     PAGES = {1173--1228},
      ISSN = {0010-3616,1432-0916},
   MRCLASS = {70F07 (37N05 70F15 70K50)},
  MRNUMBER = {3927089},
MRREVIEWER = {Eduardo\ S. G. Leandro},
       DOI = {10.1007/s00220-018-3248-z},
       URL = {https://doi.org/10.1007/s00220-018-3248-z},
}

@article {ClarkeFG24,
    AUTHOR = {Clarke, Andrew and Fejoz, Jacques and Guardia, Marcel},
     TITLE = {A counterexample to the theorem of {L}aplace-{L}agrange on the
              stability of semimajor axes},
   JOURNAL = {Arch. Ration. Mech. Anal.},
  FJOURNAL = {Archive for Rational Mechanics and Analysis},
    VOLUME = {248},
      YEAR = {2024},
    NUMBER = {2},
     PAGES = {Paper No. 19, 73},
      ISSN = {0003-9527,1432-0673},
   MRCLASS = {70F10},
  MRNUMBER = {4707559},
MRREVIEWER = {John\ A.\ Arredondo},
       DOI = {10.1007/s00205-024-01960-6},
       URL = {https://doi.org/10.1007/s00205-024-01960-6},
}

@article {ClarkeFG25,
    AUTHOR = {Clarke, Andrew and Fejoz, Jacques and Guardia, Marcel},
     TITLE = {Why are inner planets not inclined?},
   JOURNAL = {Publ. Math. Inst. Hautes \'Etudes Sci.},
  FJOURNAL = {Publications Math\'ematiques. Institut de Hautes \'Etudes
              Scientifiques},
    VOLUME = {141},
      YEAR = {2025},
     PAGES = {1--98},
      ISSN = {0073-8301,1618-1913},
   MRCLASS = {70F15 (37N05 85A05)},
  MRNUMBER = {4916233},
MRREVIEWER = {Eduardo\ S. G. Leandro},
       DOI = {10.1007/s10240-024-00151-z},
       URL = {https://doi.org/10.1007/s10240-024-00151-z},
}

@article {ClarkeT23,
    AUTHOR = {Clarke, Andrew and Turaev, Dmitry},
     TITLE = {Arnold diffusion in multidimensional convex billiards},
   JOURNAL = {Duke Math. J.},
  FJOURNAL = {Duke Mathematical Journal},
    VOLUME = {172},
      YEAR = {2023},
    NUMBER = {10},
     PAGES = {1813--1878},
      ISSN = {0012-7094,1547-7398},
   MRCLASS = {37C83 (37J40)},
  MRNUMBER = {4624370},
MRREVIEWER = {Jianlu\ Zhang},
       DOI = {10.1215/00127094-2022-0073},
       URL = {https://doi.org/10.1215/00127094-2022-0073},
}

@article {Clarke24,
    AUTHOR = {Clarke, Andrew},
     TITLE = {Geodesics with unbounded speed on fluctuating surfaces},
   JOURNAL = {Regul. Chaotic Dyn.},
  FJOURNAL = {Regular and Chaotic Dynamics. International Scientific
              Journal},
    VOLUME = {29},
      YEAR = {2024},
    NUMBER = {3},
     PAGES = {435--450},
      ISSN = {1560-3547,1468-4845},
   MRCLASS = {37D40 (37J40)},
  MRNUMBER = {4752409},
       DOI = {10.1134/S1560354724030018},
       URL = {https://doi.org/10.1134/S1560354724030018},
}

@article {DelshamsLS00,
    AUTHOR = {Delshams, Amadeu and de la Llave, Rafael and Seara, Tere M.},
     TITLE = {A geometric approach to the existence of orbits with unbounded
              energy in generic periodic perturbations by a potential of
              generic geodesic flows of {${\bf T}^2$}},
   JOURNAL = {Comm. Math. Phys.},
  FJOURNAL = {Communications in Mathematical Physics},
    VOLUME = {209},
      YEAR = {2000},
    NUMBER = {2},
     PAGES = {353--392},
      ISSN = {0010-3616,1432-0916},
   MRCLASS = {37J40 (53D25)},
  MRNUMBER = {1737988},
MRREVIEWER = {Mikhail\ B.\ Sevryuk},
       DOI = {10.1007/PL00020961},
       URL = {https://doi.org/10.1007/PL00020961},
}

@article {Treschev04,
    AUTHOR = {Treschev, D.},
     TITLE = {Evolution of slow variables in a priori unstable {H}amiltonian systems},
   JOURNAL = {Nonlinearity},
  FJOURNAL = {Nonlinearity},
    VOLUME = {17},
      YEAR = {2004},
    NUMBER = {5},
     PAGES = {1803--1841},
      ISSN = {0951-7715},
     CODEN = {NONLE5},
}

@article {DelshamsH09,
    AUTHOR = {Delshams, Amadeu and Huguet, Gemma},
     TITLE = {Geography of resonances and {A}rnold diffusion in a priori
              unstable {H}amiltonian systems},
   JOURNAL = {Nonlinearity},
  FJOURNAL = {Nonlinearity},
    VOLUME = {22},
      YEAR = {2009},
    NUMBER = {8},
     PAGES = {1997--2077},
      ISSN = {0951-7715,1361-6544},
   MRCLASS = {37J40},
  MRNUMBER = {2525822},
MRREVIEWER = {Paolo\ Perfetti},
       DOI = {10.1088/0951-7715/22/8/013},
       URL = {https://doi.org/10.1088/0951-7715/22/8/013},
}

@article {DelshamsH11,
    AUTHOR = {Delshams, Amadeu and Huguet, Gemma},
     TITLE = {A geometric mechanism of diffusion: rigorous verification in a
              priori unstable {H}amiltonian systems},
   JOURNAL = {J. Differential Equations},
  FJOURNAL = {Journal of Differential Equations},
    VOLUME = {250},
      YEAR = {2011},
    NUMBER = {5},
     PAGES = {2601--2623},
      ISSN = {0022-0396,1090-2732},
   MRCLASS = {37J40 (37J45 70H05)},
  MRNUMBER = {2756078},
MRREVIEWER = {Massimiliano\ Berti},
       DOI = {10.1016/j.jde.2010.12.023},
       URL = {https://doi.org/10.1016/j.jde.2010.12.023},
}

@article {ChenL22,
    AUTHOR = {Chen, Qinbo and de la Llave, Rafael},
     TITLE = {Analytic genericity of diffusing orbits in {\it a priori}
              unstable {H}amiltonian systems},
   JOURNAL = {Nonlinearity},
  FJOURNAL = {Nonlinearity},
    VOLUME = {35},
      YEAR = {2022},
    NUMBER = {4},
     PAGES = {1986--2019},
      ISSN = {0951-7715,1361-6544},
   MRCLASS = {37J40},
  MRNUMBER = {4399498},
MRREVIEWER = {Lufang\ Mi},
       DOI = {10.1088/1361-6544/ac50bb},
       URL = {https://doi.org/10.1088/1361-6544/ac50bb},
}

@article{GelfreichT17,
	author={Gelfreich, Vassili and Turaev, Dmitry},
	title={Arnold Diffusion in A Priori Chaotic Symplectic Maps},
	journal={Communications in Mathematical Physics},
	year={2017},
	month={07},
	day={01},
	volume={353},
	number={2},
	pages={507--547},
	issn={1432-0916},
	doi={10.1007/s00220-017-2867-0},
	url={https://doi.org/10.1007/s00220-017-2867-0},
}

@article {GideaLlaveSeara20-CPAM,
    AUTHOR = {Gidea, Marian and de la Llave, Rafael and M-Seara, Tere},
     TITLE = {A general mechanism of diffusion in {H}amiltonian systems:
              qualitative results},
   JOURNAL = {Comm. Pure Appl. Math.},
  FJOURNAL = {Communications on Pure and Applied Mathematics},
    VOLUME = {73},
      YEAR = {2020},
    NUMBER = {1},
     PAGES = {150--209},
      ISSN = {0010-3640,1097-0312},
   MRCLASS = {37J06 (70H08)},
  MRNUMBER = {4033892},
MRREVIEWER = {Arturo\ Vieiro},
       DOI = {10.1002/cpa.21856},
       URL = {https://doi.org/10.1002/cpa.21856},
}

@book{Szebehely67,
	AUTHOR = {Szebehely, V.},
	TITLE = {Theory of Orbits},
	SERIES = {},
	EDITION = {},
	NOTE = {},
	PUBLISHER = {Academic Press},
	ADDRESS = {},
	YEAR = {1967},
	PAGES = {},
	ISBN = {},
	MRCLASS = {},
	MRNUMBER = {},
}

@article {GideaL18,
    AUTHOR = {Gidea, Marian and de la Llave, Rafael},
     TITLE = {Global {M}elnikov theory in {H}amiltonian systems with general
              time-dependent perturbations},
   JOURNAL = {J. Nonlinear Sci.},
  FJOURNAL = {Journal of Nonlinear Science},
    VOLUME = {28},
      YEAR = {2018},
    NUMBER = {5},
     PAGES = {1657--1707},
      ISSN = {0938-8974,1432-1467},
   MRCLASS = {37J40 (34C37 34D35 37C29 37C50 37C60 70H08)},
  MRNUMBER = {3846872},
MRREVIEWER = {Mikhail\ B.\ Sevryuk},
       DOI = {10.1007/s00332-018-9461-2},
       URL = {https://doi.org/10.1007/s00332-018-9461-2},
}

@article {KPC07,
    AUTHOR = {Ko\c cak, H. and Palmer, K. and Coomes, B.},
     TITLE = {Shadowing in ordinary differential equations},
   JOURNAL = {Rend. Semin. Mat. Univ. Politec. Torino},
  FJOURNAL = {Rendiconti del Seminario Matematico. Universit\`a{} e
              Politecnico Torino},
    VOLUME = {65},
      YEAR = {2007},
    NUMBER = {1},
     PAGES = {89--113},
      ISSN = {0373-1243,2704-999X},
   MRCLASS = {37C50 (34C28 37C10)},
  MRNUMBER = {2339601},
MRREVIEWER = {Kazuhiro\ Sakai},
}

@article{zgliczynski2004covering,
  title={Covering relations for multidimensional dynamical systems},
  author={Zgliczy{\'n}ski, Piotr and Gidea, Marian},
  journal={Journal of Differential Equations},
  volume={202},
  number={1},
  pages={32--58},
  year={2004},
  publisher={Elsevier}
}

@article{Reach:1997,
  author =        {W. T. Reach and B. A. Franz and J. L. Weiland},
  journal =       {Icarus},
  number =        {IS975704},
  pages =         {461--484},
  title =         {The Three-Dimensional Structure of the Zodiacal Dust
                   Bands},
  volume =        {127},
  year =          {1997},
}

@phdthesis{May:2007,
  address =       {London},
  author =        {B. H. May},
  school =        {Imperial College},
  title =         {A Survey of Radial Velocities in the Zodiacal Dust
                   Cloud},
  year =          {2007},
}

@article{Wenzel:1995,
  author =        {K. Wenzel},
  journal =       {Nuclear Physics B - Proceedings SUpplements},
  pages =         {59--68},
  title =         {The Ulysses Mission - A Voyage to the Poles of the
                   Sun},
  volume =        {39},
  year =          {1995},
}

@inproceedings{Stuart:2019,
  address =       {Portland, ME},
  author =        {Jeffrey Stuart and Rodney L. Anderson and
                   Christopher Sullivan and Natasha Bosanac},
  booktitle =     {AAS/AIAA Astrodynamics Specialist Conference},
  month =         {August 11--15},
  number =        {AAS 19-728},
  title =         {Accessing Highly Out-of-Ecliptic Science Orbits via
                   Low-Energy, Low-Thrust Transport Mechanisms},
  year =          {2019},
}

@inproceedings{Koon:2000b,
  address =       {Clearwater, Florida},
  author =        {Koon, W.~S. and Martin W. Lo and J. E. Marsden and
                   S. D. Ross},
  booktitle =     {AAS/AIAA Astrodynamics Specialist Conference},
  month =         {January 23-26},
  number =        {AAS 00-166},
  title =         {Shoot the Moon},
  year =          {2000},
}

@article{Lo:2002b,
  author =        {Martin W. Lo},
  journal =       {Spaceflight Now},
  month =         {July},
  title =         {Making Space Travel Simpler With The Interplanetary
                   Highway},
  year =          {2002},
}

@article{Lo:2001c,
  author =        {Martin W. Lo and Bobby G. Williams and
                   Williard E. Bollman and Dongsuk Han and Yungsun Hahn and
                   Julia L. Bell and Edward A. Hirst and
                   Robert A. Corwin and Philip E. Hong and
                   Kathleen C. Howell and Brian Barden and Roby Wilson},
  journal =       {The Journal of the Astronautical Sciences},
  month =         {January - March},
  number =        {1},
  pages =         {169--184},
  title =         {Genesis Mission Design},
  volume =        {49},
  year =          {2001},
}

@techreport{Wilson:2003,
  author =        {Wilson, R. S.},
  institution =   {Jet Propulsion Laboratory},
  number =        {312.I-03-002},
  type =          {IOM},
  title =         {Derivation of Differential Correctors Used in GENESIS
                   Mission Design},
  year =          {2003},
}

@mastersthesis{Pernicka:1986,
  author =        {Pernicka, H. J.},
  month =         {December},
  school =        {Purdue University},
  title =         {The Numerical Determination of Lissajous Orbits in
                   the Circular Restricted Three-Body Problem},
  year =          {1986},
}

@book{Parker:2014,
  address =       {Hoboken, New Jersey},
  author =        {Jeffrey S. Parker and Rodney L. Anderson},
  edition =       {1\textsuperscript{st}},
  month =         {June},
  publisher =     {John Wiley \& Sons, Inc.},
  series =        {JPL Deep Space Communications and Navigation Series},
  title =         {Low-Energy Lunar Trajectory Design},
  volume =        {12},
  year =          {2014},
}

@article{Parker:2015,
author = {Jeffrey S. Parker and Timothy P. McElrath and Rodney L. Anderson and Theodore H. Sweetser},
journal = {Journal of the Astronautical Sciences},
month = {March},
number = {1},
pages = {44--72},
title = {Trajectory Design for MoonRise: A Proposed Lunar South Pole-Aitken Basin Sample Return Mission},
volume = {62},
year = {2015}}

@article{Anderson:2018,
author = {Rodney L. Anderson and Stefano Campagnola and Brent B. Buffington},
journal = {Journal of Guidance, Control, and Dynamics},
month = {April},
number = {4},
pages = {827--840},
title = {Analysis of Petal Rotation Trajectory Characteristics},
volume = {41},
year = {2018}}
 \bibliographystyle{alphaabbrv}
\end{document}